\documentclass[11pt]{article} 
\usepackage{amsmath,amssymb,amsthm}
\usepackage{a4wide}
 
\newtheorem{theorem}{Theorem}

\newtheorem{corollary}[theorem]{Corollary}
\newtheorem{proposition}{Proposition}
\newtheorem{lemma}{Lemma} 
\newtheorem{claim}{Claim}
\newtheorem{remark}{Remark}

\numberwithin{claim}{section}
\numberwithin{equation}{section}
\numberwithin{lemma}{section}
\numberwithin{proposition}{section}

\newcommand{\inch}{0}
\newcommand{\RR}{\mathbb{R}}

\newcommand{\JJ}{\mathcal{J}}

\newcommand{\px}{\partial_x}

\newcommand{\yy}{y}
\newcommand{\OTT}{\mathit{O}}
\newcommand{\OO}{\mathcal{O}}

\newcommand{\FF}{F}
\newcommand{\FTF}{\widetilde F}
\newcommand{\qun}{R_1}
\newcommand{\qde}{R_2}
\newcommand{\qtud}{\widetilde R_j}
\newcommand{\qtun}{\widetilde R_1}
\newcommand{\qtde}{\widetilde R_2}
\newcommand{\qoud}{\overline R_j}
\newcommand{\qoun}{\overline R_1}
\newcommand{\qode}{\overline R_2}

\newcommand{\qud}{R_j}

\newcommand{\sud}{\sum_{j=1,2}}
\newcommand{\summu}{{\bar \mu}}
\newcommand{\YYzz}{Y_0}
\newcommand{\YYrr}{Y}
\newcommand{\pyy}{y}
\newcommand{\pzz}{z}
\newcommand{\pyyb}{{\bar y}}
\newcommand{\fell}{\widetilde E}
\newcommand{\iscal}{J}
\newcommand{\TSR}{T}

\begin{document}

\title{Inelastic interaction of nearly equal solitons for the\\ BBM equation\footnote{
This research was supported in part by the Agence Nationale de la Recherche
(ANR ONDENONLIN).}}
\author{Yvan Martel$^{(1)}$ \and Frank Merle$^{(2)}$}
\date{ \small (1) Universit\'e de Versailles Saint-Quentin-en-Yvelines and IUF
\\
(2) Universit\'e de Cergy-Pontoise and IHES
}
\maketitle

\begin{abstract}
This paper is concerned with the interaction of two solitons of nearly equal speeds for the \eqref{eq:BBMi} equation.
This work is an extension of \cite{MMkdv4} addressing the same question for the quartic (gKdV) equation.
We consider the \eqref{eq:BBMi} equation, for   $\lambda \in [0,1)$,
\begin{equation} \label{eq:BBMi} \tag{BBM}
	(1-\lambda \partial_x^2) \partial_t u + \partial_x (\partial_x^2 u - u + u^2) =0.\end{equation}
Solitons are solutions of the form
$
	R_{\mu,x_0}(t,x)=Q_{\mu}(x-\mu t -x_0),
$
for $\mu>-1$, $x_0\in \RR$.

For    $\mu_0>0$ small, let $U(t,x)$ be  the unique solution   of \eqref{eq:BBMi} such that
$$
\lim_{t\to -\infty} \|{U}(t) -  Q_{-\mu_0}(.+\mu_0 t) - Q_{\mu_0}(.-\mu_0 t )\|_{H^1} = 0.
$$
First, we prove that ${U}(t)$ remains close to the sum of two solitons, for all time $t\in \RR$, 
$${U}(t,x) = {Q}_{\mu_1(t)}(x-y_1(t)) + {Q}_{\mu_2(t)}(x-y_2(t)) + {\varepsilon}(t) \quad \text{where} \quad
 \|\varepsilon(t)\|\leq \mu_0^{2^-},$$
with
$
y_1(t)-y_2(t)> 2 |\ln  \mu_0| + O(1),
$ which means that  at the main order the situation is similar to the integrable KdV case.
However, we show that the collision is perfectly elastic if and only if $\lambda=0$
(i.e. only in the integrable   case).
\end{abstract}

\section{Introduction}

We consider the so-called Benjamin-Bona-Mahony equation, for $\lambda \in (0,1)$,
\begin{equation}\label{eq:BBM} \tag{BBM}
	(1-\lambda \partial_x^2) \partial_t u + \partial_x (\partial_x^2 u - u + u^2) =0,\quad t,x\in \RR. 
\end{equation}
We refer to Appendix \ref{reduction} below for obtaining  {\eqref{eq:BBM}} from 
the following more standard form of the equation
\begin{equation}\label{eq:BBMc}
	(1-\partial_x^2) \partial_t \mathcal{U} + \partial_x( \mathcal{U} + \mathcal{U}^2)=0.
\end{equation}
Recall that \eqref{eq:BBMc} was originally introduced  by Peregrine \cite{Pe} and Benjamin, Bona and Mahony \cite{BBM} as
an alternate model to the standard integrable \eqref{eq:KdV} equation, corresponding to $\lambda=0$,
\begin{equation}\label{eq:KdV} \tag{KdV}
	\partial_t u + \partial_x( \partial_x^2 u  + u^2)=0.
\end{equation}

The Cauchy problem for \eqref{eq:BBM} is globally well-posed in $H^1$ (see \cite{BBM}), and   any $H^1$ solution $u(t,x)$ of \eqref{eq:BBM} satisfies for all $t\in \RR$,
\begin{align}
	&	\int \left(\lambda (\partial_x u)^2 + u^2 \right)(t)  = M(u(t)) = M(u(0)) \qquad \text{(mass)}	\label{eq:i5}\\
	&	\int \left((\partial_x u)^2 + u^2 -\frac 23 u^3\right)(t)  = \mathcal{E}(u(t)) = \mathcal{E}(u(0))\qquad \text{(energy)}	\label{eq:i6}
\end{align}
It is also well-known that the \eqref{eq:BBM} equation has soliton solutions : for   $\mu>-1$,  set 
$$ Q_{\mu}(x) = (1+\mu) Q\left(\sqrt{\frac {1+\mu}{1+\lambda \mu}} \, x\right)$$
where
$$	Q(x)=\frac 32 \frac 1 {\cosh^2 \left(\frac x2\right)} \quad \text{solves} \quad
	Q''+Q^2=Q.
$$
Then, for any $\mu>-1$, $y\in \RR$, 
$
	R_{\mu,y}(t,x)=Q_{\mu}(x-\mu t-y)$ is solution of {\eqref{eq:BBM}}.

\subsection{Review on the collision problem for (KdV) type equations}

We briefly review some results concerning the problem of collision of solitons for (KdV) type models and we  refer to the introduction of \cite{MMkdv4} for  more details.

\medskip

 First, it is very well-known that  the \eqref{eq:KdV}  equation has explicit pure
$N$-soliton solutions ( \cite{HIROTA}, \cite{WT}, \cite{Miura}): for any given $c_1>\ldots>c_N>0$, $y_1^-,\ldots,y_N^-\in \mathbb{R}$, there exists an explicit multi-soliton solution $u(t,x)$ of \eqref{eq:KdV} which satisfies
\[ 
\lim_{t\to \pm \infty}	\biggr\|u(t)- \sum_{j=1}^N c_j Q(\sqrt{c_j}(.-c_jt-y_j^\pm))\biggr\|_{H^1(\RR)}   =  0,
\]
for some $y_j^+$ (such  solutions were found using the inverse scattering transform).

Stability and asymptotic stability of $N$-solitons were studied by Maddocks and Sachs \cite{MS} in $H^N$ by variational techniques   and in the energy space $H^1$ by Martel, Merle and Tsai \cite{MMT}.

\medskip

Second, recall that LeVeque \cite{Le} further investigated the behavior of the explicit $2$-soliton solution $u$ above in the asymptotic $\mu = \frac {c_1-c_2} {c_1+c_2}$  small i.e. for nearly equal solitons. 
It is proved that for some explicit functions $c_j(t)$, $y_j(t)$
\begin{equation}\label{eq:Le}
\sup_{t,x\in \RR} \left| u(t,x) - 
c_1(t) Q(\sqrt{c_1(t)}(x-y_1(t))) - c_2(t) Q(\sqrt{c_2(t)}(x- y_2 (t)))\right| \leq C \mu^2,
\end{equation}
Moreover,
$
	\min_{t\in \RR} (y_1(t)-y_2(t)) =  2 | \ln \mu | + O(1),
$
which means in particular that the minimum separation between the two solitons goes to $\infty$ as $\epsilon\to 0$.
See \cite{Le} for a more precise statement.

Collision problems for (gKdV) and (BBM) have also been studied since the 60's from both experimental and numerical points of view (see \cite{FPU}, \cite{KZ}, \cite{Z}, \cite{ABM}, \cite{E}, \cite{BPS}, \cite{KB}, \cite{MMC}, \cite{SHIH}, \cite{Craig}, \cite{WM}, \cite{HHGY}).

\medskip

However, except for some integrable equations for which special explicit solutions are known,
the problem of describing rigorously the collision of two solitons is mainly open.
Now, we review  some recent rigorous works related to the interaction of two solitons in the nonintegrable situation for the generalized KdV equations 
\begin{equation}\label{eq:gKdV} \tag{gKdV}
	\partial_t u + \partial_x( \partial_x^2 u + u^p)=0,\quad t,x\in \RR. 
\end{equation}
Recall that solitons of \eqref{eq:gKdV} write
$R_{c,y}(t,x)= c^{\frac 1{p-1}} Q(\sqrt{c} (x-ct-y))$, for $c>0$, $y\in \RR$ where $Q$ satisfies
$Q'' + Q^p = Q$. 

\medskip

Mizumachi \cite{Mi}   studied rigorously the interaction of two solitons of nearly equal speeds for \eqref{eq:gKdV} for $p=3$ and $p=4$. For initial data $u_0$  close to  $Q(x)+ c^{\frac 1{p-1}} Q(\sqrt{c}(x+L))$, where
$L>0$ is large and $c$, close to $1$, satisfies $c-1 \leq  e^{-\frac L2}$, Mizumachi proved 
the two solitons remain separated for all positive time and that eventually the corresponding solution $u(t)$ behaves as 
\begin{equation}\label{eq:mi}
u(t)= (c_1^+)^{\frac 1{p-1}} Q\left(\sqrt{c_1^+}(.-c_1^+ t -y_1^+)\right) + (c_2^+)^{\frac 1{p-1}} Q\left(\sqrt{c_2^+}(.-c_2^+ t -y_2^+) \right)+ \varepsilon(t,x),
\end{equation}
for large time, for some $c_1^+<c_2^+$ close to $1$ and $\varepsilon$ small in some space.
The analysis part in \cite{Mi} relies on scattering results due to Hayashi and Naumkin \cite{HN1,HN2} and on the use of spaces of exponentially decaying functions (introduced in this context by Pego and Weinstein \cite{PW}). 

From \cite{Mi}, the situation is roughly speaking similar to the one described in the integrable case by LeVeque \cite{Le}.
However, two main questions were left   open in this work in this regime

\emph{Is the  $2$-soliton structure stable globally in time   in the energy space $H^1$?}

\emph{Does there exist a pure $2$-soliton in this regime?}

\noindent As in the integrable case, we call \textit{pure $2$-solitons}, solutions of   \eqref{eq:gKdV}   satisfying
\begin{equation}\label{eq:pure}
u(t) - \sud (c_j^\pm)^{\frac 1{p-1}} Q\left(\sqrt{c_j^\pm}(.-c_j^\pm t -y_j^\pm\right)\to 0
\quad \text{as $t\to \pm \infty$ in $H^1(\RR)$.}
\end{equation}
Note that if \eqref{eq:pure} holds both at $-\infty$ and $+\infty$, then necessarily
$c_j^-=c_j^+$ for $j=1,2$ (see \cite{MMcol2}, pp. 68, 69).

\medskip

These two questions have been answered  in a recent work by the authors \cite{MMkdv4}. Indeed,  in the   context of two solitons of almost equal speeds for the quartic (gKdV) equation, by constructing an approximate solution to the problem, we were able to prove first the global stability of the two soliton structure in $H^1$ and second, the inelastic character of the interaction. See Theorems 1 and 2 in \cite{MMkdv4}. 

\medskip

We also point out some other recent works of the authors  (\cite{MMcol1}, \cite{MMcol2}) concerning the problem of collision of two solitons of \eqref{eq:gKdV}  for a general nonlinearity $g(u)$ in the case where one soliton,  is supposed to be large with respect to the other soliton, i.e. assuming $0<c_1 \ll c_2$. 
See also \cite{MMM}, with T. Mizumachi, extending these results to the (BBM) equation.

\subsection{Main results}

In the present paper, we extend the results of \cite{MMkdv4} to the (BBM) model.

There are two main motivations to consider these questions for the \eqref{eq:BBM} model:
first, the   structure of the \eqref{eq:BBM} equation is  close to the one of the \eqref{eq:KdV} equation but it cannot be considered as a perturbation of the \eqref{eq:KdV} equation.
Second, the  present paper on (BBM) proves that our techniques extend to quadratic nonlinearity, unlike \cite{Mi}, based on scattering techniques critical for $p=3$.

\begin{theorem}[Inelastic interaction of two solitons with nearly equal speeds]\label{TH:1}
	Let $\lambda \in (0,1)$.
There exist $C,c, \sigma, \mu_*>0$ such that the following holds.
For $0<\mu_0<\mu_*$, let $U(t)$ be the unique solution of {\eqref{eq:BBM}} such that
\begin{equation}\label{eq:th1}
	\lim_{t\to -\infty}\|U(t) - Q_{-\mu_0}(. + \mu_0 t + \tfrac 12 Y_0 + \ln 2) - Q_{\mu_0}(.-  \mu_0 t -\tfrac 12 Y_0-  \ln 2)\|_{H^1} = 0,
\end{equation}
where 
$Y_0 = |\ln(\mu_0^2/\alpha)|$  and $\alpha=240/(15+10\lambda - \lambda^2).$
Then 
\begin{enumerate}
	\item [\rm (i)] Global stability of  $2$-solitons. There exist $\mu_1(t),$ $\mu_2(t),$ $y_1(t),$ $y_2(t)$ such that,
	\begin{equation}
	w(t,x)= U(t,x) -  Q_{\mu_1(t)}(x-y_1(t)) - Q_{\mu_2(t)}(x-y_2(t))
	\end{equation}
	satisfies,  for all $t\in \RR$,
\begin{align}
	&	 \|w(t)\|_{H^1(\RR)} \leq C |\ln \mu_0|^{1/2} \mu_0^{2},\quad	\Big| \min_{t\in \RR} (y_1(t)-y_2(t)) -  Y_0 \Big|\leq   C |\ln \mu_0|^\sigma \mu_0^{3/2},\label{eq:th5} \\
	& \sum_{j=1,2}| \mu_j(t) + (-1)^j \mu_0 \tanh(\mu_0 t)|+\sum_{j=1,2} | \dot y_j(t) - \mu_j(t)|   \leq  C |\ln \mu_0|^2 \mu_0^{2}.\label{eq:th6}
	\end{align}
	\item [\rm (ii)] Asymptotics and defect.
	The limits $\mu_1^+=\displaystyle\lim_{ +\infty} \mu_1$, $\displaystyle \mu_2^+=\lim_{+\infty} \mu_2$ exist and
		\begin{align} 
			&	\lim_{t\to +\infty}\| w(t)\|_{H^1(x>-(99/100) t)}=0,
			\quad \liminf_{t\to +\infty }  \|w(t)\|_{H^1(\RR)}\geq c \mu_0^3,\label{eq:th8}\\  
					&		 c  \mu_0^5	 \leq   {\mu_1^+}  - \mu_0 \leq  C |\ln \mu_0|^{2\sigma} \mu_0^4 , \quad  c  \mu_0^5	 \leq   -{\mu_2^+}  - \mu_0 \leq  C |\ln \mu_0|^{2\sigma} \mu_0^4 . \label{eq:th9}
	\end{align}
	\end{enumerate}
\end{theorem}
It follows immediately from the lower bound \eqref{eq:th8} that
\emph{no pure $2$-soliton exists}, which is a new result for the (BBM) equation in this regime.

\begin{theorem}[Stability result in the energy space for (KdV) and (BBM) equations]\label{TH:2}
Let $\lambda \in [0,1)$.
There exists $\mu_*>0$, $C,\sigma>0$, such that the following holds.
Let $\tilde\mu_{\inch}\in \RR$ and $\tilde Y_{\inch}>0$ be such that
\begin{equation}\label{eq:th2a}
	\mu_0  = \left(\tilde\mu_{\inch}^2 + 4 \alpha e^{-\tilde Y_{\inch}}\right)^{1/2} <\mu_*.
\end{equation}
Let $u_{\inch}\in H^1$ be such that
\begin{equation}\label{eq:th2b}
	\| u_{\inch} - Q_{ -\tilde\mu_{\inch}}(.- \tfrac 12 \tilde Y_{\inch}) - Q_{\tilde \mu_{\inch}}(. +\tfrac 12  \tilde Y_{\inch}) \|_{H^1(\RR)}\leq \omega \mu_0,
\end{equation}
where $0<\omega<|\ln \mu_0|^{-2}$, and let $u(t)$ be the solution of \eqref{eq:BBM} such that $u(0)=u_{\inch}$. Then,
there exist $T(t)$, $X(t)$ of class $C^1$ such that, for all $t\in \RR$,
\begin{equation}\label{eq:th2c}
	\|u(t + T(t), . + X(t)) - U(t)\|_{H^1(\RR)} + |\dot X(t)|+ \mu_0 |\dot T(t)| \leq C \omega \mu_0 + C |\ln \mu_0|^\sigma \mu_0^{ 3 / 2 },
\end{equation}
where $U(t)$ is the solution defined in Theorem \ref{TH:1}.
\end{theorem}

\noindent\emph{Comments on the results:}

\medskip

1. The (KdV) case in Theorems \ref{TH:1} and \ref{TH:2}.
\quad 
The value $\lambda = 0$ in the BBM equation corresponds to the integrable  KdV equation. In this case,  estimates \eqref{eq:th5}--\eqref{eq:th6} still hold. Estimate \eqref{eq:th5} corresponds to \eqref{eq:Le} but from the proofs in the present paper,  we improve the main result in \cite{Le} in this case by computing  explicitely the term of size $\mu_0^2$, see Remark \ref{re:mm}. Note also that for  $\lambda=0$, the existence of pure $2$-soliton solutions corresponds to $\mu_1^+=\mu_0$ and $\mu_2^+ = - \mu_0$ in \eqref{eq:th8} and \eqref{eq:th9}.

Moreover,   Theorem \ref{TH:2} holds for $\lambda = 0$ and it is also a new global stability result for the \eqref{eq:KdV} equation  in the energy space.
This kind of result cannot be proved by scattering theory.

\medskip

2. Except for the value of the constant $\alpha>0$, Theorems \ref{TH:1} and \ref{TH:2} are exactly the same as for the quartic (gKdV) equation. In particular, the orders of size in $\mu_0$ in the various estimates do not depend on the power of the nonlinearity. Moreover, the function
\begin{equation}\label{eq:HHi}
\text{$\YYrr(t)= \YYzz + 2 \ln(\cosh ( \sqrt{ \alpha} e^{-\frac 12 {\YYzz}} t))$    solution of 
	$\ddot \YYrr= 2 \alpha e^{-\YYrr}, \ \lim_{t\to -\infty} \dot \YYrr(t) = 2 \mu_0$,  $\dot Y(0)=0$.}
\end{equation}
appears in both problems and has a universal character in this problem.
Note that Theorems \ref{TH:1} and \ref{TH:2} can be extended to  any two solitons of \eqref{eq:BBMc} of almost equal sizes using a simple scaling argument. See Appendix \ref{reduction}.

Finally, from the present paper and \cite{MMkdv4}, it is clear that the   results can be extended to  (gKdV) equations with general nonlinearities.

\medskip

3. As in \cite{MMkdv4}, the lower bounds in \eqref{eq:th8} and \eqref{eq:th9} measur the inelastic character of the collision. Moreover, the different exponents of $\mu_0$ in \eqref{eq:th8} and \eqref{eq:th9} denotes a gap in the estimates which is an open problem.

 \subsection{Strategy of the proofs}
 
We describe briefly the strategy of the proofs of Theorems \ref{TH:1} and \ref{TH:2}, which is  the same as in \cite{MMkdv4}. We point out the analogies and the main technical differences between the (BBM) and the quartic (gKdV) case.
The proof of Theorem \ref{TH:2} is a consequence of the proof of Theorem \ref{TH:1} and so we focus on the proof of Theorem \ref{TH:1}. Let $U(t)$ be the solution of (BBM)  defined in Theorem \ref{TH:1}.

\medskip

(1) The first step is the construction of an approximate solution in terms of a series in $e^{-y(t)}$ where $y(t)=y_1(t)-y_2(t)$ is the distance between the two solitons, using the exponential decay of the solitons. 
From Proposition~\ref{PR:ap}, the approximate solution contains a tail of order $e^{-y(t)}$ between the two solitons, which is relevant in the description of the exact solution, see Remark~\ref{re:mm}. 
This tail of order $e^{-y(t)}$  is not related to inelasticity since it appears also in the integrable case  $\lambda =0$.
Moreover, it does not prevent the approximate solution to be in the energy space at this order, since it is localized in space between the two solitons.

In contrast,  for $\lambda \neq 0$, one cannot build an approximate solution at order $e^{-\frac32 y(t)}$ in the energy space, whereas it is possible for $\lambda=0$. The presence of a nonzero tail at $-\infty$ in space at this order is   related to nonintegability and inelasticity.

The construction of the approximate solution for (BBM) in Section 2 is more involved that in  the quartic (gKdV) case mainly because   the nonlinearity is  quadratic rather that quartic. 

\medskip

(2) After the approximate solution is constructed, we introduce the following decomposition of the solution $U(t)$:
$$
U(t,x) =  Q_{c_1(t)}(x-y_1(t))+ Q_{c_2(t)}(x-y_2(t)) + W(t,x) + \varepsilon(t,x),
$$
where $Q_{c_1(t)}(x-y_1(t))+ Q_{c_2(t)}(x-y_2(t)) + W(t,x)$ is the modulated approximate solution and 
$\varepsilon(t)$ is a rest term.
To prove stability of the two soliton structure, we have  to control both  the parameters $c_j(t)$ and $y_j(t)$  and the rest term $\varepsilon(t)$.

From the construction of the approximate solution, the parameters $c_j(t)$ and $y_j(t)$ have to satisfy an approximate dynamical system. Remarkably, it is    exactly the same dynamical system as for the quartic (gKdV) equation  (except the values of the numerical constants).
This dynamical system, and the related solution $Y(t)$ of the ODE
$$ \ddot Y =  2 \alpha e^{-Y} , \quad Y(0)= Y_0, \quad \dot Y(0)=0,$$
seem to be  universal in this type of problems.
The control of the dynamical system satisfied by the parameters is thus exactly the same as in \cite{MMkdv4} and we will not repeat the arguments in the present paper (see Section 4). 

Concerning the control of the rest term $\varepsilon(t)$, as in \cite{MMkdv4},  we use variants of techniques developed  for  large time stability and asymptotic stability of solitons and multi-solitons for the (gKdV) equations in the energy space, \cite{We2}, \cite{MM1}, \cite{MMT} and \cite{Ma2}, extended to the (BBM) case in \cite{We3}, \cite{Mi2}, \cite{Di}, \cite{DI2}, \cite{DiMa} and \cite{Ma}.
At this point, we need some new refined arguments and   the proofs are more involved than in the quartic (gKdV) case.
Note that since the nonlinearity is quadratic, one cannot use scaterring theory from \cite{HN1}, \cite{HN2} as in \cite{Mi}.

\medskip

(3) Finally, in Section 5, we prove that for $\lambda \neq 0$, the defect due to the interaction of two solitons  is bounded from below, which implies in particular that the collision is not elastic.

Assuming for the sake of contradiction   that the lower bound in \eqref{eq:th8} is not satisfied for any positive value of $c$,  we obtain first some symmetry
 properties ($x\to -x$, $t\to -t$) on the parameters $c_j(t)$, $y_j(t)$
  at a certain order.

Second,  using space decay properties of $U(t,x)$, we obtain a gain in the control of the error term in the dynamical system satisfied by $c_j(t)$, $y_j(t)$.
Using this refined version of the dynamical system which is not  symmetric for $\lambda\neq 0$ (as a consequence of the tail of order $e^{-\frac 32 y(t)}$ in the approximate solution),  we find a contradiction.

\section{Construction of an approximate solution}\label{sec:ap}
We denote by ${\cal Y}$ the set of functions $f\in C^\infty(\mathbb{R},\mathbb{R})$ such that
$$
\forall j\in \mathbb{N}, \ \exists C_j, \ r_j>0, \ \forall x\in \mathbb{R},\quad
\left|f^{(j)}(x)\right|\leq C_j (1+|x|)^{r_j} e^{-|x|}.
$$

\begin{proposition}\label{PR:ap} Let $\lambda \in [0,1)$.
There exist unique $ A_j(x) $, $ B_j(x) $, $ D_j(x)$, $\alpha$, $\beta$, $\delta$, $a$, $b_{j}$, $d_{j}$  $(j=1,2)$, $\sigma\geq 3$ and $0<\mu_*<1/10$ such that for any $0<\mu_0<\mu_*$,  the following hold.
\begin{itemize}
\item[\rm (i)] Properties of $A_j,$ $B_j,$ $D_j$ and $b_j$.
\begin{equation}\label{eq:ABD}\begin{split}
&  A_j,B_j,D_j \in L^\infty(\mathbb{R}), \quad A_j',B_j',D_j'\in {\cal Y},\\
& -\lim_{\pm\infty} A_1= \lim_{\pm \infty} A_2=\pm \theta_A,\quad
\theta_A= \frac {36(5-\lambda^2)}{15+10\lambda -\lambda^2},\quad 
\lim_{+\infty} D_1=\lim_{+\infty} D_2=0,\\
& \lim_{+\infty} B_1=\lim_{+\infty} B_2=0,\quad
\lim_{-\infty} B_1=- \lim_{-\infty} B_2= -\frac {288 \lambda^2}{15+10 \lambda - \lambda^2},
\end{split}\end{equation}
and $A_j$, $B_j$ and $D_j$ satisfy the orthogonality conditions of Lemmas \ref{LE:wA}, \ref{LE:wB} and \ref{LE:wD}.
Moreover,
\begin{equation}\label{eq:b2b1}
	b_1=b_2 \quad \text{if and only if $\lambda = 0$.}
\end{equation}
\item[\rm (ii)]  Definition of the approximate solution.
For $\Gamma= (\mu_1,\mu_2,\yy_1,\yy_2)$,  define
\begin{equation}\label{eq:ap0}
\begin{split}
  V_0(x;\Gamma)	 &	=	Q_{\mu_1}(x-\yy_1)+ Q_{\mu_2}(x-\yy_2)\\ &	+ e^{-(y_1-y_2) } \left(A_1(x-\yy_1)+A_2(x-\yy_2)\right) \\ & 	+    {\theta}  (\mu_1-\mu_2)  \, x \, Q(x-\yy_1)Q(x-\yy_2)  \\
		&	+  (y_1-y_2) e^{-(y_1-y_2)} \left(\mu_1B_1(x-\yy_1)+\mu_2 B_2(x-\yy_2)\right) \\ &
			+  e^{-(y_1-y_2)} \left(\mu_1D_1(x-\yy_1)+\mu_2 D_2(x-\yy_2)\right),
\end{split}
\end{equation}
where
\begin{equation}\label{eq:dQ}
  {\theta}= \frac {(1+\lambda)(5-10 \lambda + \lambda^2)} {15+10 \lambda - \lambda^2} .
\end{equation}
\item[\rm (iii)]  Equation of $V_0(x;\Gamma(t))$.
Let $I$ be some time interval and  $\Gamma(t)=(\mu_1(t),\mu_2(t),\yy_1(t),\yy_2(t))$ be a $C^1$ function defined on $I$ such that, for some constant $K>1$, 
\begin{align} 
	\forall t\in I,\quad 
	& \YYzz-1\leq y_1(t)-y_2(t)\leq K \YYzz,
	\quad |\mu_1(t)|\leq 2 \mu_0,\quad |\mu_2(t)|\leq 2 \mu_0, \label{eq:mh1}\\
	& |\mu_1(t)+\mu_2(t)|\leq   \YYzz^2 e^{-\YYzz},\quad
	|y_1(t)+y_2(t)|\leq \YYzz^4 e^{-\frac 12 \YYzz}, \label{eq:mh2}
\end{align} 
where 	$$\YYzz=|\ln (\mu_0^2/\alpha)|\quad \text{and}\quad \alpha= {240}/{(15+10\lambda - \lambda^2)}.$$

Let 
\begin{equation}\label{eq:ap}
\begin{split}
V_0(t,x) = V_0(x;\Gamma(t)),\quad \pyy(t)=y_1(t)-y_2(t).
\end{split}
\end{equation}
Then, on $I$, $V_0(t,x)$ solves
\begin{equation}\label{eq:VV0} 
	(1-\lambda \partial_x^2) \partial_t V_0 + \partial_x (\partial_x^2 V_0 - V_0 + V_0^2) = \fell(V_0) +E_0(t,x) 
\end{equation}
where 
$$
\fell(V_0)=
 \sum_{j=1,2} (\dot \mu_j-{\cal M}_j)
(1-\lambda\partial_x^2) \frac {\partial V_0} {\partial \mu_j} 
- \sum_{j=1,2} (\mu_j - \dot \yy_j - {\cal N}_j) (1-\lambda\partial_x^2)
\frac {\partial V_0} {\partial y_j}
$$
\begin{equation}\label{eq:pa}
\begin{split}
& {\cal M}_1(t)= 
\alpha\,  e^{-\pyy(t)} +  \beta\,  \mu_1(t) \pyy(t) e^{-\pyy(t)} + \delta  \,\mu_1(t) e^{-\pyy(t)} ,\\
& {\cal M}_2(t)=  
- \alpha\,  e^{-\pyy(t)} -  \beta\,  \mu_2(t) \pyy(t) e^{-\pyy(t)} - \delta  \,\mu_2(t) e^{-\pyy(t)},\\
& {\cal N}_1(t)= 
a \,e^{-\pyy(t)} + b_1 \,\mu_1(t) \pyy(t) e^{-\pyy(t)} + d_1 \,\mu_1(t) e^{-\pyy(t)},\\
& {\cal N}_j(t)= 
a \,e^{-\pyy(t)} + b_2 \,\mu_2(t) \pyy(t) e^{-\pyy(t)} + d_2 \,\mu_2(t) e^{-\pyy(t)},
\end{split}
\end{equation}
and  for some $C=C(K)>0$,
\begin{equation}\label{eq:E0}
	\forall t\in I, \quad \sup_{x\in \RR} \left\{\left(1+e^{\frac 12(x-\yy_1(t))} \right)| E_0(t,x)|\right\} \leq C \YYzz^\sigma e^{-\YYzz} e^{-\pyy(t)}.
\end{equation}
\end{itemize}
\end{proposition}
Sections 2.1--2.5 are devoted to the proof of Proposition \ref{PR:ap}.

\medskip

Note that the function $V_0$ is not in $L^2$ since $B_j$ have non zero limits at $-\infty$.
We now introduce an $L^2$ approximation of $V_0$, using a suitable cut-off function. Let $\psi:\mathbb{R}\to [0,1]$ be a $C^\infty$ function such that 
\begin{equation}\label{def:psi}
\hbox{$\psi'\geq 0$, $\psi\equiv 0$ on $\mathbb{R}^-$, $\psi\equiv 1$ on $[\frac 12 , +\infty)$,}
\end{equation}
As a consequence of Proposition \ref{PR:ap}, we obtain the following result.

\begin{proposition}[$L^2$ approximate solution]\label{PR:app}
 Under the assumptions of Proposition \ref{PR:ap} (i)--(iii),   
let 
\begin{equation}\label{eq:V}
	V(x;\Gamma)= V_0(x;\Gamma) \psi\left( e^{-\frac 12 \YYzz}  x +1 \right),\quad
	V(t,x)=V(x;\Gamma(t)).
\end{equation}
Then,
\begin{itemize}
	\item[\rm (i)]  Closeness to the sum of two solitons.
	\begin{equation}\label{eq:ap2bis}
		\|V 	- \left\{Q_{\mu_1}(.-y_1) +Q_{\mu_2}(.-y_2) \right\}\|_{L^\infty}
		\leq  C  e^{-\pyy},
	\end{equation}
	\begin{equation}\label{eq:ap2}
		\|V 	- \left\{ Q_{\mu_1}(.-y_1)+ Q_{\mu_2}(.-y_2) \right\}\|_{H^1}
		\leq  C \sqrt{\pyy} e^{-\pyy}.
	\end{equation}
	\item[\rm (ii)]  Equation of $V(t,x)$. 
\begin{equation}\label{eq:VV}
	(1-\lambda \partial_x^2) \partial_t V + \partial_x (\partial_x^2 V - V + V^2) = \fell(V) + E(t,x) \\
\end{equation}
where 
\begin{equation}\label{eq:defell}
\fell(V) = \sum_{j=1,2} (\dot \mu_j-{\cal M}_j) (1-\lambda\partial_x^2)
\frac {\partial V}{\partial \mu_j} 
- \sum_{j=1,2}
 (\mu_j - \dot \yy_j - {\cal N}_j) (1-\lambda\partial_x^2)
 \frac {\partial V}{\partial y_j},
\end{equation}
and   for some $C=C(K)>0$,
\begin{equation}\label{eq:eE}\begin{split}
\forall t \in I,\quad &	\sup_{x\in \RR} \{\left(1+e^{\frac 12(x-\yy_1(t))} \right)| E(t,x)|\} \leq C \YYzz^\sigma e^{-\YYzz} e^{-\pyy(t)},\\
&	\|E(t)\|_{L^2} \leq C \YYzz^{\sigma} e^{-\frac 34 \YYzz} e^{-\pyy(t)}.
\end{split}\end{equation}
\end{itemize}
\end{proposition}

The proof of Proposition \ref{PR:app} being very similar to the one of Proposition 2.2 in \cite{MMkdv4}, it is omitted.

\subsection{Preliminary expansion}
We set
\begin{equation}\label{A.23}
\begin{split}
&    \qtud(t,x)= Q_{\mu_j(t)}(x-y_j(t)), \quad \qud(t,x)= Q (x-y_j(t)), \\
&
	 \Lambda \qtud(t,x)= \Lambda Q_{\mu_j(t)}(x-y_j(t)), \quad \Lambda \qud(t,x)  = \Lambda Q(x-y_j(t)),
\end{split}
\end{equation}
and similarly for $\Lambda^2 \qud$,
where $\Lambda Q_{\mu}$, $\Lambda^2 Q_\mu$  are defined in Claim \ref{LE:A1}.

We introduce the   notation
\begin{equation}\label{eq:19}\begin{split}
& \hbox{$r(t)=\OTT_{k}$, for $k\geq 1$, if  $\exists \sigma\geq 0$ s.t. } 
\sup_{t \in I}\{ e^{\pyy(t)} |r(t)|\}\leq C ( 1+ \YYzz^\sigma )  e^{-(k-1) \YYzz} ,\\
& \text{$f(t,x)=\OO_{k}$, for $k\geq 1$, if $  \sup_{x\in \RR} \left\{\left(1+e^{\frac 12 (x-y_1(t))}\right) |f(t,x)| \right\}=  O_k$.}
\end{split}\end{equation}

Define
\begin{equation}\label{eq:S}
{\cal S}(v) = 	(1-\lambda \partial_x^2) \partial_t v + \partial_x (\partial_x^2 v - v + v^2),
\end{equation}
and $\mathcal{M}_j$, $\mathcal{N}_j$ as in \eqref{eq:pa}, for $\alpha$, $\beta$, $\delta$ and $a$, $b_j$, $d_j$ to be determined.

We look for an approximate solution of ${\cal S}(v)=0$ under the form $v(t,x)=v(x;\Gamma(t))$,
\begin{equation}\label{70bis}
	v=\qtun + \qtde + w,
\end{equation}
where $w(t,x)=w(x;\Gamma(t))$ so that using the equation of $Q_\mu$ (see \eqref{eq:Qm}) and
$\frac {\partial}{\partial \mu_j} \qtud= \Lambda \qtud$, 
 $\frac {\partial}{\partial y_j} \qtud= -\px \qtud$,
\begin{equation}\label{eq:dc}
	{\cal S}(v)= \fell (v) + \FF + \FTF + G(w) + H(w),
\end{equation}
where
\begin{align*}
  \fell(v) & =  \sum_{j=1,2}  (\dot\mu_j-\mathcal{M}_j) (1-\lambda \partial_x^2) \frac {\partial v}{\partial \mu_j}
					-\sum_{j=1,2}  (\mu_j-\dot y_j-\mathcal{N}_j) (1-\lambda \partial_x^2) \frac {\partial v}{\partial y_j}
\\		  \FF & = 2 \px \left(\qtun \qtde\right),\\
	  \FTF &  = \mathcal{M}_1 (1-\lambda \partial_x^2)\Lambda \qtun + \mathcal{M}_2 (1-\lambda \partial_x^2)\Lambda \qtde + \mathcal{N}_1 (1-\lambda \partial_x^2)\px \qtun
	+ \mathcal{N}_2 (1-\lambda \partial_x^2)\px \qtde,
\end{align*}
and
\begin{align*}
G(w) & = \px \left[ \px^2 w - w + 2 \left(\qtun + \qtde\right) w\right]  + \sud \mu_j \frac {\partial w}{\partial y_j} 
\\
H(w) & = \px \left[ w^2\right] + \sud \mathcal{M}_j  (1-\lambda \partial_x^2)\frac {\partial w}{\partial \mu_j} 
- \sud \mathcal{N}_j  (1-\lambda \partial_x^2)\frac {\partial w}{\partial y_j} .
\end{align*}
In the rest of this section, we give preliminary expansions of $\FF$ and $\FTF$.

\begin{lemma}[Expansion of $\FF $]\label{LE:tF}
Under the assumptions of Proposition \ref{PR:ap},
\begin{equation*} 
 \FF  = 2 \px\left(\qtun \qtde\right) = \FF_A + \FF_Q + \FF_B + \FF_D +   \OO_{2},
 \end{equation*}
where
\begin{align*}
 \FF_A  & 
= - 12 e^{-\pyy} (\partial_x \qun+\qun)+ 
 12 e^{-\pyy} (-\partial_x \qde+\qde),\\
  \FF_Q& = (1-\lambda) (\mu_1-\mu_2) \qun \qde ,
 \\ 
\FF_B &=  - 6  (1-\lambda) \mu_1 \pyy e^{-\pyy}   (\partial_x \qun+\qun)
+ 6  (1-\lambda) \mu_2 \pyy e^{-\pyy}   (-\partial_x \qde+\qde)\\ 
\FF_D &= e^{-y} [ \mu_1 S_{F,1} (x-y_1) + \mu_2 S_{F,2}(x-y_2)],
\end{align*}
with $S_{F,1}\in \mathcal{Y}$ and $S_{F,2}(x) = - S_{F,1}(-x)$.
\end{lemma}

Note that the term $F_Q$ does not exist in the quartic case (see Lemma 2.1 in \cite{MMkdv4}).
The proof of Lemma  \ref{LE:tF} is given in Appendix \ref{AP:ZZ}.

\begin{lemma}[Expansion of $\FTF $]\label{LE:tl}
Under the assumptions of Proposition \ref{PR:ap},
\begin{equation*}
				\FTF =\FTF_A + \FTF_B  + \FTF_D +  \OO_2,
\end{equation*}
where
\begin{align*}
\FTF_A & =(1-\lambda \px^2)\big[ \alpha e^{-y} \Lambda \qun 
+ ae^{-y} \px \qun - \alpha e^{-y} \Lambda \qde + a e^{-y} \px \qde\big],\\
\FTF_B &= (1-\lambda \px^2)\big[ \beta \mu_1 y e^{-y} \Lambda \qun 
+ b_1 \mu_1 y e^{-y} \px \qun  - \beta \mu_2 y e^{-y} \Lambda \qde + b_2 \mu_2 y e^{-y} \qde\big]\\
\FTF_D & =(1-\lambda \px^2)\big[  \delta \mu_1   e^{-y} \Lambda \qun + \alpha \mu_1 e^{-y} \Lambda^2 \qun  + d_1 \mu_1   e^{-y} \px \qun+ a \mu_1 e^{-y}  \px \Lambda \qun 
 \\
& - \delta \mu_2  e^{-y} \Lambda \qde  - \alpha \mu_2 e^{-y} \Lambda^2 \qde + d_2 \mu_2   e^{-y} \px \qde + a \mu_2 e^{-y}  \px \Lambda \qde\big]
			 .\end{align*}
\end{lemma}
The proof of this result is the same as the one of Lemma 2.2 in \cite{MMkdv4}, thus it is omitted.

\subsection{Determination of $A_j$}
\begin{lemma}\label{LE:wA}
 Let
\begin{equation}\label{eq:aal}
	\alpha= \frac {240}{15+10 \lambda -\lambda^2},\quad
	\theta_A=  \frac {36 (5 -\lambda^2)}{15+10 \lambda - \lambda^2}.
\end{equation}
\begin{itemize}
\item[\rm (i)] 
There exist $a$ and  $\hat A_1 \in \mathcal{Y}$  such that 
$A_1 =\hat A_1 + \theta_A \frac {Q'}Q$ solves
\begin{equation*}
(-LA_1)' + 2 \theta_A Q' +  \alpha (1-\lambda \px^2)\Lambda Q 
+ a (1-\lambda \px^2) Q'= 12 (Q+Q'),
\end{equation*}
$$
\int A_1 (1-\lambda \px^2) Q' = \int (A_1+\theta_A) (1-\lambda \px^2)Q= 0.
$$
\item[\rm (ii)] Set $A_2(x)=A_1(-x)$ and
$$
w_A(t,x)=e^{-\pyy(t)} ( A_1(x-y_1(t))+ A_2(x-y_2(t))).
$$
Then,
\begin{equation}\label{eq:bab}\begin{split}
\FF_A+ \FTF_A + G(w_A)  &= 
   - \frac {\theta_A}{18} (\mu_1-\mu_2) \qun \qde\\
& +  e^{-\pyy} [\mu_1 S_1(x-y_1)+ \mu_2 S_2(x-y_2)] + \OO_2
\end{split}\end{equation}
where  $S_1\in \mathcal{Y}$, $S_2(x)=-S_1(-x)$.

Moreover,
\begin{equation}\label{ortho}
\sud\left|\int w_A (1-\lambda \px^2) \qud\right|+
\sud\left|\int w_A (1-\lambda \px^2) \px \qud \right|
= \OO_2.
\end{equation}  
\end{itemize}
\end{lemma}

\begin{proof}
Proof of (i). First, we determine  $\alpha$. Multiplying the equation of $A_1$ by $Q$, integrating and using $L (Q')=0$, we obtain by \eqref{eq:58} and \eqref{eq:qint3}
\begin{equation}\label{eq:a1}
	\alpha  \int  [ (1-\lambda \partial_x^2) \Lambda Q]Q = 12 \int Q^2 \quad
	\text{so that} \quad
	\alpha  = \frac {240}{15 + 10\lambda -\lambda^2}.
\end{equation}

Second, we find the value of $\theta_A$.
For $\theta_A$ to be chosen, set 
$
A_1= \theta_A  \frac {Q'}{Q} + \hat  A_1$, so that from \eqref{eq:qpq}, $\hat A_1$ has to satisfy
$$
(- L\hat A_1)'= -\alpha  (1-\lambda \px^2)\Lambda Q + 12 (Q+Q') - \theta_A (2 Q -\frac 53 Q^2)
- 2 \theta_A Q' - a (1-\lambda \px^2) Q'.
$$
To find $\hat A_1$ in $\mathcal{Y}$, we need
\begin{equation}\label{eq:gl}
 \theta_A = \frac 12 \int  \left(-\alpha  (1-\lambda \px^2) \Lambda Q+ 12 Q\right)=
\frac 12 \left(-\frac \alpha 2 (1+\lambda) + 12\right)  \int  Q = \frac {36 (5 -\lambda^2)}{15+10 \lambda - \lambda^2}.
\end{equation}
 
For this choice of  $\theta_A$, there exists $Z\in \mathcal{Y}$, $\int Z (1-\lambda \px^2) Q'=0$   such that
$$
Z'=
-\alpha  (1-\lambda \px^2)\Lambda Q + 12 (Q+Q') - \theta_A (2 Q -\frac 53 Q^2)
- 2 \theta_A Q'   .
$$
By Claim \ref{LE:A2}, there exists $A\in \mathcal{Y}$ such that $\int A(1-\lambda \px^2)Q'=0$ and
$-LA=Z$. Let  $\hat A_1 = A - a \Lambda Q \in \mathcal{Y}$. Then,  $\int \hat A_1 (1-\lambda \px^2) Q'=0$
and $-L \hat A_1 = Z - a (1-\lambda \px^2) Q.$
Finally, we uniquely choose $a$ such that  $\int (\hat A_1+\theta_A) (1-\lambda \px^2) Q = 0$

\medskip

Proof of (ii).  
 First, by the parity properties of $Q$, $A_2(x)=A_1(-x)$ satisfies
$$
		(-L A_2)' + 2 \theta_A Q' - \alpha (1-\lambda\px^2) \Lambda Q + a  (1-\lambda\px^2)Q' = 
		-12 (Q-Q') .
$$
Now, we compute $\FF_A+\FTF_A+G(w_A)$. Using \eqref{eq:sqQQ}, we have
\begin{align*}
	G(w_A) & = \px\left( \px^2 w_A - w_A + 2\left(\qun   + \qde \right) w_A\right) \\
	& +  2 \px \left(\left( \mu_1  \Lambda \qun + \mu_2  \Lambda \qde\right) w_A \right) + \mu_1 \frac {\partial w_A}{\partial y_1} + \mu_2 \frac {\partial w_A}{\partial y_2} + \OO_2.
\end{align*}

First, 
\begin{align*}
	& \px\left( \px^2 w_A - w_A + 2 \left(\qun  + \qde \right) w_A\right)\\
	& = e^{-y} \left(- L A_1 + 2 \theta_A Q  \right)'(x-y_1) 
	+ e^{-y} \px\left(2 \qun  (A_2(x-y_2)-\theta_A) \right) \\
	& + e^{-y} \left(- L A_2 + 2 \theta_A Q \right)'(x-y_2)
	+ e^{-y} \px\left(2 \qde  (A_1(x-y_1)-\theta_A) \right).
\end{align*}
Using the estimate 
\begin{equation}\label{eq:29}
	|A_2(x-y_2)-\theta_A|\leq C (1+|x-y_2|^\omega) e^{-(x-y_2)}
\quad \text{for $x>y_2$}
\end{equation}
and \eqref{eq:sqQt}, we have
$$
	e^{-y} \qun (A_2(x-y_2)-\theta_A) = \OO_2 \quad \text{and similarly} \quad
	e^{-y} \qde (A_1(x-y_1)-\theta_A) = \OO_2.
$$
Thus, using the expressions of $\FF_A$ and $\FTF_A$ in Lemmas \ref{LE:tF} and \ref{LE:tl}
and the equations of $A_1$ and $A_2$, we find
$$
	\FF_A+\FTF_A+\px\left( \px^2 w_A - w_A + 2 \left(\qun   + \qde \right) w_A\right)=\OO_2.
$$

Second, by similar arguments,
\begin{align*}
	   2 \px \left(\left( \mu_1  \Lambda \qun + \mu_2   \Lambda \qde\right) w_A \right) 
	 &=  2 \mu_1 e^{-y} \px \left(    \Lambda \qun (A_1(x-y_1)+\theta_A)\right)\\ &
	+ 2  \mu_2 e^{-y}  \px \left(    \Lambda \qde (A_2(x-y_2) + \theta_A) \right)  +\OO_2.
\end{align*}

Finally, we compute $\mu_1 \frac {\partial w_A}{\partial y_1} + \mu_2 \frac {\partial w_A}{\partial y_2}$.
We have
$$
	\frac {\partial w_A}{\partial y_1} = -w_A - e^{-y} A_1'(x-y_1),\quad 
	\frac {\partial w_A}{\partial y_2} = w_A - e^{-y} A_2'(x-y_2).
$$
Thus, using \eqref{eq:mh2},
\begin{align*}
	 \mu_1 \frac {\partial w_A}{\partial y_1} + \mu_2 \frac {\partial w_A}{\partial y_2}
	 & = -(\mu_1-\mu_2)  w_A - \mu_1 e^{-y} A_1'(x-y_1)  - \mu_2 e^{-y} A_2'(x-y_2)\\
	& = - \theta_A (\mu_1 - \mu_2) e^{-y} \left(\frac {\px \qun}\qun -\frac {\px \qde} \qde\right)
	\\ & - \mu_1 e^{-y} (2\hat A_1+A_1')(x-y_1)  - \mu_2 e^{-y} (-2 \hat A_2+ A_2')(x-y_2) +\OO_2.
\end{align*}
For this term,  we use Claim \ref{CL:pp} (see Appendix A.2), i.e.
\begin{equation*}
		 e^{-\pyy} \left( \frac {\px \qun}{\qun}   - \frac {\px \qde}{\qde} \right)  = \frac 1{18} \qun\qde+  \frac 13 e^{-\pyy} \left(   \qun  +  \qde\right) + \OO_{3/2}.
\end{equation*}

We obtain
\begin{align*}
	& \mu_1 \frac {\partial w_A}{\partial y_1} + \mu_2 \frac {\partial w_A}{\partial y_2}
	  = - \frac 1 {18} \theta_A (\mu_1 - \mu_2) \qun\qde
	\\ & - \mu_1 e^{-y} (\tfrac 23  \theta_A Q + 2\hat A_1+A_1')(x-y_1)  - \mu_2 e^{-y} (- \tfrac 23 \theta_A Q - 2 \hat A_2+ A_2')(x-y_2)+\OO_2.
\end{align*}
Combining these computations, we obtain 
\begin{align*}
	&\FF_A + \FTF_A + G(w_A)  = -\frac 1 {18} \theta_A (\mu_1-\mu_2) \qun \qde\\
	& + \mu_1 e^{-y} \left(  2 \px \left( \Lambda Q (A_1+\theta_A)\right) - \tfrac 2 3 \theta_A Q -( 2 \hat A_1 + A_1')\right)(x-y_1)\\
	& + \mu_2 e^{-y} \left(  2 \px \left(  \Lambda Q (A_2+\theta_A)\right) +\tfrac  2 3 \theta_A Q + (2\hat A_2-   A_2')\right)(x-y_2),
\end{align*}
so that 
$$
S_1= 2  (\Lambda Q(A_1+\theta_A))'  - \tfrac 23 \theta_A Q - 2 \hat A_1 - A_1'.
$$
Using \eqref{eq:29}, and $\int (A_1 + \theta_A) (1-\lambda \px^2) Q=0$, we have
$$
\int w_A (1-\lambda \px^2) \qun = e^{-y} \int \left[ A_1 (1-\lambda \px^2)Q + \theta_A (1-\lambda \px^2) Q\right] + \OO_2 = \OO_2,
$$
and similarly for the other scalar products in \eqref{ortho}.
\end{proof}

\subsection{Nonlocalized term of order $\OO_{3/2}$}
\begin{lemma}[Approximate solution at order $\OO_{3/2}$ with localized error tem]\label{LE:wALOC}
Let 
$$	w_Q=  \theta  (\mu_1-\mu_2)   x\qun\qde,\quad
\theta = 1-\lambda - \frac {\theta_A}{18} 
=\frac {(1+\lambda)(5-10 \lambda + \lambda^2)}{15+10\lambda -\lambda^2}.
$$
Then
\begin{equation*}\begin{split}
  G(w_Q)  &=-  \theta  (\mu_1-\mu_2) \qun \qde\\
&-		 	 		  
		  6  \theta  \mu_1  y e^{-y}    e^{-(x-y_1)} \left( 3  (\qun-\px \qun) -\px(\qun^2)  - \qun^2\right) \\
		  & 
		  +6  \theta  \mu_2  y e^{-y}   e^{(x-y_2)} \left( 3 (\qde+\px \qde) +\px(\qde^2) 
		  - \qde^2 \right) \\ &
		  + e^{-y}\left( \mu_1 \widetilde S_1(x-y_1)+\mu_2 \widetilde S_2(x-y_2)\right) + \OO_2,
\end{split}\end{equation*}
where $\widetilde S_1 \in \mathcal{Y}$ and $\widetilde S_1(x)= - \widetilde S_2(-x)$.
\end{lemma}
\begin{proof}
The proof is based on Claim \ref{LE:Q1Q2} in Appendix A.

First, arguing as in the proof of Lemma \ref{LE:wA}, we have
$$
	G(w_Q)= \px\left(\px^2 w_Q - w_Q + 2 \left(\qun  + \qde \right) w_Q\right) +\OO_2.
$$
Moreover,
since 
$x = \frac 12 (x-y_1 + x - y_2) + \frac 12 (y_1+y_2)$, using \eqref{eq:mh2},
we have 
$$w_Q =   \frac 12 \theta  (\mu_1-\mu_2) (x-y_1 + x - y_2)\qun\qde + \OO_2.$$

Therefore, using Claim \ref{LE:Q1Q2} and the asymptotics of $Q$ from \eqref{eq:as}, we get
\begin{align*}
	 G(w_Q) & 	   =  - \frac 12 \theta  (\mu_1-\mu_2)\px \big\{ - \px^2 ((x-y_1 + x - y_2)\qun\qde)  + (x-y_1 + x - y_2)\qun\qde 
	 \\ & \quad - 2 (\qun  + \qde) (x-y_1 + x - y_2)\qun\qde \big\}+\OO_2 \\
		 &= - \theta   (\mu_1-\mu_2) \qun \qde\\
&-		 	 		  
		  6  \theta   \mu_1  y e^{-y}    e^{-(x-y_1)} \left( 3  (\qun-\px \qun -\px(\qun^2)) - \qun^2\right) \\
		  & 
		  +6  \theta   \mu_2  y e^{-y}   e^{(x-y_2)} \left( 3 (\qde+\px \qde+\px(\qde^2))
		  - \qde^2 \right) \\&
		  + e^{-y}\left( \mu_1 \widetilde S_1(x-y_1)+\mu_2 \widetilde S_2(x-y_2)\right) + \OO_2, 		  
\end{align*}
where
$\widetilde S_1$ and $\widetilde S_2$ satisfy the desired conditions.

\end{proof}
  
\subsection{Determination of $B_j$ and $D_j$}

\begin{lemma}\label{LE:wB}
Let
\begin{equation*}\begin{split}
& Z(x) =   6 (1-\lambda) e^{-x} (Q-Q') + 6 \theta e^{-x}\left(3 (Q-Q' - (Q^2)' ) - Q^2\right)
\\&	\beta=  \frac {120 (1-\lambda)}{15 + 10 \lambda - \lambda^2},	\qquad
	\theta_B=  \frac {144 \lambda^2}{ 15 + 10 \lambda - \lambda ^2} .
\end{split}\end{equation*}
\begin{itemize}
\item[\rm (i)]
There exist unique $b_1$ and $\hat B_1\in \mathcal{Y}$ such that
$B_1= \hat B_1 + \theta_B \left(1+\frac {Q'}Q\right)$ satisfies
$$(- L B_1)' +  \beta (1-\lambda \partial_x^2) \Lambda Q
+  b_1  (1-\lambda \partial_x^2) Q' = Z $$
$$
\int B_1 (1-\lambda \px^2) Q' = \int B_1 (1-\lambda \px^2) Q= 0.
$$
\item[\rm (ii)]
There exist unique $b_2$ and $\hat B_2\in \mathcal{Y}$ such that
$B_2 = \hat B_2 - \theta_B \left(1+\frac {Q'}Q\right)$ satisfies
$$(- L B_2)'  -4 \theta_B Q' - \beta (1-\lambda \partial_x^2) \Lambda Q
+  b_2  (1-\lambda \partial_x^2) Q' = -Z(-x) $$
$$
\int B_2 (1-\lambda \px^2) Q' = \int (B_2-2 \theta_B) (1-\lambda \px^2) Q= 0.
$$
Moreover, 
\begin{equation}\label{eq:defect}
	 b_2-b_1 = 0 \quad \Leftrightarrow \quad \lambda = 0. 
\end{equation}
\item[\rm (iii)]
Set
$$
w_B(t,x) = \pyy e^{-\pyy(t)} (\mu_1(t) B_1(x-y_1(t))+ \mu_2(t)B_2(x-y_2(t)).
$$
Then,
\begin{align}
& \FF_A + \FTF_A + G(w_A) + G(w_Q) + \FF_B + \FTF_B + G(w_B)
\nonumber\\ 
& = e^{-y} \left[\mu_1 (S_1+ \widetilde S_1) (x-y_1) +\mu_2 (S_2 +\widetilde S_2)(x-y_2)\right] + \OO_2,
\label{eq:2.27}
\end{align}
\begin{equation}\label{eq:2.28}
\sud \left|\int w_B (1-\lambda \px^2) \qud\right|+
\sud \left|\int w_B (1-\lambda \px^2) \px\qud\right|= \OO_{5/2}.
\end{equation}
\end{itemize}
\end{lemma}

\begin{proof}
We follow the strategy of the proof of Lemma \ref{LE:wA}.
The only difference is that we now look for solutions $B_1$, $B_2$ both with  limit $0$ at $+\infty$.

\medskip

Proof of (i). We find the  value of $\beta$ from the equation of $B_1$ multiplied by $Q$,
using \eqref{eq:58} and \eqref{eq:A777},
\begin{align*}
& \beta    \int Q(1-\lambda \px^2) \Lambda Q   = \frac 3{10} \beta (15 +10 \lambda -\lambda^2) \\ &= \int ZQ  =  
   ( 3 (1-\lambda ) +9 \theta) \int e^{-x} Q^2 - 10 \theta \int e^{-x} Q^3   =  36 (1-\lambda) .
\end{align*}
Next, from \eqref{eq:56}, \eqref{eq:58}, \eqref{eq:qint3}, 
we have
$$
\int Z = 
6(1-\lambda) \int Q + 6 \theta \left(3 \int Q -  2 \int e^{-x} Q^2\right)
= 36\left( (1-\lambda) - \theta\right) = 2 \theta_A,
$$
and
we find $\theta_B$ by integrating   the equation of $B_1$ ($2\theta_B = \int (-LB_1)'$)
\begin{align*}
2 \theta_B & = - 3 \beta (1+\lambda) + 2 \theta_A = \frac {288 \lambda^2}{15 + 10 \lambda - \lambda^2}.\end{align*}
We now obtain the existence of  $\hat B_1\in \mathcal{Y}$ as in the proof of Lemma \ref{LE:wA}, with $b_1$ uniquely chosen so that
$\int B_1 (1-\lambda \px^2) Q=0$ and $\int B_1 (1-\lambda \px^2) Q'$=0.

\medskip

Proof of (ii). We solve the equation of $B_2$ exactly in the same way. We check that the values of $\beta$ and 
$\theta_B$ are suitable to solve the problem, and we obtain  unique $\hat B_2\in \mathcal{Y}$ and $b_2$ so that $\int B_2 (1-\lambda \px^2)Q' =\int (B_2 - 2 \theta_B) (1-\lambda \px^2)Q=0$. 

\medskip

We now check that $b_1\neq b_2$ for $\lambda \neq 0$.
Let $B(x)=\hat B_1(x) - \hat B_2(-x)= B_1(x)-B_2(-x) - 2 \theta_B$. Then $B\in \mathcal{Y}$ and
$$
- L B +(b_1-b_2) (1-\lambda \px^2) Q  + 8 \theta_B Q = 0,\quad
\int (B+ 4 \theta_B) (1-\lambda \px^2) Q=0. 
$$
Then, multiplying the equation of $B$ by $\Lambda Q$, integrating and using
 $L\Lambda Q = - (1-\lambda \px^2) Q $ (see Claim \ref{LE:A1} (ii)), we obtain 
$$
\int B (1-\lambda \px^2) Q + (b_1-b_2) \int (1- \lambda \px^2) Q \Lambda Q
+ 8 \theta_B \int Q \Lambda Q= 0,
$$
and so,  by $\int Q\Lambda Q = \frac 14 (\lambda +3) \int Q $,
$$
(b_1-b_2) \int (1- \lambda \px^2) Q \Lambda Q =  4 \theta_B \left( \int Q - 2 \int Q\Lambda Q\right) = -2 \theta_B (\lambda +1),
$$
and thus, in view of the expression of $\theta_B$, we obtain
$$b_1=b_2 \quad \Leftrightarrow \quad \lambda = 0.$$

Proof of (iii). We finish the proof of Lemma \ref{LE:wB} as the one of Lemma \ref{LE:wA}.
In particular, using the limits of $B_1$ and $B_2$ at $\pm \infty$ 
\begin{align*}
	G(w_B) & = \mu_1 ye^{-y} (-LB_1)'(x-y_1) + \mu_2 y e^{-y} (-LB_2 - 4 \theta_B Q)' (x-y_2) +\OO_2.
\end{align*}
This, combined with the equations of $B_1$ and $B_2$ and Lemmas \ref{LE:tF}, \ref{LE:tl}, \ref{LE:wA} and \ref{LE:Q1Q2} proves \eqref{eq:2.27}. Note that $w_B$ is not in $L^2$ since it has a nonzero limit at $-\infty$.
However, it has exponential decay as $x\to +\infty$. This allows us to prove that all rest terms are indeed of the form $\OO_2$ (see notation $\OO_2$ in \eqref{eq:19}).

The control of the various scalar products is easily obtained as in Lemma \ref{LE:wA} from the properties of $B_1$, $B_2$.
\end{proof}

Finally, we claim without proof the following result.

\begin{lemma}[Definition and equation of $w_D$]\label{LE:wD}
Let
$$
S= - S_{F,1} - S_1 - \widetilde S_1 - (1-\lambda \px^2) \left( \alpha \Lambda^2 Q + a (\Lambda Q)'\right).
$$
\begin{itemize}
\item[\rm (i)] There exist unique $\delta$, $\theta_D$, $d_1$ and $\hat D_1\in \mathcal{Y}$ such that
$D_1= \hat D_1 + \theta_D \left(1+\frac {Q'}Q\right)$ satisfies
$$
	(-L D_1)' +  \delta (1-\lambda \partial_x^2) \Lambda Q
+  d_1  (1-\lambda \partial_x^2) Q' = S(x)
$$
\item[\rm (ii)] There exist unique $d_2$ and $\hat D_2\in \mathcal{Y}$ such that
$D_2 = \hat D_2 - \theta_D \left(1+\frac {Q'}Q\right)$ satisfies
$$ (-L D_2)' - \delta  (1-\lambda \partial_x^2) \Lambda Q
- 4 \theta_D Q'+ d_2 (1-\lambda \partial_x^2)Q' = - S(-x)$$
\item[\rm (iii)]
Set
$$
w_D(t,x) =  e^{-\pyy(t)} (\mu_1(t) D_1(x-y_1(t))+ \mu_2(t)D_2(x-y_2(t)).
$$
Then,
\begin{equation*}
 \FF_A + \FTF_A + G(w_A) + G(w_Q) + \FF_B + \FTF_B + G(w_B)
+ \FF_D + \FTF_D + G(w_D) = \OO_2,
\end{equation*}
\begin{equation*}
\sud \left|\int w_D (1-\lambda \px^2) \qud\right|+
\sud \left|\int w_D (1-\lambda \px^2) \px\qud\right|= \OO_{5/2}.
\end{equation*}
\end{itemize}
\end{lemma}

We do not need to compute $d_1-d_2$, this is the reason why the exact expression of $S_1$ and $\widetilde S_1$ are not needed.

\subsection{End of the proof of Proposition \ref{PR:ap}}
Set 
$$
V_0 = \qtun + \qtde + W_0, \quad W_0=w_A+w_Q+w_B+w_D.
$$
From the preliminary expansion \eqref{eq:dc}, we have
$$
\mathcal{S}(V_0) = \fell(V_0) + E_0, \quad E_0=\FF + \FTF + G(W_0) + F(W_0).
$$
In view of  notation \eqref{eq:19}, estimate \eqref{eq:E0} holds true for some $\sigma>0$ provided that 
$E_0=\OO_2.$
From Lemmas \ref{LE:wA}, \ref{LE:wALOC}, \ref{LE:wB} and \ref{LE:wD}, we have $\FF + \FTF + G(W_0) =\OO_2$.
Thus, we only have to check that $H(W_0)=\OO_2$.

First, 
$
 \px \left[ W_0^2\right]=\OO_2.
 $
Second, since $|\mathcal{M}_j |+ |\mathcal{N}_j|\leq C e^{-y}$,
we also obtain 
$$ \sud {\cal M}_j \frac {\partial W_0}{\partial \mu_j}	- \sud   {\cal N}_j \frac {\partial W_0}{\partial y_j}
 =\OO_2.
$$
Thus, Proposition \ref{PR:ap} is proved.

\section{Preliminary long time stability arguments}

\subsection{Stability of the $2$-soliton structure in the interaction region}

We start with the decomposition of any solution of \eqref{eq:BBM} close the approximate solution $V$ (introduced in Proposition \ref{PR:app}) by modulation theory.
See Appendix B for the proof.

\begin{lemma}[Decomposition around the approximate solution]\label{PR:de}
There exists $\omega_0>0$, $C>0$, $\bar y_0>0$ such that if
$u(t)$ is a solution of {\eqref{eq:BBM}} on some time interval $I$ satisfying
for $0<\omega<\omega_0$, $y_0>\bar y_0$
\begin{equation}\label{close}
	\forall t\in I,\quad 
	\inf_{y_1-y_2>y_0} \|u(t) - V(.;(0,0,y_1,y_2))\|_{H^1}\leq\omega,
\end{equation}
then there exists a unique decomposition $(\Gamma(t),\varepsilon(t))$ of  $u(t)$ on $I$,
\begin{equation}
	u(t,x)=V(x;\Gamma(t)) + \varepsilon(t,x), \quad
	\Gamma(t)=(\mu_1(t),\mu_2(t),y_1(t),y_2(t)) \text{ of class $C^1$},
\end{equation}
such that $\forall t\in I$,
\begin{equation}\label{eq:or}
\begin{split}
	&	\int \varepsilon(t,x) (1-\lambda \partial_x^2) \qtud(t,x) dx  = \int \varepsilon(t,x) (1-\lambda \partial_x^2) \partial_x \qtud(t,x) dx = 0,\\
	& \pyy(t)=y_1(t)-y_2(t)> y_0- C\omega, \quad \|\varepsilon(t)\|_{H^1}+|\mu_1(t)|+|\mu_2(t)|\leq C\omega ,
\end{split}
\end{equation}
\begin{equation}\label{eq:ep}\begin{split}
	& (1-\lambda \partial_x^2) \partial_t \varepsilon 
	+ \partial_x(\partial_x^2 \varepsilon - \varepsilon + 2 V \varepsilon + \varepsilon^2 ) + E(t,x) + \fell(V) = 0,
\end{split}\end{equation}
where $\qtud(t,x) = Q_{\mu_j(t)} (x-y_j(t))$ and $V$, $E(t,x)$, $\fell(V)$ are defined in Proposition \ref{PR:app}.

Moreover, assuming
\begin{equation}\label{eq:gga}
\forall t\in I,\quad 	(|\mu_1(t)|+|\mu_2(t)|) \pyy(t) \leq 1,
\end{equation}
 $\dot \Gamma(t)$ satisfies the following estimates
\begin{equation}\label{eq:ga}
\begin{split} 
& |  \dot \mu_j - {\cal M}_j   | \leq C  \Big[\|\varepsilon\|_{L^2}^2 + \pyy e^{-\pyy} \|\varepsilon\|_{L^2} +  \int |E| (\qtun+ \qtde) \Big], \\
&|\mu_j - \dot \yy_j - {\cal N}_j | \leq 
C \Big[\|\varepsilon\|_{L^2}+ \int |E| (\qtun+ \qtde) \Big].
\end{split}\end{equation}
\end{lemma}

The next proposition presents almost monotonicity laws which are essential in proving long time stability  results in the interaction region. They will allow us to compare the approximate solution $V(t,x)$ with exact solutions. The functional is different depending on whether $\mu_1(t)>\mu_2(t)$ or $\mu_1(t)<\mu_2(t)$.

The constant $0<\rho<1/32$  to be fixed later, set
\begin{equation}\label{eq:ph}
\begin{split}
	&\varphi(x)=\frac 2 \pi \arctan(\exp(8 \rho x)), \quad \text{so that }
        \lim_{-\infty} \varphi=0,\ \lim_{\infty} \varphi=1,\\
        & \forall x\in \mathbb{R},\quad \varphi(-x)=1-\varphi(x),\quad
 \varphi'(x)=\frac {8\rho}{  \pi\cosh( 8 \rho x)},\\
& |\varphi''(x)|\leq 8 \rho |\varphi'(x)|, \quad  |\varphi'''(x)|\leq (8 \rho)^2 |\varphi'(x)|.
\end{split}\end{equation}

\begin{proposition}[Almost monotonicty laws]\label{PR:cFG} For $\rho>0$ small enough,
 and under the assumptions of Lemma \ref{PR:de},
let
\begin{equation}\label{eq:en}
	{\cal F}_+(t) =
		\int \left[(\partial_x \varepsilon)^2 + \varepsilon^2 
		- \frac 23 \left((\varepsilon+V)^3 - V^3 - 3 V^2\varepsilon \right)\right]
		+\int \left[\lambda (\partial_x \varepsilon)^2 + \varepsilon^2\right] \Phi(t,x),
\end{equation} 
 where $\Phi(t,x)=\mu_1(t) \varphi(x)+ \mu_2(t) (1-\varphi(x));$
\begin{equation}\label{eq:eG}
	{\cal F}_-(t) =
		\int \left[(\partial_x \varepsilon)^2 + \varepsilon^2 - \frac 23 \left((\varepsilon+V)^3 - V^3 - 3 V^2\varepsilon \right)\right]\Phi_1(t,x)+\int \left[\lambda (\partial_x \varepsilon)^2 + \varepsilon^2\right]\Phi_2(t,x),
\end{equation}
where 
\begin{equation}\label{eq:cG0} 
    \Phi_1(t,x)=\frac {\varphi(x)}{(1+\mu_1(t))^2} + \frac {1-\varphi(x)}{(1+\mu_2(t))^2} ,\quad
  \Phi_2(t,x)=\frac {\mu_1(t)  \varphi(x)}{(1+\mu_1(t))^2} + \frac {\mu_2(t) (1-\varphi(x))}{(1+\mu_2(t))^2}.
\end{equation}
There exists $C>0$ such that
\begin{equation}\label{eq:FGcoer}
\|\varepsilon(t)\|_{H^1}^2 \leq C {\cal F}_+(t),\qquad \|\varepsilon(t)\|_{H^1}^2 \leq C {\cal F}_-(t).
\end{equation}
Moreover,

\noindent {\rm (i)}  If $t\in I$ is such that 
\begin{equation}\label{eq:en1}
\mu_1(t)\geq \mu_2(t)\quad  \text{and} \quad y_2(t)\leq - \frac 14 \pyy(t),\quad
	y_1(t)\geq \frac 14 \pyy(t),
\end{equation}
then
\begin{equation}\label{eq:cF}
\begin{split}
\frac d{dt}{\cal F}_+(t) & \leq  C \|\varepsilon\|_{L^2}^2 
\left[e^{- \frac 34  \pyy}   + (|\mu_1|+|\mu_2|+\|\varepsilon\|_{L^2} )( e^{-2 \rho y} +  \|\varepsilon\|_{L^2} )  \right]
+ C \|\varepsilon\|_{L^2} \|E\|_{L^2}.
\end{split}\end{equation}
 {\rm (ii)}  
	If $t\in I$ is such that 
	\begin{equation}\label{eq:en2}
		\mu_2(t)\geq \mu_1(t) \quad \text{and}\quad  y_2(t)\leq -\frac 14 \pyy(t),\quad
		y_1(t)\geq \frac 14 \pyy(t),
	\end{equation}
	then
	\begin{equation}\label{eq:cG}
	\begin{split}
	\frac d{dt}{\cal F}_-(t) & \leq  C \|\varepsilon\|_{L^2}^2 
	\left[  e^{- \frac 34  \pyy} +(|\mu_1|+|\mu_2|+\|\varepsilon\|_{L^2} ) (e^{-2 \rho \pyy}       
	+  \|\varepsilon\|_{L^2}) \right]
	+ C \|\varepsilon\|_{L^2} \|E\|_{L^2}.
	\end{split}\end{equation}
\end{proposition}

See proof of Proposition \ref{PR:cFG} in Appendix \ref{AP:B}.

\begin{remark}
 The introduction  of almost monotone variants of the energy  and mass is related to Weinstein's approach for stability of one soliton \cite{We3} and to   Kato identity for the \eqref{eq:gKdV} equation (see \cite{Kato2}). These techniques have been developed in \cite{MM1}, \cite{MMT}  and then extended in \cite{MMas2}, \cite{Di}, \cite{Mi} and \cite{DiMa}.
\end{remark}

\subsection{Stability of the two soliton structure for large time}

In this section, we present a stability result for the two soliton structure for large time, i.e. far away from the interaction time. The argument, similarly to the one of Propositions \ref{PR:cFG}, is based on almost monotone variant of energy and mass.
As a corollary, we obtain a sharp estimate for large negative time on the pure two solution solution considered in Theorem \ref{TH:1}.

\begin{proposition}[Stability for large time]\label{PR:STAB}
For $0<\rho<1/32$ small enough, there exist  $C>0$ and  such that for
$ \mu_0>0$ and $\omega>0$ small enough, if $u(t)$ is an $H^1$ solution of {\eqref{eq:BBM}} satisfying
\begin{equation}\label{eq:clo}
\	\left\|u(t_0) - Q_{-\mu_0}(. + \mu_0 t_0) - Q_{\mu_0}(.-\mu_0 t_0)\right\|_{H^1(\mathbb{R})} \leq \omega \mu_0,
 \end{equation}
 for some $t_0<- (\rho \mu_0)^{-1} |\log \mu_0|$, then 
there exist $y_1(t)$, $y_2(t)$ and $\mu_1^+$, $\mu_2^+$ such that
\begin{enumerate}
 \item[\rm (i)]  For all $t_0 \leq  t \leq -  (\rho \mu_0)^{-1} |\log \mu_0|$,
\begin{equation}\label{eq:stab+}\begin{split}
	& \left\|u(t)- Q_{-\mu_0}(. - y_1( t ) ) - Q_{\mu_0}(.- y_2( t ) )\right\|_{H^1(\mathbb{R})}
	\leq   C \omega\mu_0 +  C \exp\left({-4 \rho { \mu_0  }   |t|}\right),\\
	& y_1(t)-y_2(t)\geq \frac 32 \mu_0 |t|, \\
	& |-\mu_0 - \dot y_1(t)| + |\mu_0 - \dot y_2(t)|\leq C  \omega \mu_0+  C \exp\left({-4 \rho { \mu_0  }   |t|}\right).
\end{split} \end{equation}
 \item[\rm (ii)]  For all $t\leq t_0$,
\begin{equation}\label{eq:stab-}\begin{split}
	& \left\|u(t)- Q_{-\mu_0}(. - y_1( t ) ) - Q_{\mu_0}(. - y_2( t ) )\right\|_{H^1(\mathbb{R})}
	\leq   C \omega\mu_0 +  C \exp\left({-4 \rho { \mu_0  }   |t_0|}\right),
	\\
	& y_1(t)-y_2(t)\geq \frac 32 \mu_0 |t|, \\
	& |-\mu_0 - \dot y_1(t)| + |\mu_0 - \dot y_2(t)|\leq C  \omega \mu_0+  C \exp\left({-4 \rho { \mu_0  }   |t_0|}\right).
\end{split} \end{equation}
 \item[\rm (iii)]  Asymptotic stability.
 \begin{equation}\label{eq:as2}\begin{split}
&	\lim_{t\to -\infty} \left\|u(t)- Q_{\mu_1^+}(. - y_1(t)) - Q_{\mu_2^+}(.- y_2(t) )\right\|_{H^1(x<\frac {99}{100} |t|)}=0,\\
&	\lim_{t\to -\infty} \dot y_1(t) = \mu_1^{+}, \quad \lim_{t\to -\infty} \dot y_2(t)=\mu_2^{+},\\
&	|\mu_1^+ + \mu_0|+|\mu_2^+ -\mu_0|\leq C \omega\mu_0 +  C \exp\left({- 4 \rho { \mu_0  }|t_0|}\right).
	\end{split}
 \end{equation}
 \end{enumerate}
\end{proposition}

See the proof of this result in Appendix \ref{AP:B}.

\begin{remark}\label{RE:STAB}
Using the invariance of the BBM equation by the transformation
\begin{equation}\label{eq:inv}
x\to -x,\quad t\to -t,
\end{equation}
it follows that a statement similar to Proposition \ref{PR:STAB} holds for $t_0>(\rho \mu_0)^{-1} |\log \mu_0|$.
\end{remark}

\begin{corollary}\label{COR:pur}
Let $u(t)$ be the  unique solution of {\eqref{eq:BBM}} satisfying 
$$
\lim_{t\to -\infty}
\left\|u(t)- Q_{-\mu_0}( . + \mu_0t  ) - Q_{\mu_0}( . - \mu_0t  )\right\|_{H^1} = 0.
$$
Then, for all $	 t\leq  - (\rho \mu_0)^{-1} |\log \mu_0|, $
\begin{equation}\label{eq:pp}
 \left\|u(t)- Q_{-\mu_0}( . + \mu_0t  ) - Q_{\mu_0}( . - \mu_0t  )\right\|_{H^1}
	\leq    \exp\left({- 4 \rho \mu_0 |t| }\right). 
 \end{equation}
\end{corollary}
We refer to Theorem 1 in \cite{DiMa} for the existence and uniqueness of the solution $u(t)$.

\begin{proof}[Proof  of Corollary \ref{COR:pur} assuming Proposition \ref{PR:STAB}]
For fixed $t$, we  can pass to the limit $\omega\to 0$, $t_0\to -\infty$ in \eqref{eq:stab+}.
Then, we integrate the estimates on $\dot y_1(t)$ and $\dot y_2(t)$ (see \eqref{eq:stab+}) from $-\infty$ to $t$.
\end{proof}

\section{Stability of  the $2$-soliton structure}

In this section, using the approximate solution constructed in Propositions \ref{PR:ap} and \ref{PR:app} and 
the asymptotic arguments of Section 3, we prove the stability part of Theorem \ref{TH:1} and Theorem~\ref{TH:2}.

\subsection{Description of the global behavior of the asymptotic $2$-soliton solution}

Let  $0<\rho<1/32$ being fixed as in Propositions \ref{PR:cFG} and \ref{PR:STAB}.
Recall that $\sigma\geq 3$ is defined in Propositions \ref{PR:ap} and \ref{PR:app}.

We recall the following notation from the introduction
\begin{equation}\label{eq:defH}
\YYzz=|\ln(\mu_0^2/ \alpha)|\quad \text{or equivalently}\quad 
\mu_0 = \sqrt{\alpha} e^{-\frac 12{\YYzz}},
\end{equation}
\begin{equation}\label{eq:HH}
\text{$\YYrr(t)= \YYzz + 2 \ln(\cosh (\mu_0 t))$    solution of 
	$\ddot \YYrr= 2 \alpha e^{-\YYrr}, \ \lim_{t\to -\infty} \dot \YYrr(t) = -2 \mu_0$,  $\dot Y(0)=0$.}
\end{equation}
Note  that $\dot Y(t) = 2 \mu_0 \tanh(\mu_0 t)$ and,  for all $t\in \RR$,
\begin{equation}\label{eq:eH}
	0 \leq  \YYrr( t ) - \left(\YYzz + 2\mu_0  |t|  - 2 \ln 2 \right)\leq 2 \exp\left(-2\mu_0  |t|\right).
\end{equation}

\begin{proposition}[Description of the $2$-soliton solution in the interaction region]\label{pr:st}  \quad\\
Let  $U(t)$ be the  unique solution of {\eqref{eq:BBM}} such that
\begin{equation}\label{eq:pur2}
	\lim_{t\to -\infty} 
	\left\|U(t) - Q_{-\mu_0}(. + \tfrac 12 \YYrr(t) ) - Q_{\mu_0}(. - \tfrac 12 \YYrr(t) )\right\|_{H^1(\mathbb{R})} =0.
\end{equation}
Let $\TSR>0$ be such that $\YYrr(\TSR)=  400  \rho^{-2} \YYzz$.
Then, for $\mu_0>0$ small enough, there exists   $(\Gamma(t),\varepsilon(t)) \in C^1$ such that for all $t\in [-\TSR,\TSR],$
$$
	U(t,x)= V(t;\Gamma(t)) + \varepsilon(t,x),\quad \Gamma(t)=(\mu_1(t),\mu_2(t),y_1(t),y_2(t)),
$$
and
\begin{align}&
	|\bar \mu(t)|\leq   \YYzz^2 e^{- \YYzz}, \quad	|\bar y(t)| \leq   \YYzz^4 e^{-\frac 12 \YYzz},
\label{eq:J11} \\
& 
  	|\mu(t)-\dot \YYrr(t)|\leq C \YYzz^{\sigma+1} e^{- \frac 54  \YYzz},\quad
	|\pyy(t)-\YYrr(t)|\leq C \YYzz^{\sigma+2} e^{- \frac 34  \YYzz}, \label{eq:U+-} \\
	& \|\varepsilon(t)\|_{H^1} \leq C \YYzz^{\sigma} e^{-\frac 54 \YYzz},\label{eq:U++} 
	\end{align}
where 	
	\begin{equation}\label{eq:mth}\begin{split}&
	\mu(t)=\mu_1(t)-\mu_2(t), \quad \pyy(t)=y_1(t)-y_2(t),\\ & \summu(t)=\mu_1(t)+\mu_2(t),\quad
	\pyyb(t)=y_1(t)+y_2(t).
	\end{split}\end{equation}
Moreover, there exists $t_0$ such that
\begin{equation}\label{eq:t0}
	|t_0|\leq C \YYzz^{\sigma}e^{- \frac 14 \YYzz}, \quad \text{$\mu(t_0)=0;$ \ $\forall t\in [-\TSR,t_0)$,\  $\mu(t)<0;$ \ $\forall t\in (t_0,\TSR]$,\  $\mu(t)>0$.} 
\end{equation}
\end{proposition}

The proof of Proposition \ref{pr:st} is omitted since it is exactly the same as the one of Proposition~4.1 in \cite{MMkdv4}, using Sections 2 and 3.

\subsection{Conclusion of the proof of the stability of the $2$-soliton structure}

In this section, we finish the proof of the stability part of Theorem \ref{TH:1}.

\medskip

\noindent\emph{Proof of \eqref{eq:th5}--\eqref{eq:th6} and partial proof of \eqref{eq:th8} and \eqref{eq:th9}.}
Let $T>0$ be defined as in Proposition \ref{pr:st}. We prove the existence of $\mu_j(t)$ and $y_j(t)$ and estimates  \eqref{eq:th5}--\eqref{eq:th6} separately on $(-\infty,-T]$, $[-T,T]$ and $[T,+\infty)$. It is straightforward that the functions $\mu_j(t)$ and $y_j(t)$ can be ajusted to have $C^1$ regularity on $\RR$.

\smallskip

For $t<-T$, Corollary \ref{COR:pur} clearly implies \eqref{eq:th5}--\eqref{eq:th6}.

\smallskip

On $[-T,T]$, \eqref{eq:th5}--\eqref{eq:th6} are direct consequences of \eqref{eq:J11}--\eqref{eq:U++} and \eqref{eq:ap2} (comparing in $H^1$ the approximate solution with the sum of two solitons). 

\begin{remark}\label{re:mm}
By \eqref{eq:U++} and the definition of $V$
(see \eqref{eq:ap0} and \eqref{eq:V}), for $t\in [-T,T]$,
\begin{equation}\label{eq:rel}
\|U(t) - \qtun(t) - \qtde(t) - e^{-y(t)} (A_1(.-y_1(t))+A_2(.-y_2(t)))\|_{H^1}
\leq C Y_0^\sigma e^{-\frac 54 Y_0},
\end{equation}
where for $t$ close to $0$, the term
$e^{-y} (A_1(x-y_1)+A_2(x-y_2))$ is indeed relevant as a correction term in the computation of $U(t)$. In view of the behavior at $\pm \infty$ of the functions $A_1$ and $A_2$  (see Lemma \ref{LE:wA}), this term decays exponentially for $x>y_1(t)$ and $x<y_2(t)$ but contains a tail for $y_2(t)<x<y_1(t)$. Note that this tail also appears in the integrable case i.e. for $\lambda=0$, and thus it is not related to the lack of integrability.
\end{remark}

Now, we consider the region $t\geq T$.
By \eqref{eq:eH}, 
$\TSR > \frac {10  } { \sqrt{\alpha}}\rho^{-1} \YYzz e^{\frac 12 \YYzz} > 10 (\rho\mu_0)^{-1} |\ln \mu_0|$. 
From \eqref{eq:J11}, \eqref{eq:U+-}, \eqref{eq:U++} and \eqref{eq:ap2} written at $t=T$,  
$$
\left\| U(T) - Q_{ \mu_1(T)} (.-y_1(T)) -Q_{ \mu_2(T)}(.-y_2(T))\right\| \leq C Y_0^\sigma e^{-\frac 54 Y_0}
\leq C' Y_0^\sigma e^{-\frac 34 Y_0} \mu_0 ,
$$
where $|\mu_1(T)-\mu_0|+|\mu_2(T)+\mu_0|\leq C Y_0^2 e^{-Y_0}$.

Therefore, we can apply Proposition \ref{PR:STAB} backwards (i.e. for  $t\geq T$ -- see Remark \ref{RE:STAB}), with $\omega=C' Y_0^\sigma e^{-\frac 34 Y_0}$. There exist
$y_1(t)$, $y_2(t)$ and $\displaystyle\mu_1^+=\lim_{+\infty} \mu_1$, $\displaystyle\mu_2^+=\lim_{+\infty} \mu_2$, such that
$$
w(t)= U(t) - Q_{\mu_1(T)} (.-y_1(t)) -Q_{\mu_2(T)}(.-y_2(t))
$$
satisfies
\begin{equation}\label{eq:final}\begin{split}
	\sup_{t\in [T,+\infty)} \|w(t)\|_{H^1}\leq C \YYzz^{\sigma} e^{-\frac 54 \YYzz},\quad  \lim_{t\to +\infty} \|w(t)\|_{H^1(x>-(99/100)t)}=0, 
	\\ |\mu_1^+-\mu_0|+ |\mu_2^+ +\mu_0|\leq  C Y_0^2 e^{-Y_0},
	\quad 
\lim_{+\infty} \dot y_j=\mu_j^+\quad (j=1,2).
\end{split}\end{equation}

Finally, using the conservation laws and the above asymptotics for $w(t)$, we claim the following refined estimates on the limiting scaling parameters:
\begin{equation}\label{eq:am}  	0 \leq   {\mu_1^+}  - \mu_0 \leq  C \YYzz^{2\sigma} e^{-2 \YYzz},\quad
0 \leq   -{\mu_2^+}  - \mu_0 \leq  C \YYzz^{2\sigma} e^{-2 \YYzz},
\end{equation}
which is a consequence of \eqref{eq:final} and the following lemma.

\begin{lemma}[Monotonicity of the speeds by conservation laws]\label{le:am}
	There exists $C>0$ such that
	\begin{equation*}
		\begin{split}
		&	\frac 1 C e^{\YYzz}\limsup_{t\to +\infty} \|w(t)\|_{H^1}^2 \leq \frac {\mu_1^+}{\mu_0} - 1 \leq  C e^{\YYzz}\liminf_{t\to +\infty} \|w(t)\|_{H^1}^2\leq C \YYzz^{2\sigma} e^{-\frac 32 \YYzz}, \\
		&	\frac 1 C e^{\YYzz}\limsup_{t\to +\infty} \|w(t)\|_{H^1}^2 \leq  - \frac {\mu_2^+}{\mu_0} - 1 \leq  Ce^{\YYzz} \liminf_{t\to +\infty} \|w(t)\|_{H^1}^2\leq C \YYzz^{2\sigma} e^{-\frac 32 \YYzz}.
		\end{split}
	\end{equation*}
\end{lemma}

\begin{proof}
We first write the conservation of mass and energy for $U(t)$ (see \eqref{eq:i5} and \eqref{eq:i6}) and then pass to the limit
$t\to -\infty$, $t\to +\infty$, using \eqref{eq:final}. It follows that the limits
$\lim_{+\infty} M(w)$ and $\lim_{+\infty} \mathcal{E}(w)$ exist and
\begin{align}
&	M(Q_{\mu_0})+M(Q_{-\mu_0}) = M(Q_{\mu_1^+})+M(Q_{\mu_2^+}) + \lim_{+\infty} M(w),	 \label{eq:am1}\\
&	\mathcal{E}(Q_{\mu_0})+\mathcal{E}(Q_{-\mu_0}) = \mathcal{E}(Q_{\mu_1^+})+\mathcal{E}(Q_{\mu_2^+}) + \lim_{+\infty} \mathcal{E}(w).	 \label{eq:am2}
\end{align}
Let
\begin{equation*}
	\nu_1 = \frac {\mathcal{E}(Q_{\mu_0})-\mathcal{E}(Q_{\mu_1^+})}{M(Q_{\mu_1^+})-M(Q_{\mu_0})},\quad
	\nu_2 = \frac {\mathcal{E}(Q_{-\mu_0})-\mathcal{E}(Q_{\mu_2^+})}{M(Q_{\mu_2^+})-M(Q_{-\mu_0})}.
\end{equation*}
so that by \eqref{as:sc} and \eqref{eq:am},
\begin{equation*}
	\left| \frac {\nu_1}{\mu_0}- 1 \right|\leq \frac 14, \quad 
		\left| \frac {\nu_2}{-\mu_0}- 1 \right|\leq \frac 14.
\end{equation*}
We combine \eqref{eq:am1} and \eqref{eq:am2} to get
\begin{align*}
	\lim_{+\infty} \mathcal{E}(w) 	& = \nu_1 (M(Q_{\mu_1^+})-M(Q_{\mu_0})) + \nu_2 (M(Q_{\mu_2^+})-M(Q_{-\mu_0})), 	\\
    		& = (\nu_1 - \nu_2)  (M(Q_{\mu_1^+})-M(Q_{\mu_0})) - \nu_2 \lim_{+\infty} M(w),	\\
		& = (\nu_1 - \nu_2)  (M(Q_{-\mu_0})-M(Q_{\mu_2^+})) - \nu_1 \lim_{+\infty} M(w).
\end{align*}
Since $\|w\|_{L^\infty}\leq C \|w\|_{H^1}\leq C \YYzz^{\sigma} e^{-\frac 54 \YYzz}$, we have
$$\frac 1 2 \limsup_{+\infty} \|w\|_{H^1}^2 \leq \lim_{+\infty} \mathcal{E}(w) \leq 2 \liminf_{+\infty} \|w\|_{H^1}^2,$$
and by $|\nu_1|+|\nu_2|\leq \frac 14$, for $\mu_0$ small enough, we obtain
$\lim_{+\infty} \mathcal{E}(w) + \nu_2 \lim_{+\infty} M(w) > \frac 14 \limsup_{+\infty} \|w\|_{H^1}^2$, 
$\lim_{+\infty} \mathcal{E}(w) + \nu_1 \lim_{+\infty} M(w) < 2 \limsup_{+\infty} \|w\|_{H^1}^2$.

By ${\frac d{d\mu} Q_{\mu}}_{| \mu=0} >0$, and so for all $|\mu|\leq 2 \mu_0$, ${\frac d{d\mu'} Q_{\mu'}}_{| \mu'=\mu} >C>0$,
we get \eqref{eq:am1} and \eqref{eq:am2}.
\end{proof}

For future reference, we   observe that the following hold for $t\in [-T,T]$
\begin{equation}\label{eq:J2} 
  	\sup_{x\in \RR} \{1+e^{\frac 12(x-\yy_1(t))} | E(t,x)|\} \leq C \YYzz^\sigma e^{-\YYzz} e^{-\YYrr(t)}, 
	\quad	\|E(t)\|_{L^2} \leq C \YYzz e^{-\frac 34 \YYzz} e^{-\YYrr(t)},
\end{equation}
\begin{equation}\label{eq:J3} 
	|\dot \mu_j - {\cal M}_j| \leq C \YYzz^{\sigma} e^{-2 \YYzz},\quad 
	|\mu_j - \dot y_j - {\cal N}_j| \leq 	C \YYzz^{\sigma} e^{-\frac 54 \YYzz}.
\end{equation}

\subsection{Proof of Theorem \ref{TH:2}}

First, we claim a refined stability result around the family of asymptotic $2$-soliton solutions (defined in the next claim) in the spirit of Proposition \ref{PR:STAB} but without the exponential error term (see \eqref{eq:stab+}--\eqref{eq:stab-}). The proof is given in Appendix B.

\begin{proposition}[Sharp stability]\label{PR:SHARP}
Let $U$ be defined as in Theorem~\ref{TH:1}.
For $\mu_0>0$  small enough, if $u(t)$ is a solution of {\eqref{eq:BBM}} such that
\begin{equation}\label{eq:sh1}
	\|u({{\mathcal T}_1})-U({{\mathcal T}_1})\|_{H^1} = \omega \mu_0
\end{equation}
for some ${\mathcal T}_1$, where $0<\omega<|\ln \mu_0|^{-2}$,
then there exist $t\in \RR \mapsto (T(t), X(t))\in \RR^2$ such that
\begin{equation}\label{eq:sh2}
	\forall t\in \RR ,\quad 
 \|u(t+T(t),x+X(t))-U(t)\|_{H^1} + |\dot X(t)|+ e^{-\frac 12 \YYzz}|\dot T(t)|\leq C \omega \mu_0 .
\end{equation}
\end{proposition}

Now, we prove Theorem \ref{TH:2}.
Let $\tilde \mu_0\in \RR$ and $\tilde Y_0>0$ be such that 
$$\mu_0 = \sqrt{\tilde \mu_0^2 +  4 \alpha e^{-\tilde Y_0}}$$ is small enough. Let
$u_0\in H^1$ satisfy \eqref{eq:th2b} and let $u(t)$ be the corresponding solution of \eqref{eq:BBM}. We assume that $\tilde \mu_0\leq 0$, the proof being the same in the case $\tilde \mu_0>0$
by using the transformation $x\to -x$, $t\to -t$ and translation in space invariance.

For this value of $\mu_0$, 
let  $U(t)$ and $Y(t)$ be defined as in Theorem \ref{TH:1} and Sections 4.1 and 4.2.
Recall that for all $t$, $\dot Y^2(t) + 4 \alpha e^{-Y(t)} = 4 \mu_0^2$.
Since ${\tilde Y_0} \geq Y_0$, there exists $\tilde T_0<0$ such that 
$Y(\tilde T_0)=\tilde Y_0$, so that $\dot Y(\tilde T_0)=2 \tilde \mu_0$. We claim that for some
$\mathcal{X}_1\in \RR$,
\begin{equation}\label{eq:tbp}
	\|U(\tilde T_0,.+\mathcal{X}_1) - Q_{-\tilde \mu_0}(. - \tfrac {\tilde Y_0} 2) - Q_{\tilde \mu_0}(.+\tfrac {\tilde Y_0} 2)\|_{H^1} \leq C |\ln \mu_0|^{\sigma+2} \mu_0^{3/2}.
\end{equation}
Indeed, if $\tilde T_0 <-T$, thenwe use Corollary \ref{COR:pur}. Otherwise, by Proposition \eqref{pr:st}, we have
\begin{equation*}
	\begin{split}
		& | y(\tilde T_0) - Y(\tilde T_0) | \leq C |\ln \mu_0|^{\sigma+2} \mu_0^{3/2},
		\quad \|\varepsilon(\tilde T_0)\|_{H^1} \leq C |\ln \mu_0|^{\sigma} \mu_0^{5/2},	\\
		& | \mu_1(\tilde T_0) + \tfrac 12 \dot Y(\tilde T_0) | \leq C |\ln \mu_0|^{2} \mu_0^{2}, \quad
		| \mu_2(\tilde T_0) - \tfrac 12 \dot Y(\tilde T_0) | \leq C |\ln \mu_0|^{2} \mu_0^{2}.
	\end{split}
\end{equation*}
Using in addition \eqref{eq:ap2}, we get \eqref{eq:tbp}.

By \eqref{eq:tbp} and \eqref{eq:th2b}, we obtain
$$
\|u_0  -  U(\tilde T_0,.+\mathcal{X}_1)\|_{H^1} \leq C \omega \mu_0 + C |\ln \mu_0|^{\sigma+2} \mu_0^{3/2},
$$
and by Proposition \ref{PR:SHARP}, we obtain \eqref{eq:th2c}.

\section{Nonexistence of a pure $2$-soliton and interaction defect}\label{sec:5}

In this section, we complete the proof of Theorem \ref{TH:1} by proving the lower bounds in \eqref{eq:th8} and  \eqref{eq:th9}.

\subsection{Refined control of the translation parameters}

Now, we introduce specific functionals $\mathcal{J}_j(t)$ related to 
the translation parameters $y_j(t)$ to obtain a refined version of  the dynamical system.

\begin{lemma}\label{PR:J}
Under the assumptions of Proposition \ref{pr:st}, for $j=1,2$, let
\begin{equation}\label{eq:J} 
  {\cal J}_j(t)= \frac 1{\int [(1-\lambda \partial_x^2) \Lambda Q] Q} \int \varepsilon(t,x)(1-\lambda \partial_x^2) \iscal_j(t,x) dx,  \end{equation}
   where $\iscal_j(t,x)=  \int_{-\infty}^x \Lambda \qtud(t,y) dy.$ 
Then $\JJ_j(t)$ is well-defined and the following hold
\begin{enumerate}
	\item[\rm (i)]  Estimates on $\JJ_j$.
	\begin{equation}\label{eq:JJ} \forall t\in [-\TSR,\TSR],\quad 
		|\JJ_1(t)|+|\JJ_2(t)|\leq C \YYzz^{\sigma+1} e^{-\frac 54 \YYzz}.
	\end{equation}
	\item[\rm (ii)]  Equation of $\JJ_j$. For $j=1,2$,
	\begin{equation}\label{eq:JJJ} \forall t\in [-\TSR,\TSR],\quad 
		\left| \frac d{dt}{\cal J}_j(t) - (\mu_j - \dot \yy_j - {\cal N}_j)  \right| \leq  C \YYzz^{\sigma+2} e^{-\frac 74 \YYzz}.
	\end{equation}
\end{enumerate}
\end{lemma}

\begin{remark}\label{rk:J}
The constant $\int ((1-\lambda \partial_x^2) \Lambda Q) Q$ is not zero (see \eqref{eq:58}).

Note also that   $\int \Lambda \neq 0$ (see \eqref{eq:qint3}), and so the functions $\iscal_j(x)$ are bounded but have no decay at $+\infty$ in space.
Therefore, $\JJ_j(t)$ is not well-defined for a general $\varepsilon\in H^1$.
Part of the proof of Lemma \ref{PR:J} consists on obtaining decay in space for $\varepsilon(t)$ in order to give a rigorous sense to $\mathcal{J}_j$.
\end{remark}

\begin{remark}
	Estimate \eqref{eq:JJJ} says formally that $\mu_j - \dot y_j - {\cal N}_j$ is of order $O_{7/4}$, which is 
	an decisive improvement with respect to \eqref{eq:J3} (gain of a factor $e^{-\frac 12 \YYzz}$).
\end{remark}

\begin{proof}
\emph{Preliminary estimates.}
We work under the assumptions of Proposition \ref{pr:st}, and on the interval $[-\TSR,\TSR]$.
First,   we claim exponential decay properties of $U(t)$ on the right ($x>y_1(t)$).

\begin{claim}[Decay estimate on $u(t)$]\label{cl:J1}
	There exist $C>1$ and $\rho_0 >1$ such that for all $t\in [-\TSR,\TSR]$, for all $X_0>1$,  
	\begin{equation}\label{eq:J4}
		\int_{x>X_0+y_1(t)} (U^2+(\partial_x U)^2) (t,x)  dx\leq C e^{- \rho_0 {X_0}}.
	\end{equation}
\end{claim}

Recall that the proof of Claim \ref{cl:J1} is  obtained by  
$$
	\lim_{t\to -\infty} \|U(t)\|_{H^1(x>\frac 12 |t|)}=0
$$
combined with     monotonicity arguments, see e.g. \cite{DI2} for the case of the \eqref{eq:BBM} equation.

\medskip

\emph{Estimate of $\JJ_j$.}
Note that $\iscal_j$ does not belong to $L^2$ (see Remark \ref{rk:J}) but   satisfies
\begin{equation}\label{eq:J5}
	\sup_{x\in \RR} \left\{\left(1+e^{- \frac 12 (x-\yy_j(t))} \right) | \iscal_j(t,x)|\right\} \leq C.
\end{equation}
It follows from \eqref{eq:J4}, \eqref{eq:J5}, and the decomposition of $U(t)$ in Lemma \ref{PR:de} that
\begin{equation*}\begin{split}
	|\JJ_1(t)| & \leq   C \int_{x<y_1(t)} |\varepsilon(t,x)| |(1-\lambda \partial_x^2) \iscal_j(t,x)| dx+  C \int_{y_1(t)\leq x<y_1(t)+10\rho_0^{-1} \YYzz} |\varepsilon(t,x)|  dx\\ &
	+  C \int_{x>y_1(t)+10\rho_0^{-1} \YYzz} |\varepsilon(t,x)| |(1-\lambda \partial_x^2) \iscal_j(t,x)| dx \\
	& \leq C (1+ \YYzz) \|\varepsilon(t)\|_{L^\infty}+ C \int_{x>y_1(t)+10\rho_0^{-1} \YYzz} |U(t,x)| + C e^{- 5  \YYzz }
	\\ & \leq C \YYzz^{\sigma+1} e^{-\frac 54 \YYzz}.
\end{split}\end{equation*}

Moreover, using  $\yy_1(t)-\yy_2(t)=\pyy(t)\leq   \YYrr(T) \leq C \YYzz$, one gets  by similar arguments $|\JJ_2(t)|\leq C \YYzz^{\sigma +1} e^{-\frac 54 \YYzz}$.

\medskip

\emph{Equation of $\JJ_1$.}
To prove \eqref{eq:JJJ}, we make use of the equation of $\varepsilon$ (see \eqref{eq:ep}),
and of the special algebraic structure of the approximate solution $V(t,x)$ introduced in Propositions \ref{PR:ap} and \ref{PR:app}.
We have 
$$\left(\int [(1-\lambda \partial_x^2) \Lambda Q] Q\right) \frac d{dt} \JJ_1(t)=\int (1-\lambda \partial_x^2)\partial_t \varepsilon \iscal_1 + \int \varepsilon \partial_t (1-\lambda \partial_x^2) \iscal_1.$$

First observe that
$$
\partial_t  \iscal_1  (x)=  \int_{-\infty}^x \partial_t \Lambda \qtun(y) dy
= \int_{-\infty}^x \left\{ \dot \mu_1 \frac {\partial \Lambda \qtun} {\partial \mu_1} 
+ \dot \yy_1 \frac {\partial \Lambda \qtun} {\partial \yy_1}\right\}(y) dy.
$$
Thus, by  $|{\cal M}_j|+|{\cal N}_j|\leq C e^{-\YYzz}$, $|\mu_j(t)|\leq C e^{-\frac 12 \YYzz}$, \eqref{eq:J3} and  \eqref{eq:J4}, we have
\begin{equation}\label{eq:J6}
	\left|\int \varepsilon \partial_t  (1-\lambda \partial_x^2)\iscal_j   \right|\leq C \YYzz^\sigma e^{-\frac 74 \YYzz}.
\end{equation}

Next, using \eqref{eq:ep} and $\partial_x \iscal_1= \Lambda \qtun$, we have
\begin{equation*}
\int (1-\lambda \partial_x^2) \partial_t \varepsilon \iscal_1 = \int (\partial_x^2 \varepsilon - \varepsilon + 2 V \varepsilon + \varepsilon^2 ) \Lambda \qtun
-\int E \iscal_1  +\int \fell(V) \iscal_1		.
\end{equation*}

For the term $\int (\partial_x^2 \varepsilon - \varepsilon + 2 V \varepsilon + \varepsilon^2 ) \Lambda \qtun$, we argue as the proof of Lemma \ref{PR:de}.
Using $L_{\mu_j} \Lambda Q_{\mu_j} = (1-\lambda \partial_x^2) Q_{\mu_j}$ (see \eqref{eq:Qm}), \eqref{eq:ap2bis}, \eqref{eq:sQQ}, $\int \varepsilon (1-\lambda \partial_x^2)  \qtud=0$, 
Proposition \ref{pr:st} and the definition of $V$ (see Proposition \ref{PR:app}), we obtain
\begin{equation*}
	\left| \int (\partial_x^2 \varepsilon - \varepsilon + 2 V \varepsilon + \varepsilon^2 ) \Lambda \qtun \right|
	\leq C \YYzz^2 e^{-\YYzz} \|\varepsilon\|_{L^2}+ C \|\varepsilon\|_{L^2}^2
	\leq C \YYzz^{\sigma+2} e^{-\frac  94\YYzz} .
\end{equation*}

By \eqref{eq:J2} and \eqref{eq:J5}, we have
\begin{equation*}
	\left| \int E \iscal_1 \right| \leq C \YYzz^\sigma e^{- 2 \YYzz} .
\end{equation*}

Next, we consider the term $\int \fell(V) \iscal_1$. From the definition of $\fell(V)$ in \eqref{eq:defell}, the structure of
$V_0$ and $V$, see \eqref{eq:ap} and \eqref{eq:V} (see also \eqref{eq:de1}), and \eqref{eq:J3}, we have
\begin{equation}\label{eq:J7}
	\sup_{x\in \RR}\Bigg\{(1+e^{\frac 12(x-y_1(t))})\bigg| \fell(V) -\sum_{j=1,2} (\mu_j-\dot y_j - {\cal N}_j) (1-\lambda \partial_x^2 ) \partial_x \qtud \bigg|\Bigg\}
	\leq C \YYzz^{\sigma} e^{- 2 \YYzz}.
\end{equation}
Thus, by \eqref{eq:J5}, we obtain
\begin{equation}\label{eq:J8}
	\left| \int \fell(V) \iscal_1 -  (\mu_1-\dot y_1 - {\cal N}_1) \int [(1-\lambda \partial_x^2 ) \partial_x \qtun] \iscal_1 \right|
	\leq C \YYzz^\sigma e^{- 2 \YYzz} .
\end{equation}

Finally, using $\left|\int [(1-\lambda \partial_x^2 ) \partial_x \qtun] \iscal_1 + \int [(1-\lambda \partial_x^2 )  Q] \Lambda Q \right|
\leq C |\mu_1(t)|\leq C e^{-\frac 12 \YYzz} $ and \eqref{eq:J3}, we obtain \eqref{eq:JJJ}.

The proof for $\frac d{dt} \JJ_2$ is exactly the same.
\end{proof}

\subsection{Preliminary symmetry arguments}

First, we claim the following additional information obtained on the parameters of the solution $U(t)$, under the assumptions of Proposition \ref{pr:st}.

\begin{claim} For all $t\in [-T,T],$
\begin{equation}\label{eq:J9}
	|\mu_1(t)-\mu_2(-t)|\leq C \YYzz^{\sigma+3} e^{-\frac  54 \YYzz}.
\end{equation}
\end{claim}

\begin{proof}[Proof of \eqref{eq:J9}]
From \eqref{eq:U+-} and $\dot \YYrr(0)=0$, we have
$|\mu_1(0)- \mu_2(0)|\leq C \YYzz^{\sigma+1} e^{-\frac 54 \YYzz}$.
From \eqref{eq:J11}, $|\pyy(t)-\pyy(-t)|\leq |\pyy(t)-\YYrr(t)|+|\YYrr(t)-\pyy(-t)|\leq  C \YYzz^{\sigma+2} e^{-\frac 34 \YYzz}$.
Thus, by \eqref{eq:J3} and the expression of ${\cal M}_j$ in \eqref{eq:pa}, we obtain
\begin{equation*}\begin{split}
&	 \left| \dot \mu_1(t) - \left\{\alpha\,  e^{-\pyy(t)} +  \beta\,  \mu_1(t) \pyy(t) e^{-\pyy(t)} + \delta  \,\mu_1(t) e^{-\pyy(t)}\right\}\right|\leq
C \YYzz^{\sigma} e^{- 2 \YYzz},	\\
&	\left| -\dot \mu_2(-t) - \left\{\alpha\,  e^{-\pyy(-t)} +  \beta\,  \mu_2(-t) \pyy(-t) e^{-\pyy(-t)} + \delta  \,\mu_2(-t) e^{-\pyy(-t)}\right\}\right|\leq
C \YYzz^{\sigma} e^{- 2 \YYzz}, \\
\text{and so}\quad
&	\left| -\dot \mu_2(-t) - \left\{\alpha\,  e^{-\pyy(t)} + \beta\,  \mu_2(-t) \pyy(t) e^{-\pyy(t)} + \delta  \,\mu_2(-t) e^{-\pyy(t)}\right\}\right|\leq
C \YYzz^{\sigma+2 } e^{- \frac 7 4  \YYzz}.
\end{split}\end{equation*}
It follows that for all $t\in [-\TSR,\TSR]$, $|\mu_1(t)-\mu_2(-t)|\leq C \YYzz^{\sigma+3} e^{-\frac  54 \YYzz}$.
\end{proof}

The next lemma claims that if the asymptotic $2$-soliton solution $U(t)$ considered in Proposition~\ref{pr:st} has an approximate symmetry property (i.e. $U(t,x) -  U(-t + {t_0}, -x +  x_0) $ is small for some ${t_0}$, $x_0$) then the corresponding decomposition parameters (i.e. $\Gamma(t)$ in Proposition \ref{pr:st}) also have some symmetry properties, despite the fact that the decomposition itself is not symmetric (see the definition of $V(t,x)$ in Propositions \ref{PR:ap}--\ref{PR:app}).

\begin{lemma}\label{le:id}
	Let $t_0,$ $x_0$ be such that $|t_0|\leq 1$, $|x_0|\leq 1$.
	Under the assumptions of Proposition \ref{pr:st}, for all $t\in [-\TSR,\TSR]$,
	\begin{equation}\label{eq:idd}\begin{split}
		 & |\mu_1(t)-\mu_2(-t+t_0)|+ |y_1(t)+y_2(-t+t_0)-x_0| \\ 
		 & \leq C (\| U(t,x) -  U(-t + t_0, -x +  x_0) \|_{H^1} + \YYzz^5e^{-2\YYzz}).
	\end{split}\end{equation}
	In particular, assume that $U(t,x) =  U(-t + t_0, -x +  x_0)$ for some $t_0,$ $x_0$, then
	\begin{equation}\label{eq:iddrrr} 
	  |\mu_1(t)-\mu_2(-t+t_0)|+ |y_1(t)+y_2(-t+t_0)-x_0|  \leq C \YYzz^5e^{-2\YYzz}.
	 \end{equation}
\end{lemma}

\begin{remark}
	Assuming $\| U(t,x) -  U(-t + t_0, -x +  x_0) \|_{H^1}\leq C \YYzz^{\sigma+3} e^{-\frac 54 \YYzz},$ it follows from \eqref{eq:idd} that the following hold
	\begin{equation}\label{eq:iddd}
		|x_0|\leq C \YYzz^{2} e^{-\frac 12 \YYzz }\quad |t_0|\leq C \YYzz^{\sigma+3} e^{-\frac 14 \YYzz}.
	\end{equation}
	Indeed,
	on the one hand, estimate $|x_0|\leq C \YYzz^2 e^{-\frac 12 \YYzz}$ follows from \eqref{eq:J11} taken at time $t=t_0/2$ and \eqref{eq:idd} taken at $t=t_0/2$.

	On the other hand, from \eqref{eq:J9} and \eqref{eq:iddd}, we have $|\mu_1(t)-\mu_1(t-t_0)|\leq C \YYzz^{\sigma+3} e^{-\frac 54 \YYzz}$. Since
	$\dot \mu_1(t)\geq C e^{-\YYzz}$ for $|t|$ close to $0$, we obtain $|t_0|\leq C \YYzz^{\sigma+3} e^{-\frac 14 \YYzz}$.
\end{remark}

Using Section 2, the proof is similar to the one of Lemma 5.2 in \cite{MMkdv4} and it is omitted.

\subsection{Lower bound on the defect}\label{sec:5.5}

In this section, we prove the following result.

\begin{proposition}\label{PR:qualit}
	Let $\lambda \in (0,1)$.
	Under the assumptions of Proposition \ref{pr:st}, for a possibly smaller $\mu_0>0$, there exists a constant $c>0$ such that,   
	\begin{equation}\label{eq:lower}\begin{split}
		& \liminf_{t\to +\infty} \|w(t)\|_{H^1} \geq  c \YYzz e^{-\frac 32 \YYzz} ,
		\\ &
		\mu_1^+ - \mu_0 \geq c \YYzz^2 e^{-\frac 52 \YYzz} ,\quad
		-\mu_2^+ - \mu_0\geq  c \YYzz^2  e^{-\frac 52 \YYzz}.
	\end{split}\end{equation}
\end{proposition}

\begin{proof}
The proof is the same as the one of Proposition 5.2 in \cite{MMkdv4} but we repeat it here since
the argument is the key of the nonexistence of a pure $2$-soliton solution.

It suffices to prove the estimate on $w(t)$. The estimates on the final parameters then follow from Lemma \ref{le:am}.

Let $\epsilon>0$ arbitrary, and suppose for the sake of contradiction that 
\begin{equation}\label{eq:def2}
	\liminf_{t\to +\infty} \| w(t)\|_{H^1} \leq   \epsilon \YYzz e^{ - \frac 32 \YYzz}.
\end{equation}
 
\emph{Step 1.}
We claim that for some $\tilde T(t)$, $\tilde X(t)$,  for all $t\in \RR$,
 \begin{equation}\label{eq:def5}
	\|U(-t+ \tilde T(t),-x + \tilde X(t))-U(t,x)\|_{H^1} 
	+ |\dot {\tilde X}(t)| + e^{-\frac 12 \YYzz} |\dot {\tilde T}(t)| \leq C  \epsilon  \YYzz e^{-\frac 32 \YYzz}.
\end{equation}
\emph{Proof of \eqref{eq:def5}.}
By Lemma \ref{le:am}, it follows that
$$
0\leq \mu_1^+ -\mu_0 \leq C \epsilon^2 \YYzz^2 e^{- \frac 52 \YYzz},\quad
0\leq -(\mu_2^+ + \mu_0) \leq C \epsilon^2  \YYzz^2 e^{- \frac 52 \YYzz}.
$$
In particular, for all $t$
$$
\| Q_{\mu_1^+}(.-y_1(t)) + Q_{\mu_2^+}(.-y_2(t))
- (Q_{\mu_0}(.-y_1(t)) + Q_{-\mu_0}(.-y_2(t)))\|_{H^1}\leq C \epsilon^2 \YYzz e^{- \frac 52 \YYzz}.
$$ 
From \eqref{eq:def2} and the behavior of $U$, it follows that  there exist $\mathcal{T}_1$, $\mathcal{T}_2>T$ and $X$
 such that
\begin{equation}\label{eq:def3}
	\|U(\mathcal{T}_1,x)-U(-\mathcal{T}_2,-x+X)\|_{H^1} \leq  
	2 \epsilon Y_0e^{ - \frac 32 \YYzz} +C \epsilon^2 Y_0^2 e^{-\frac 52 \YYzz}
	\leq 3  \epsilon \YYzz e^{ - \frac 32 \YYzz},
\end{equation}
for $\YYzz$ large enough.
From  Proposition \ref{PR:SHARP}, it follows that there exist $\tilde T(t)$ and $\tilde X(t)$ such that
 \eqref{eq:def5} holds.

\medskip

\emph{Step 2.}  {Conclusion of the proof of Proposition \ref{PR:qualit}.}
Take $0<t_1<t_2$ such that $\YYrr(t_1) = \YYzz+1$ and
$\YYrr(t_2)=\YYzz+2$. Note that $t_2-t_1< C e^{\frac 12 \YYzz}$ since $\dot \YYrr > c_0 e^{-\frac 12 \YYzz}$ on $[t_1,t_2]$, for some $c_0>0$.

Note that for $t\in [-T,T]$, $\tilde T(t)$ and $\tilde X(t)$ are small by Proposition \ref{pr:st}.
Applying Lemma~\ref{le:id} at $t_1$ and  $t$, for all $t\in [t_1,t_2]$, we obtain (for $\YYzz$ large enough depending on $\epsilon$)
\begin{equation*} |\mu_1(t)-\mu_2(-t+\tilde T(t))|+ | y_1(t)+y_2(-t+\tilde T(t)) - \tilde X(t)|\leq C  \epsilon Y_0 e^{-\frac 32 \YYzz}
	 .
\end{equation*}
By \eqref{eq:def5}, for all $t\in [t_1,t_2]$, we have $| \tilde T(t)-\tilde  T(t_1)|\leq C \epsilon Y_0 e^{-\frac 12 \YYzz}$, $| \tilde  X(t)- \tilde  X(t_1)|\leq
C \epsilon Y_0 e^{- \YYzz}$ and thus,
\begin{equation*}\begin{split}
 \forall t\in [t_1,t_2],\quad
	& | \mu_2(-t+ \tilde T(t_1))- \mu_2(-t+ \tilde T(t))|\leq C \epsilon e^{-\frac 32 \YYzz} ,\\
	& | y_2(-t+ \tilde T(t_1)) - \tilde X(t_1) - ( y_2(-t+ \tilde T(t))- \tilde X(t)) |\leq C \epsilon Y_0 e^{-  \YYzz} 
\end{split}\end{equation*}
Therefore, setting 
$$
	\nu(t) = \mu_1(t) - \mu_2(-t+ \tilde T(t_1)), \quad
	z(t) = y_1(t) + y_2(-t+ \tilde T(t_1)) - \tilde X(t_1),
$$
we obtain
\begin{equation}\label{eq:def6}
	|\nu(t)|\leq C \epsilon \YYzz e^{-\frac 32 \YYzz}, \quad |z(t)|\leq C \epsilon Y_0 e^{-  \YYzz}.
\end{equation}

We claim
\begin{align} 
	& 	|\JJ_1(t)|+|\JJ_2(t)|\leq C \YYzz^{\sigma+1} e^{-\frac 54 \YYzz},		\label{eq:stp3bis}
\\	& \left| \frac d{dt} \left( \pzz(t) -  \{(b_- \pyy(t)+(b_-+d_-)\} e^{-\pyy(t)} \right) +   (\dot \JJ_1(t) - \dot \JJ_2(-t+\tilde T(t_1)))\right|
	\leq C \epsilon \YYzz e^{-\frac 32 \YYzz},	\label{eq:stp3}
\end{align}
where
$$
b_- = -\frac 12 (b_1-b_2),\quad d_- = -\frac 12 (d_1-d_2).
$$
Note that estimate \eqref{eq:stp3bis} is just \eqref{eq:JJ} from Lemma \ref{PR:J}. Assuming estimate \eqref{eq:stp3} for the moment, we complete the proof of the Proposition.

Integrating  \eqref{eq:stp3} on $[t_1,t_2]$, using \eqref{eq:stp3bis} and $t_2-t_1< C e^{\frac 12 \YYzz}$, we obtain
\begin{equation}\label{eq:c16}\begin{split}
&	\left| (\pzz(t_1)- \{b_- \pyy(t_1) + (b_-+d_-)\} e^{-\pyy(t_1)}) - (\pzz(t_2)- \{b_- \pyy(t_2) + (b_-+d_-)\} e^{-\pyy(t_2)}) \right| 
\\ & \leq 	C \epsilon \YYzz  e^{-\YYzz}.
\end{split}\end{equation}
Thus, by \eqref{eq:def6}, for $k= 1 + \frac {d_-}{b_-}$, ($b_-\neq 0$ for $\lambda \neq 0$),
\begin{equation}\label{eq:c17}
	|(\pyy(t_1)+k) e^{-\pyy(t_1)} - (\pyy(t_2)+k) e^{-\pyy(t_2)}|\leq C \epsilon \YYzz  e^{-\YYzz}.
\end{equation}
But since $\YYrr(t_1) = \YYzz+1$ and
$\YYrr(t_2)=\YYzz+2$,  \eqref{eq:c17} is a contradiction for $\epsilon$ small enough and $\YYzz$ large enough.

\medskip

Let us now prove \eqref{eq:stp3}. By \eqref{eq:JJJ} and the expression of $\mathcal{N}_j$ in \eqref{eq:pa},
we have
\begin{equation}\label{eq:c3}\begin{split}
	&	\dot y_1 = \mu_1 - a e^{-\pyy} - b_1 \mu_1 \pyy e^{-\pyy} - d_1 \mu_1 e^{-\pyy} - \dot \JJ_1 +O(\YYzz^{\sigma+2} e^{-\frac 74 \YYzz}) ,	\\
	&	\dot y_2 = \mu_2 - a e^{-\pyy} - b_2 \mu_2 \pyy e^{-\pyy} - d_2 \mu_2 e^{-\pyy} - \dot \JJ_2 +O(\YYzz^{\sigma+2} e^{-\frac 74 \YYzz}).
\end{split}\end{equation}
Moreover, by \eqref{eq:def6}, we check
\begin{equation}\label{eq:c4}
	\left |e^{-\pyy(t)}-e^{-\pyy(-t+\tilde T(t_1))}\right|\leq C \epsilon\YYzz e^{-2\YYzz}.	
\end{equation}
Thus, using again \eqref{eq:def6}, we obtain
\begin{equation*}\begin{split}
	\dot{\pzz}(t) &= \dot y_1(t) -  \dot y_2(-t+\tilde T(t_1))
	\\ &=  \nu(t) - (b_1-b_2) \mu_1 \pyy e^{-\pyy} - (d_1-d_2) \mu_1 e^{-\pyy} - (\dot \JJ_1(t)- \dot \JJ_2(-t+\tilde T(t_1))) + O(\YYzz^{\sigma+2} e^{-\frac 74 \YYzz}),
\end{split}\end{equation*}
Since $|\mu_1 + \mu_2 |\leq C \YYzz^2 e^{-\YYzz}$ (see \eqref{eq:J11}), we have
$|\mu_1 -\frac 12 \mu|\leq C \YYzz^2 e^{-\YYzz}$ and thus, by $|\mu  - \dot \pyy|\leq C e^{-\YYzz}$, we obtain
$|\mu_1 e^{-\pyy} - \frac 12 \dot \pyy e^{-\pyy}|\leq C \YYzz^2 e^{-2 \YYzz}.$ Therefore,
$$
\dot {\pzz}=  \nu + b_- \dot \pyy \pyy e^{-\pyy} + d_- \dot \pyy e^{-\pyy} -(\dot \JJ_1(t)- \dot \JJ_2(-t+\tilde T(t_1))) + O(\YYzz^{\sigma+2} e^{-\frac 74 \YYzz}),
$$
where  $| \nu(t)|\leq C \epsilon e^{-\frac 32 \YYzz}$ from \eqref{eq:def6}.
Estimate \eqref{eq:stp3} then follows.
\end{proof}

\appendix

\section{Appendix to the construction of an approximate solution}\label{AP:ZZ}

\subsection{Linearized operator, identities and asymptotics for solitons}
Recall that we set
\begin{equation}\label{eq:QMM}
	Q_{\mu}(x) = (1+\mu) Q\left(\sqrt{\frac {1+\mu}{1+\lambda \mu}} \, x\right)
	\quad \text{where}\quad
	Q(x)=\frac 32 \frac 1 {\cosh^2 \left(\frac x2\right)}
\end{equation}
satisfies
\begin{equation}\label{eq:QQ}
	Q''+Q^2=Q\quad  \text{and}\quad   (Q')^2 +\frac 23 Q^3 = Q^2.
	\end{equation}

We   recall the following well-known spectral properties of $L$
(see \cite{We3} and  Lemma 2.2 from \cite{MMcol1})

\begin{claim}[Properties of the operator $L$]\label{LE:A2}
	The operator ${L}$ defined in $L^2(\mathbb{R})$ by
    $$
        {L} f= - f'' + f - 2Q f
    $$
    is self-adjoint and satisfies the following properties:
    \begin{itemize}
        \item[{\rm (i)}] First eigenfunction : ${L} Q^{\frac 3 2} = - \frac 54  Q^{\frac 3 2}$;
        \item[{\rm (ii)}] Second eigenfunction : ${L} Q'=0$; the kernel of ${L}$ is 
        $\{c_1 Q', c_1 \in \mathbb{R}\}$;
        \item[{\rm (iii)}] For any   function $h \in L^2(\mathbb{R})$ orthogonal to $Q'$ for the $L^2$ scalar product, 
        there exists a unique function $f \in H^2(\mathbb{R})$ orthogonal to $(1-\lambda \px^2)Q'$ such that ${L} f=h$; moreover,
        if $h$ is even (respectively, odd), then $f$ is even (respectively, odd).
        \item[{\rm (iv)}] Suppose that $f\in H^2(\mathbb{R})$ is such that ${L} f \in \mathcal{Y}$.
            Then, $f\in \mathcal{Y}$.
	\item[{\rm (v)}]  
	There exists $c_1>0$ such that for all $f\in H^1(\mathbb{R})$,
	$$
	\int (1-\lambda\partial_x^2)Qf = \int (1-\lambda \partial_x^2)Q' f =0 \quad
	\Rightarrow \quad (L f,f ) \geq c_1 \|f\|_{H^1}^2.
	$$
    \end{itemize}
\end{claim}

\begin{claim}[Preliminary computations on solitons]\label{LE:A1}
\begin{itemize}
\item[\rm (i)] Scaling.
$$
(1+\lambda \mu) Q_{\mu}'' 
- (1+\mu) Q_{\mu} + Q_{\mu}^2 =0.
$$
Set
$$
	\Lambda Q_\mu =  \left(\frac d{d\mu'} Q_{\mu '}\right)_{|\mu'=\mu},\quad
	\Lambda^2 Q_\mu= \left(\frac {d^2}{{d\mu '}^2} Q_{\mu'} \right)_{|\mu'=\mu}.
$$
Then,
\begin{equation}\label{eq:ST}\begin{split}
	& \Lambda Q =\Lambda Q_0=Q+\frac 12 (1-\lambda) x Q', 
	\\ & \Lambda^2 Q= \Lambda^2 Q_0=\frac 34 (1-\lambda) (1-\lambda + \lambda^2) x Q' +\frac 14 (1-\lambda)^2 x^2 Q - \frac 14 (1-\lambda)^2 x^2 Q^2.
\end{split}\end{equation}

\item[\rm (ii)] Linearized operator.
Let
\begin{equation}
\label{eq:dL}
	L_\mu v = - (1+\lambda \mu) v'' + (1+\mu) v - 2Q_\mu v , \quad Lv=L_0 v = - v'' + v - 2Qv .
\end{equation}
Then,
\begin{equation}
\label{eq:Qm} 
  L_{\mu} Q_\mu = - Q_\mu^2,  \quad  L_{\mu} \Lambda Q_\mu = - (1-\lambda \partial_x^2) Q_\mu,\quad L_\mu Q_\mu'=0,
  \end{equation}
\begin{equation}\label{eq:qpq}
	L \frac {Q'}Q = - \frac 53 Q' + \frac {Q'}{Q},\quad
\left( L \frac {Q'}Q\right)' = - 2 Q +  \frac 53 Q^2,\quad
\lim_{\pm \infty} \frac {Q'}{Q} = \mp 1.
\end{equation}
\item[\rm (iii)] Integral identities.
\begin{align}
  & \int Q=\int Q^2,\quad 
\int Q^3=\frac{6}{5}\int Q^2, \quad \int Q'^2=\frac15\int Q^2,\label{eq:qint1}\\
& \int Q^2=6,  \quad \int \Lambda Q = \frac 12 (1+\lambda)\int Q=3(1+\lambda),
\label{eq:qint3}
\end{align}
\begin{equation}\label{eq:58}
\int  [(1-\lambda \partial_x^2) (\Lambda Q ) ] Q= \frac 3{10} (15+10 \lambda - \lambda^2),\quad
\int Q\Lambda Q= \frac 32 (3+\lambda),
\end{equation}
\begin{equation}\label{eq:sc2}\begin{split}
	&\int Q_\mu^2 = (1+\mu)^{\frac 32}(1+\lambda\mu)^{\frac 12}\int Q^2,\quad
	\int Q_\mu^3 = (1+\mu)^{\frac 52} (1+\lambda \mu)^{\frac 12} \int Q^3,\\
	&\int (  Q_\mu')^2 = (1+\mu)^{\frac 52} (1+\lambda \mu)^{-\frac 12} \int (Q')^2,
	\end{split}\end{equation}
\begin{equation}\label{as:sc} 
 	  - \frac d{d\mu} \mathcal{E}(Q_\mu) = \mu \frac d{d\mu} M(Q_\mu),\quad
	\frac d{d\mu} M(Q_\mu) >0 .
\end{equation}
\item[\rm (iv)] Pointwise identities.
\begin{equation}\label{eq:56}
Q -Q' = e^{x} (Q+Q') =  \frac {12 e^{2x}} {(e^{x}+1)^3} ,\quad e^{-x} Q^2 = - Q^2 + 3 (Q'+Q),
\end{equation}
\begin{equation}\label{eq:57}  e^{-x}((\Lambda Q)'-\Lambda Q-\frac 12(1-\lambda) Q)=
- \frac 12 (3-\lambda) (Q'+Q) 
+ \frac 12 (1-\lambda) x(Q^2 -2  (Q'+Q)).
\end{equation}
\item[\rm (v)] Asymptotics
\begin{equation}\label{eq:as}
Q(x)=6 e^{-x}-12 e^{-2x} + O(e^{-3x}) \quad \text{at $+\infty$}.
\end{equation}
\end{itemize}
\end{claim}
\begin{proof}
(i) 
First, we check that $Q_\mu(.+x_0)$  solves the following equation
\begin{equation}\label{eq:69bis}
(1+\lambda \mu) Q_{\mu}''(.+x_0)-(1+\mu) Q_{\mu}(.+x_0) + Q_{\mu}^2(.+x_0) = 0.
 \end{equation}
Indeed, we have
\begin{align*}
& (1+\lambda \mu) Q_{\mu}''(x+x_0)   = (1+\mu)^2 Q''\left(\sqrt{\frac {1+\mu}{1+\lambda \mu}}  (x+x_0)\right)
\\ &=(1+\mu)^2 Q\left(\sqrt{\frac {1+\mu}{1+\lambda \mu}} (x+x_0)\right) - (1+\mu)^2 Q^2\left(\sqrt{\frac {1+\mu}{1+\lambda \mu}} (x+x_0)\right) \\ & = (1+\mu) Q_\mu(x+x_0) - Q_\mu^2(x+x_0).
\end{align*}
 
We have by direct computations
$$
\Lambda Q_\mu			=Q\left(\sqrt{\frac {1+\mu}{1+\lambda \mu}} \, x\right) +
			\frac 12 \frac {(1+\mu)^{\frac 12}(1-\lambda)}{(1+\lambda \mu)^{\frac 32}} x Q' \left(\sqrt{\frac {1+\mu}{1+\lambda \mu}} \, x\right).
$$
The expressions of $\Lambda Q$ and $\Lambda^2 Q$ then follow.

\medskip

(ii) Differentiating \eqref{eq:69bis} with respect to $\mu$ and then with respect to $x_0$, we obtain
\begin{equation}\label{as:sc1}\begin{split}
	&	(1+\lambda \mu ) (\Lambda Q_{\mu})''-(1+\mu) \Lambda Q_{\mu} + 2 Q_{\mu} \Lambda Q_{\mu} = -\lambda Q_{\mu}'' + Q_{\mu},	\\
	&	(1+ \lambda \mu) (Q_{\mu}')''-(1+\mu) Q_{\mu}' + 2 Q_{\mu} Q_{\mu}' =0.
\end{split}\end{equation}

Let us check \eqref{eq:qpq}. First, $\lim_{\pm \infty} \frac {Q'}Q = \mp 1$ is clear from the expression of $Q$. 
Next, we have, using \eqref{eq:QQ},
\begin{equation}\label{eq:69}
	\left(\frac {Q'}{Q}\right)' = \frac {Q''Q-(Q')^2}{Q^2} = -\frac 13 Q.
\end{equation}
Thus,
$$
- L \frac {Q'}Q = -\frac 13 Q' - \frac {Q'}{Q} + 2Q' =  \frac 53 Q' - \frac {Q'}{Q},\quad
\left( - L \frac {Q'}Q\right)' = 2 Q - \frac 53 Q^2.
$$

(iii) These identities are readily obtained from \eqref{eq:QQ}. Note for example:
\begin{align*}
& \int  [(1-\lambda \partial_x^2) \Lambda Q ] Q= \int (Q+ \frac 12 (1-\lambda) xQ')
(Q-\lambda Q + \lambda Q^2)\\
&= (1-\lambda) \int Q^2 + \lambda \int Q^3 -\frac 14 (1-\lambda)^2 \int Q^2
-\frac 16 (1-\lambda) \lambda \int Q^3\\
& = (1-\lambda) (1-\frac 14 + \frac \lambda 4) \int Q^2 
+ \frac 65 (\lambda - \frac \lambda 6 + \frac {\lambda^2}6) \int Q^2
=\frac 3{10} (15+10 \lambda - \lambda^2).
\end{align*}

The identities on $\int Q_\mu^2$, $\int (  Q_\mu ')^2$ and $\int Q_{\mu}^3$ follows directly from 
\eqref{eq:QMM}. Now, we prove \eqref{as:sc}. We first observe that 
\begin{align*}
	  \frac 12 \frac d{d\mu} \mathcal{E}(Q_{\mu}) & =  \int ( Q_\mu')( \Lambda Q_\mu)'  +  Q_\mu \Lambda Q_\mu -  Q_\mu^2 \Lambda Q_\mu	\\
&	= - \mu \int [\lambda (  Q_\mu') (  \Lambda Q_\mu)' + Q_\mu \Lambda Q_\mu]  = - \frac 12  \mu \frac d{d \mu} M(Q_{\mu}),
\end{align*}
is a consequence of \eqref{eq:69bis},   multiplied by $\Lambda Q_\mu$ and integrated over $\RR$.

Then, we check $\frac d{d \mu} M(Q_{\mu})>0$. For $\mu$ small, the result is true by \eqref{eq:58} and a perturbation argument.
In fact, it is true for all $\mu>-1$ (see Weinstein \cite{We3}). Indeed,
by the expressions of $\int (\partial_x Q_\mu)^2$ and $\int Q_\mu^2$, we have
$$
M(Q_\mu) = \int \lambda ( Q_\mu')^2 + Q_\mu^2 =  \frac 15 (\int Q^2) (1+\mu)^{\frac 32} (1+\lambda \mu)^{-\frac 12}  (5+\lambda (1+6 \mu)).
$$
Differentiating with respect to $\mu$, we find
\begin{align*}
&	\frac d{d\mu} M(Q_\mu)   = \frac 1{10} (\int Q^2) (1+\mu)^{\frac 12} (1+\lambda \mu)^{-\frac 32}  (15+10\lambda - \lambda^2 + 40 \lambda \mu + 8 \lambda^2 \mu + 24 \lambda^2 \mu^2)\\
& 	= \frac 1{10} (\int Q^2) (1+\mu)^{\frac 12} (1+\lambda \mu)^{-\frac 32}   (15 (1-\lambda)^2 + 8 \lambda (1+\mu) (5 (1-\lambda) + 3 \lambda (1+\mu)))>0.
\end{align*}

(iv)--(v)
These identities and asymptotic properties are easily obtained from the explicit expression of $Q$:
$$Q(x)=  \frac {6e^{x}} {(e^{x}+1)^2}=\frac {6e^{-x}} {(e^{-x}+1)^2},\quad 
Q'(x)=  \frac {6(e^{x}-e^{2x})} {(e^{x}+1)^3}=\frac {6(e^{-2x}-e^{-x})} {(e^{-x}+1)^3}.
$$
In particular, we observe that
\begin{equation*}
	Q^2 =36 \frac {e^{-2x}}{(e^{-x} +1)^4}
	= 3 \frac {Q'+Q}{e^{-x} +1}
\quad 
\text{and so}
\quad
e^{-x} Q^2 = - Q^2 + 3 (Q'+Q).
\end{equation*}
We obtain in particular
\begin{equation}\label{eq:A777}
10 \int e^{-x} Q^3 = 9 \int e^{-x} Q^2.
\end{equation}
Moreover,
\begin{equation*}\begin{split}
e^{-x}( (\Lambda Q)'-\Lambda Q-\frac 12(1-\lambda) Q) &=
\frac 12 (3-\lambda) e^{-x} (Q'-Q) + \frac 12 (1-\lambda) x e^{-x}(Q''-Q')
\\ &=- \frac 12 (3-\lambda) (Q'+Q) + \frac 12 (1-\lambda) x (Q'+Q-e^{-x}Q^2).
\end{split}\end{equation*}
\end{proof}

\subsection{Technical claim \ref{CL:pp}}

\begin{claim}\label{CL:pp}
	\begin{equation}\label{eq:qq}
		\frac 1{18} \qun\qde=  e^{-\pyy} \left( \frac {\px \qun}{\qun}- \frac 13 \qun  - \frac {\px \qde}{\qde}- \frac 13 \qde\right) + \OO_{3/2}.
	\end{equation}
\end{claim}
\begin{proof}
	We distinguish the two regions
	$x-y_2> \frac \pyy 2$ and $x-y_2<\frac \pyy2$.
		
	- Case $x-y_2>\frac \pyy2$.
	For $x>y_2 + \frac \pyy  2$, we have
	$\qde(t,x)= 6 e^{-(x-y_1)-\pyy }+ O(e^{-2(x-y_1)-2\pyy })$ from \eqref{eq:as} in  Claim \ref{LE:A1}.	Thus, we obtain for such $x$,
	\begin{align*}
	& \qun\qde -  18 e^{-\pyy} \left( \frac {\px \qun}{\qun}- \frac 13 \qun  - \frac {\px \qde}{\qde}- \frac 13 \qde\right) \\ & = -18 e^{-\pyy} \left( \frac {\px \qun}{\qun}- \frac 13 \qun -\frac 13 e^{-(x-y_1)} \qun  - \frac {\px \qde}{\qde}- \frac 13 \qde\right) + \qun O(e^{-2(x-y_1)-2\pyy }).
	\end{align*}
	Note that $\qun O(e^{-2(x-y_1)-2\pyy })=\OO_{3/2}$.
	
	Next, it is a direct consequence of \eqref{eq:56} that
	\begin{equation}\label{eq:561}
		\frac {\px \qun}{\qun} -\frac 13 \qun = -  1 +\frac 13 e^{-(x-y_1)} \qun
		\quad \text{and}
	\quad 		\frac {\px \qde}{\qde} +\frac 13 \qde =   1 -\frac 13 e^{x-y_2} \qde.
	\end{equation}
	In particular, for $x>y_2 + \frac \pyy  2$, we obtain
	$$
		\frac {\px \qde}{\qde} +\frac 13 \qde =  -1+ \OO_{1/2}.
	$$	
	Therefore, we also obtain
	$$
	\qun\qde -  18 e^{-\pyy} \left( \frac {\px \qun}{\qun}- \frac 13 \qun  - \frac {\px \qde}{\qde}- \frac 13 \qde\right)= \OO_{3/2}		$$ in this region.
	
	- Case $x-y_1<-\frac \pyy 2$ (or equivalently $x-y_1<-\frac \pyy 2$).
 It is treated similarly.
\end{proof}

\subsection{Proof of Lemma \ref{LE:tF}}
First, we claim the following estimates.
\begin{claim}\label{CL:100}
 The following holds ($\omega \geq 0$)
\begin{equation}\label{eq:sq} 
\sqrt{\frac {1+\mu_j}{1+\lambda \mu_j}}=
 1+\frac 12 (1-\lambda) \mu_j -\frac 14 (1-\lambda)(1+3 \lambda) \mu_j^2+O(e^{-\frac 32 \YYzz}),
 \end{equation}
 \begin{equation}\label{eq:sQQ}
\qtud(t,x)+ |\px \qtud(t,x)| \leq C e^{-(1- 2 \mu_0) |x-y_j|}, 
 \end{equation}
 \begin{equation}\label{eq:sqQQ}\begin{split}
 &	\left| \qtud(t,x) - \left\{
	\qud(t,x) + \mu_j(t) \Lambda \qud(t,x) + \frac 12 \mu_j^2 (t) \Lambda^2 \qud(t,x)\right\} \right|\\ & \leq 
	C \mu_0^3 (1+|x-y_j(t)|^3) e^{-(1- 2 \mu_0) |x-y_j(t)|},
\end{split} \end{equation}
\begin{equation}\label{eq:sqQt}
		 (1+|x-y_1|^\omega+|x-y_2(t)|^{\omega})  e^{-(1- 2 \mu_0) |x-y_1|} e^{-(1- 2 \mu_0) |x-y_2|} =\OO_1,
		 \end{equation}
		 	\begin{equation}\label{eq:sqQb}
		\int(1+|x-y_1|^\omega+|x-y_2(t)|^{\omega})  e^{-(1- 2 \mu_0) |x-y_1|} e^{-(1- 2 \mu_0) |x-y_2|} dx\leq
		C(1+|\pyy|^{1+\omega})e^{-\pyy}.
	\end{equation}
\end{claim}

\begin{proof}
 We have
\begin{equation}\label{eq:sqp}\begin{split}
\sqrt{\frac {1+\mu_j}{1+\lambda \mu_j}}
& =(1+\frac 12 \mu_j-\frac 18 \mu_j^2 ) (1-\frac 12 \lambda \mu_j+\frac 38 \lambda^2 \mu_j^2)+O(\mu_j^{3})\\
&= 1+\frac 12 (1-\lambda) \mu_j -\frac 14 (1-\lambda)(1+3 \lambda) \mu_j^2+O(\mu_j^{3}),
\end{split}\end{equation}
and \eqref{eq:sq} follows.
Estimate \eqref{eq:sQQ} is clear from $Q(x)\leq C e^{-|x|}$ and \eqref{eq:sq}.

Now, we prove \eqref{eq:sqQQ}, using the Taylor formula (in the $\mu_j$ variable):
\begin{equation*}\begin{split}
	\qtud(t,x) &=  \qud(t,x) + \mu_j(t) \Lambda \qud(t,x) + \frac 12 \mu_j^2(t) \Lambda^2 \qud(t,x) \\ & + \frac 12  \mu_j^3(t) \int_0^1  (1-s)^2  \left(\frac {\partial^3 Q_\mu}{\partial\mu^3}\right)_{|\mu=s \mu_j(t)}(x-y_j)  ds.
\end{split}\end{equation*}
Note that, from the proof of \eqref{eq:ST} and elementary computations:
$$
\left|\left(\frac {\partial^3 Q_\mu}{\partial\mu^3}\right)_{|\mu=s \mu_j(t)}(x ) \right|
\leq C (1+|x |^3) e^{-(1- 2 \mu_0) |x |};
$$
\eqref{eq:sqQQ} follow.

\medskip

Proof of \eqref{eq:sqQb}. 
For $y_2<x<y_1$, we have
$$e^{-(1- 2 \mu_0) |x-y_1|}e^{-(1- 2 \mu_0) |x-y_2|}
= e^{-(1- 2 {\mu_0}) \pyy} \leq C e^{-\pyy},$$
since \eqref{eq:mh1} implies $ 2 {\mu_0 \pyy(t)} \leq C {\mu_0 \YYzz} \leq C'$ for $\YYzz$ large enough.

For $x>y_1>y_2$,
\begin{align*}
e^{-(1- 2 \mu_0) |x-y_1|}e^{-(1- 2 \mu_0) |x-y_2|}
&= e^{-(1- 2 {\mu_0})(2 x - y_1 -y_2)}
= e^{- 2 (1- 2 \mu_0) (x-y_1)} e^{- (1- 2 \mu_0) \pyy(t)}
\\ & \leq C e^{-\frac 32 |x-y_1|} e^{-\pyy(t)}.
\end{align*}
Arguing similarly for the case $x<y_2<y_1$, we prove \eqref{eq:sqQt} and \eqref{eq:sqQb}.
\end{proof}

We continue the proof of Lemma \ref{LE:tF}.
In order to expand $ \FF$, we   perform the following preliminary decomposition:
\begin{equation}\label{eq:70}\begin{split}
\FF & 
= 2 \qtde \px \qtun + 2 \qtun \px \qtde
 \\ &= 2 \qtde \left(\partial_x \qtun - \sqrt{\frac {1+\mu_1}{1+\lambda \mu_1}} \qtun\right) 
 + 2 \qtun \left(\partial_x \qtde + \sqrt{\frac {1+\mu_2}{1+\lambda \mu_2}}
 \qtde\right)\\
& + 2 \left(\sqrt{\frac {1+\mu_1}{1+\lambda \mu_1}}-
 \sqrt{\frac {1+\mu_2}{1+\lambda \mu_2}}\right) \qtun \qtde.
  \end{split}\end{equation}

 Using \eqref{eq:sq} and \eqref{eq:sqQQ}, we have
\begin{equation}\begin{split}
  \partial_x \qtun - \sqrt{\frac {1+\mu_1}{1+\lambda \mu_1}} \qtun
& = \partial_x (\qun + \mu_1 \Lambda \qun) 
 - \left(1+\tfrac 12 (1-\lambda) \mu_1 \right) (\qun+\mu_1 \Lambda\qun)+ \tilde \Theta \\
& =
\partial_x \qun -  \qun +\mu_1 (\partial_x \Lambda \qun -\Lambda \qun  - \tfrac 12 (1-\lambda) \qun  )
+\Theta,
\end{split}\end{equation}
where
$$ |\Theta|+|\tilde \Theta|\leq C |\mu_1|^2 (1+|x-y_1|^3) e^{-(1- 2 \mu_0) |x-y_1|}.$$
Thus, by \eqref{eq:sQQ} and \eqref{eq:sqQt}, we obtain
\begin{equation*}\begin{split}
& 2 \qtde \left(\partial_x \qtun - \sqrt{\frac {1+\mu_1}{1+\lambda \mu_1}} \qtun\right)  
\\
& = 2(\qde+\mu_2 \Lambda \qde ) \left[ \px \qun -  \qun  +\mu_1 (\px (\Lambda \qun)  - \qun  - \tfrac 12 (1-\lambda) \qun  )\right]+ \OO_2.
\end{split}\end{equation*}
Using \eqref{eq:as}, \eqref{eq:ST} and $x-y_2=x-y_1+\pyy$, we obtain
\begin{equation*}\begin{split}
& 2 \qtde \left(\partial_x \qtun - \sqrt{\frac {1+\mu_1}{1+\lambda \mu_1}} \qtun\right)  \\ &
= 12e^{-\pyy}e^{-(x-y_1)}
\left[1+\mu_2 (1-\tfrac 12 (1-\lambda) (x-y_1+\pyy))\right]\times
\\ & \quad \times \left[\partial_x \qun -  \qun +\mu_1 (\partial_x \Lambda \qun -\Lambda \qun - \tfrac 12 (1-\lambda) \qun  )
 \right]  +\OO_{2}.
\end{split}\end{equation*}
 Using \eqref{eq:56}, \eqref{eq:57} and then \eqref{eq:mh2},
\begin{equation}\label{eq:FU3}\begin{split}
& 2 \qtde \left(\partial_x \qtun - \sqrt{\frac {1+\mu_1}{1+\lambda \mu_1}} \qtun\right) 
=
- 12 e^{-\pyy} (\partial_x \qun+\qun)
\\&+12\mu_1e^{-\pyy}
\left[-\tfrac 12 (3-\lambda) (\partial_x \qun+\qun) 
+ \tfrac 12 (1-\lambda) (x-y_1)(\qun^2 - 2 (\partial_x \qun+\qun))
\right]
\\ & -12\mu_2e^{-\pyy} (1-\tfrac 12 (1-\lambda) (x-y_1+\pyy)) (\partial_x \qun+\qun)+ \OO_{2}
\\
&=
- 12 e^{-\pyy} (\partial_x \qun+\qun)
- 6  (1-\lambda) \mu_1 \pyy e^{-\pyy}   (\partial_x \qun+\qun)
\\ & +12\mu_1e^{-\pyy}
\left[- \tfrac 12 (3-\lambda) (\partial_x \qun+\qun) 
+ \tfrac 12 (1-\lambda) (x-y_1)(\qun^2 -2 (\partial_x \qun+\qun))
\right]
\\ & -12\mu_2e^{-\pyy} (1-\tfrac 12 (1-\lambda) (x-y_1))(\partial_x \qun+\qun)+\OO_{2}\\
& =- 12 e^{-\pyy} (\partial_x \qun+\qun)
- 6  (1-\lambda) \mu_1 \pyy e^{-\pyy}   (\partial_x \qun+\qun)
\\ & +6\mu_1e^{-\pyy}(1-\lambda)
\left[-  (\partial_x \qun+\qun) 
+  (x-y_1)(\qun^2 -3 (\partial_x \qun+\qun))\right]+\OO_{2}.
\end{split}\end{equation}

Next, by similar computations,
\begin{equation}\begin{split}
& \partial_x \qtde + \sqrt{\frac {1+\mu_2}{1+\lambda \mu_2}} \qtde
  =     12 e^{-\pyy} (-\partial_x \qde+\qde)
+ 6  (1-\lambda) \mu_2 \pyy e^{-\pyy}   (-\partial_x \qde+\qde)
\\ & +6\mu_2 e^{-\pyy}(1-\lambda)
\left[(-\partial_x \qde+\qde) 
- (x-y_2)(\qde^2 -3 (-\partial_x \qde+\qde))\right]+\OO_{2}
\end{split}\end{equation}
 
Finally, using Claim \ref{CL:100},
\begin{equation*} 
  2 \left(\sqrt{\frac {1+\mu_1}{1+\lambda \mu_1}}-
 \sqrt{\frac {1+\mu_2}{1+\lambda \mu_2}}\right) \qtun \qtde
  =(1-\lambda) (\mu_1-\mu_2) \qun \qde +  \OO_{2}.
 \end{equation*}
\subsection{Approximate antecedent of $\qun\qde$}
\begin{claim}\label{LE:Q1Q2}
Let  $x_j=x-y_j$. Then
\begin{equation*}\begin{split}
&
\px ( - \px^2((x_1+x_2) \qun\qde) + (x_1+x_2) \qun\qde - 2 (\qun+\qde) (x_1+x_2)\qun\qde)
 \\&= 2 \qun\qde  	- \pyy(3 (\px \qun-\qun)+ \px( \qun^2)) \qde  + \pyy \qun^2  \px \qde\\ 
& 	+  \pyy(3 (\px \qde +\qde)+ \px( \qde^2) )\qun -\pyy \qde^2 \px \qun 
\\ &
+	2 [\qun^2   - x_1 \px \qun^2 - 3 x_1 (\px \qun-\qun)- 3 (\qun - \px \qun) ]\qde\\
	&  + 2 [ x_1 \qun^2 - 3  (\px \qun -\qun)] \px \qde \\
	&  + 2 [ \qde^2  - x_2 \px \qde^2- 3 x_2 (\px \qde +\qde) -3 (\qde + \px \qde)  ] \qun \\
	& + 2 [  x_2 \qde^2 - 3 (\px \qde+\qde)]   \px \qun  
	  .
\end{split}\end{equation*}
\end{claim}
\begin{proof}
 For any two functions $F_1$, $F_2$,
the following holds true
\begin{equation}\label{eq:FF}\begin{split}
	&
	\px ( - \px^2(F_1(x_1)F_2(x_2)) + F_1(x_1) F_2(x_2) - 2 (R_1+R_2) (F_1(x_1) F_2(x_2)))
	\\& = 
	\px \left[ F_2(x_2) ( L F_1)(x_1) + F_1(x_1) ( L F_2)(x_2) - 2  F_1'(x_1)F_2'(x_2) - F_1(x_1) F_2(x_2) \right] \\
	& =(L F_1)'(x_1) F_2(x_2) +  (L F_2)'(x_2) F_1(x_1)\\
	&\quad + ( L F_1-F_1- 2 F_1'')(x_1)F_2'(x_2) + 
	(L F_2-F_2- 2 F_2'')(x_2) F_1'(x_1).
\end{split}\end{equation}

Now, we apply   formula \eqref{eq:FF} with $F_1=xQ$ and $F_2=Q$.
Note that (see Claim \ref{LE:A1} (ii))
$$LQ= - Q^2,\quad 
 LQ-Q -2 Q''=-3Q+Q^2,$$
$$L (xQ)= x  LQ - 2 Q'
= - x Q^2 - 2 Q',$$
$$
L(xQ)- xQ - 2 (xQ)'' =
-x(Q^2+Q)-2 Q' - 2 x(Q - Q^2) -4 Q'
= x (-3 Q + Q^2) - 6 Q'.
$$
Thus, from \eqref{eq:FF}
\begin{equation}\label{eq:xQQ12}\begin{split}
	& \px ( - \px^2(x_1 \qun\qde) + x_1 \qun\qde - 2 (\qun+\qde) x_1\qun\qde)
	\\
 	& =   (\qun^2- x_1 \px \qun^2) \qde    + x_1 \qun^2 \px \qde- \px( \qde^2) (x_1\qun) 	+\qde^2(\qun+x_1 \px \qun) \\
	&  - 5 \qun\qde - 3 x_1 \px \qun \qde - 3  x_1 \qun \px \qde - 6 \px \qun \px \qde.
\end{split}\end{equation}
Similarly,
\begin{equation}\label{eq:xQQ21}\begin{split}
	& \px ( - \px^2(x_2 \qun\qde) + x_2 \qun\qde - 2 (\qun+\qde) x_2\qun\qde)
	\\&  = (\qde^2  - x_2 \px \qde^2) \qun + x_2 \qde^2 \px \qun - \px( \qun^2)(x_2\qde)	+ \qun^2(\qde+x_2 \px \qde)\\
	&  - 5 \qde\qun - 3 x_2 \px \qde \qun - 3  x_2 \qde \px \qun - 6 \px \qde \px \qun.
\end{split}\end{equation}

Therefore, summing up,
\begin{equation}\label{eq:tc}\begin{split}
	& 
	\px ( - \px^2((x_1+x_2) \qun\qde) + (x_1+x_2) \qun\qde - 2 (\qun+\qde) (x_1+x_2)\qun\qde)
	  \\ & =
	(\qun^2   - x_1 \px \qun^2) \qde +x_1 \qun^2 \px \qde-\px( \qun^2) (x_2\qde)	+\qun^2  (\qde+x_2 \px \qde) \\
	&  + (\qde^2  - x_2 \px \qde^2) \qun +x_2 \qde^2 \px \qun - \px( \qde^2) (x_1\qun)	+\qde^2(\qun+x_1\px \qun)\\
	&  - 10 \qun\qde  -12 \px \qun \px \qde \\& 
	  - 3   x_1 \px \qun \qde - 3  x_1 \qun \px \qde  - 3   x_2 \px \qde \qun - 3   x_2 \qde \px \qun  .
\end{split}\end{equation}
The terms in the last line of \eqref{eq:tc} are handled as follows (recall that $x_2-x_1=y$)
\begin{equation*}\begin{split}
& 3   x_1 \px \qun \qde + 3  x_1 \qun \px \qde  + 3   x_2 \px \qde \qun + 3   x_2 \qde \px \qun 
\\ 
&=3 x_1 (\px \qun -\qun) \qde + 3 x_1 \qun (\qde+\px \qde) 
+ 3 x_2 (\px \qde +\qde) \qun + 3 x_2 \qde (\px \qun-\qun)\\
&=
6 x_1 (\px \qun-\qun) \qde + 3 \pyy(\px \qun-\qun)\qde
+ 6 x_2 (\px \qde +\qde) \qun - 3\pyy(\px \qde +\qde)\qun. 
\end{split}\end{equation*}
For the term $12 \px \qun \px \qde$ in \eqref{eq:tc}, we observe
\begin{equation*}\begin{split}
12 \px \qun \px \qde&= 6 (\px \qun -\qun) \px \qde 
+ 6 (\px \qde+\qde) \px \qun  + 6 \qun \px \qde - 6 \qde \px \qun\\
&=6 (\px \qun -\qun) \px \qde 
+ 6 (\px \qde+\qde) \px \qun  \\ & + 6 (\qde + \px \qde)  \qun 
+ 6 (\qun - \px \qun) \qde  - 12 \qun \qde.
\end{split}\end{equation*}

Thus, we obtain
\begin{equation*}\begin{split}
&
\px ( - \px^2((x_1+x_2) \qun\qde) + (x_1+x_2) \qun\qde - 2 (\qun+\qde) (x_1+x_2)\qun\qde)
 \\ & =
	(\qun^2   - x_1 \px \qun^2) \qde + x_1 \qun^2 \px \qde- \px( \qun^2) (x_2\qde)	+\qun^2 (\qde+x_2 \px \qde) \\
	&  - 6 x_1 (\px \qun-\qun) \qde -	6 (\px \qun -\qun) \px \qde 
	- 6 (\qun - \px \qun)\qde 
	\\
	&  + (\qde^2  - x_2 \px \qde^2) \qun +x_2 \qde^2 \px \qun- \px( \qde^2) (x_1\qun)	+\qde^2(\qun+x_1\px \qun)\\
	&
- 6 x_2 (\px \qde +\qde) \qun 
- 6 (\px \qde+\qde) \px \qun   - 6 (\qde + \px \qde)  \qun\\
&+ 3 \pyy(-\px \qun+\qun)\qde + 3\pyy(\px \qde +\qde)\qun
	  +2 \qun\qde
	  \\&= 2 \qun\qde  	- \pyy(3 (\px \qun-\qun)+ \px( \qun^2)) \qde  + \pyy \qun^2  \px \qde\\ 
& 	+  \pyy(3 (\px \qde +\qde)+ \px( \qde^2) )\qun -\pyy \qde^2 \px \qun 
\\ &
+	2 [\qun^2   - x_1 \px \qun^2 - 3 x_1 (\px \qun-\qun)- 3 (\qun - \px \qun) ]\qde\\
	&  + 2 [ x_1 \qun^2 - 3  (\px \qun -\qun)] \px \qde \\
	&  + 2 [ \qde^2  - x_2 \px \qde^2- 3 x_2 (\px \qde +\qde) -3 (\qde + \px \qde)  ] \qun \\
	& + 2 [  x_2 \qde^2 - 3 (\px \qde+\qde)]   \px \qun  
	  ,
\end{split}\end{equation*}
	  and the proof of Claim \ref{LE:Q1Q2} is complete.
\end{proof}

\section{Modulation and monotonicity arguments}\label{AP:B}
 
  \subsection{Proof of Lemma \ref{PR:de}}
Let 
$$
	\mathcal{V}(\omega_0,y_0) = \{ u \in H^1(\RR);
	\inf_{y_1-y_2>y_0}\| u - V(x;(0,0,y_1,y_2))\|_{H^1} \leq \omega_0\},
$$
where $V(x;\Gamma)$ is defined in Proposition \ref{PR:app}.

\begin{lemma}[Time independent modulation]\label{le:huit}
There exist $\omega_0$, $\bar y_0>0$ and a unique $C^1$
map $\Gamma=(\mu_1,\mu_2,y_1,y_2): \mathcal{V}(\omega_0,\bar y_0)\to (0,\infty)^2\times \RR^2$
such that if $u\in \mathcal{V}(\omega,y_0)$ for $0<\omega\leq \omega_0$,
$y_0\geq \bar y_0$ and 
$$
	\varepsilon (x) = u(x) - V(x;\Gamma),
$$
then, for $j=1,2$,
$$
	\int \varepsilon (1-\lambda \partial_x^2)   Q_{\mu_j}(.-y_j) =
	\int \varepsilon (1-\lambda \partial_x^2)   Q_{\mu_j}'(.-y_j) = 0
$$
$$
	y_1-y_2 > y_0 - C \omega ,\quad 
	\|\varepsilon\|_{H^1} + |\mu_1 |+|\mu_2 |\leq C \omega .
$$
\end{lemma}

\begin{proof}
	The proof, based on the implicit function theorem, is  similar to the one of Lemma 8 in \cite{MMT} (see also \cite{DiMa} for the \eqref{eq:BBM} case), the only difference being that the modulation   uses the map
$(\mu_1,\mu_2,y_1,y_2)\mapsto V(x;(\mu_1,\mu_2,y_1,y_2)) $ instead of the family of sums of two solitons. By the properties of $V$ (see  \eqref{eq:ap2} and below \eqref{eq:de1}) and \eqref{eq:58}, the nondegeneracy condition is the same as in \cite{DiMa}.
\end{proof}

The existence, uniqueness and continuity of $\Gamma(t)$ is a consequence of Claim \ref{le:huit}.
	The $C^1$ regularity of $\Gamma(t)$ is obtained by standard regularization arguments and the equation of $\varepsilon(t)$ which is deduced easily from {\eqref{eq:BBM}} and \eqref{eq:VV}.

	Next, we prove the estimates on $\dot \Gamma(t)$, i.e. \eqref{eq:ga}, omitting standard regularization arguments to justify the formal computations. First, we expand $0=\frac d{dt}\int \varepsilon (1-\lambda \partial_x^2) \qtud$. Using   \eqref{eq:ep}, we obtain ($k\neq j$)
\begin{equation*}
\begin{split}
	0 &= \frac d{dt} \int \varepsilon (1-\lambda \partial_x^2) \qtud
	 \\ &= \int \varepsilon (1-\lambda \partial_x^2) \partial_t \qtud + \int (\partial_x^2 \varepsilon - \varepsilon + 2V \varepsilon) \partial_x \qtud   + \int \varepsilon^2 \partial_x  \qtud- \int E  \qtud  \\
	& - (\dot \mu_j - {\cal M}_j ) \int  (1-\lambda \partial_x^2)\frac {\partial V}{\partial \mu_j}  \qtud -
	(\dot \mu_k - {\cal M}_k) \int (1-\lambda \partial_x^2) \frac {\partial V}{\partial \mu_k}  \qtud\\
	& + (\mu_j - \dot \yy_j - {\cal N}_j) \int (1-\lambda \partial_x^2) \frac {\partial V}{\partial y_j}   \qtud
	+ (\mu_k - \dot \yy_k - {\cal N}_k) \int (1-\lambda \partial_x^2) \frac {\partial V}{\partial y_k}   \qtud.
\end{split}
\end{equation*}

We claim the following estimates.
\begin{claim}
	Assuming \eqref{eq:gga},
	\begin{equation}\label{eq:de2}
	\left|\int \qtun \qtde\right|\leq C (\pyy+1) e^{-\pyy}, 
	\end{equation}
\begin{equation}\label{eq:de1}
j=1,2,\quad 
 \left\|\frac {\partial V}{\partial \mu_j} - \Lambda \qtud\right\|_{H^1}+ \left\|\frac {\partial V}{\partial y_j} + \partial_x \qtud\right\|_{L^\infty}+ \frac 1{\sqrt{y}} \left\|\frac {\partial V}{\partial y_j} + \partial_x \qtud\right\|_{H^1}\leq C e^{-\pyy}.
\end{equation}
\end{claim}
Indeed, under assumption \eqref{eq:gga}, \eqref{eq:de2} is a consequence of \eqref{eq:sQQ}
and \eqref{eq:sqQb},
and \eqref{eq:de3} is a   consequence of \eqref{eq:ap0}, \eqref{eq:V} and the properties of $A_j$, $B_j$ and $D_j$ (see \eqref{eq:ABD}).

\medskip

By \eqref{eq:de1}, \eqref{eq:de2}, \eqref{eq:ap2bis}, $L_{\mu_j}   Q_{\mu_j}'=0$ (see Claim \ref{LE:A1}) and $\int  (1-\lambda \partial_x^2) \partial_x \qtud  \qtud =0$, we get ($j\neq k$)
\begin{equation*}
\begin{split}
	0 & = \dot \mu_j \int \varepsilon (1-\lambda \partial_x^2) \Lambda \qtud +  
	(\mu_j - \dot \yy_j) \int \varepsilon(1-\lambda \partial_x^2) \partial_x \qtud + \|\varepsilon(t)\|_{L^2}O((\pyy+1)e^{-\pyy}) \\ & + O(\|\varepsilon\|_{L^2}^2)
	   - \int E  \qtud \\
	& - (\dot \mu_j - {\cal M}_j ) \left(\int  (1-\lambda \partial_x^2)  \Lambda \qtud \qtud +O(e^{-\frac 12 \pyy})\right)+
	(\dot \mu_k - {\cal M}_k) O(y^2 e^{-\pyy})\\
	& 
	+  (\mu_j - \dot \yy_j - {\cal N}_j)O(e^{- \pyy})
	+  (\mu_k - \dot \yy_k - {\cal N}_k) 
	O(\pyy e^{-\pyy}).
\end{split}
\end{equation*}
Hence, by $\left| \int (1-\lambda \partial_x^2)\Lambda \qtud   \qtud\right|\geq c_0>0$ (see \eqref{eq:58}),
for $\pyy$ large and $\varepsilon$ small, we get
\begin{equation*}\begin{split} 
& |\dot \mu_j - {\cal M}_j| \leq C \big[\|\varepsilon\|_{L^2}^2 + \pyy e^{-\pyy} \|\varepsilon\|_{L^2} + \int |E \qtud| \big]\\
& + C e^{-\pyy}   |\mu_j - \dot \yy_j - {\cal N}_j |   + C y^2e^{-  \pyy}  |\dot \mu_k - {\cal M}_k| +  C \pyy e^{- \pyy} |\mu_k - \dot \yy_k - {\cal N}_k | .
\end{split}\end{equation*}
Similarly, expanding $0=\frac d{dt}\int \varepsilon (1-\lambda \partial_x^2)\partial_x \qtud$, we obtain
\begin{equation*}\begin{split} 
& |\mu_j - \dot \yy_j - {\cal N}_j | \leq C \big[\|\varepsilon\|_{L^2}+  \int | E \partial_x \qtud| \big]\\
&+ C  (\|\varepsilon\|_{L^2}+e^{-\frac 12 \pyy}) |\dot \mu_j - {\cal M}_j|   + C e^{- \frac 12 \pyy}  |\dot \mu_k - {\cal M}_k| + C \pyy e^{-\pyy} |\mu_k - \dot \yy_k - {\cal N}_k | .
\end{split}\end{equation*}
Combining these   estimates, for $\pyy$ large and $\varepsilon$ small,   \eqref{eq:ga} is proved.

\subsection{Proof of Proposition \ref{PR:cFG}.}
The proof of Proposition \ref{PR:cFG} is inspired by the proof of Proposition 3.1 in \cite{MMkdv4}. However, it is technically more involved in the BBM case. We refer to \cite{Mi2}, \cite{Di}, \cite{DiMa} and \cite{MMM} for previous similar arguments for the (BBM) equation.

\medskip

The proof of \eqref{eq:FGcoer} is standard, see for example Lemma 4 in \cite{MMT} and \cite{DiMa}. Recall that it is based on coercivity property of the operator $L$ under orthogonality conditions, see Claim \ref{LE:A2} (v).

\medskip

	We continue with the following claim:
\begin{claim}\label{CL:ab}
Let $\varphi(x)$ be defined by \eqref{eq:ph}.
If $a(x),b(x)\in L^2$ are such that $a-\lambda \partial_x^2 a= b$ then
\begin{equation}\label{eq:ab}
	\left(1-(8\rho)^2 \right) \int a^2 \varphi' + 2 \lambda \int (\partial_x a)^2 \varphi'
	+ \lambda^2 \int (\partial_x^2 a)^2 \varphi' \leq \int b^2 \varphi'. 
\end{equation}
\end{claim}
Indeed,  integrating by parts and then using \eqref{eq:ph},
\begin{equation*}\begin{split} 
\int b^2 \varphi' & =
\int (a-\lambda \partial_x^2 a)^2 \varphi' 
= \int a^2 \varphi' + 2 \lambda \int (\partial_x a)^2 \varphi' - \lambda \int a^2 \varphi''' + \lambda^2 \int (\partial_x^2 a)^2 \varphi'
\\ & \geq  \left( 1 - (8\rho)^2 \right) \int a^2 \varphi' + 2 \lambda \int (\partial_x a)^2 \varphi'
	+ \lambda^2 \int (\partial_x^2 a)^2 \Phi'.
\end{split}\end{equation*} 

\medskip

\noindent\textbf{$\bullet$ Case $\mu_1(t)\geq \mu_2(t)$.}\quad 
We first  claim the following technical estimates, as consequences of 
\eqref{eq:ph}, \eqref{eq:en1}, \eqref{eq:ap0} and \eqref{eq:V}.

\begin{claim}\label{CL:cl}
\begin{equation}
\|V-\qtun-\qtde\|_{L^\infty} \leq C e^{-y},
\end{equation}
\begin{equation}\label{eq:dc1}
	\|V \partial_x \Phi\|_{L^\infty} \leq C (|\mu_1|+|\mu_2|) e^{-2 \rho \pyy  },
\end{equation}
\begin{equation}\label{eq:dc10}
	\|(\Phi - \mu_j ) e^{-\frac 12 |x-\yy_j|}\|_{L^\infty} \leq C (|\mu_1|+|\mu_2|) e^{-2 \rho \pyy  },
\end{equation}
\begin{equation}\label{eq:dc2}
	\| \Phi \partial_x V  - \sum_{j=1,2} \mu_j \partial_x \qtud \|_{L^\infty}
	 \leq C (|\mu_1|+|\mu_2|) e^{- 2 \rho \pyy  } + C e^{-\pyy},
 \end{equation}
\begin{equation}\label{eq:de3}
 \bigg\| {\partial_t V}- \sum_{j=1,2}\left\{\dot \mu_j \Lambda \qtud -\dot y_j  \partial_x \qtud\right\}\bigg\|_{L^\infty}\leq
C e^{- \pyy}.
\end{equation}
\end{claim}

Let
$$
	\Theta = \|\varepsilon\|_{L^2}^2 
\left[e^{- \frac 34  \pyy}  + (|\mu_1|+|\mu_2|+\|\varepsilon\|_{L^2} )( e^{- 2 \rho \pyy }      
+  \|\varepsilon\|_{L^2})  | \right]
+  \|\varepsilon\|_{L^2} \|E\|_{L^2}.
$$
Let us compute $\frac d{dt}{\cal F}_+(t)$:
\begin{equation*}\begin{split}
\frac 12 \frac d{dt}{\cal F}_+(t) & = \int \partial_t \varepsilon \left(-\partial_x^2 \varepsilon + \varepsilon - ((\varepsilon +V)^2-V^2) +
[(1-\lambda \partial_x^2) (\varepsilon \Phi)]\right) 
\\
& - \lambda \int (\partial_x\partial_t \varepsilon)   \varepsilon \partial_x \Phi 
+ \frac 12 \int \partial_t \Phi [\lambda (\partial_x \varepsilon)^2 + \varepsilon^2]\\
& - \int \partial_t V \left((\varepsilon+V)^2 - V^2 - 2V\varepsilon \right)
=F_1+F_2+F_3+F_4.
\end{split}\end{equation*}
Observe that $\partial_x \Phi = (\mu_1-\mu_2) \varphi'\geq 0$ by \eqref{eq:en1}.

Using  \eqref{eq:ep} and then by direct computations and estimates, we claim the following estimates, which imply immediately 
\eqref{eq:cF}.
\begin{claim}\label{cl:4.1}
\begin{align}
	F_1  & \leq -\frac 32 \int (\partial_x \varepsilon)^2 \partial_x \Phi - \frac 38 \int \varepsilon^2 \partial_x \Phi
	- \int \varepsilon^2 (\mu_1 \partial_x \qtun + \mu_2 \partial_x \qtde)  , \nonumber\\
	& +\sum_{j=1,2}(\dot \mu_j - {\cal M}_j ) \int \varepsilon^2 \Lambda \qtud 
	+\sum_{j=1,2}(\mu_j -\dot \yy_j - {\cal N}_j ) \int \varepsilon^2 \partial_x \qtud +C \Theta\label{eq:F1}\\
	F_2 & \leq \frac 32 \int (\partial_x \varepsilon)^2 \partial_x \Phi + \frac 38 \int \varepsilon^2 \partial_x \Phi+C \Theta,\label{eq:F2}\\
	F_3 & \leq C\Theta,\label{eq:F3}\\
	F_4 & \leq  \int \varepsilon^2 (\mu_1 \partial_x \qtun + \mu_2 \partial_x \qtde)\nonumber \\
	& -\sum_{j=1,2}(\dot \mu_j - {\cal M}_j ) \int \varepsilon^2 \Lambda \qtud  
	-\sum_{j=1,2}(\mu_j -\dot \yy_j - {\cal N}_j ) \int \varepsilon^2 \partial_x \qtud  +C \Theta. \label{eq:F4}
\end{align}
\end{claim}

Indeed,
\begin{equation*}\begin{split} 
F_1 & = 
- \int \left(-\partial_x^2 \varepsilon + \varepsilon  - ((\varepsilon +V)^2-V^2)\right) (\Phi   \partial_x \varepsilon + \varepsilon \partial_x \Phi)
\\&-
 \int E (1-\lambda \partial_x^2)^{-1}\left(-\partial_x^2 \varepsilon + \varepsilon  +[(1-\lambda \partial_x^2) (\varepsilon \Phi)] - ((\varepsilon +V)^2-V^2)\right)
\\
& - \sum_{j=1,2}
(\dot \mu_j- {\cal M}_j) \int \frac {\partial V}{\partial \mu_j} \left(-\partial_x^2 \varepsilon + \varepsilon  +[(1-\lambda \partial_x^2) (\varepsilon \Phi)]- ((\varepsilon +V)^2-V^2)\right)\\
& + \sum_{j=1,2}
(\mu_j - \dot \yy_j - {\cal N}_j) \int \frac {\partial V}{\partial y_j}  \left(-\partial_x^2 \varepsilon + \varepsilon  +[(1-\lambda \partial_x^2) (\varepsilon \Phi)]- ((\varepsilon +V)^2-V^2)\right)\\
&= F_{1,1}+F_{1,2}+F_{1,3}+F_{1,4}.
\end{split}\end{equation*} 
By \eqref{eq:ph}, Claim \ref{CL:cl} and several integration by parts,
 we get 
\begin{equation*}\begin{split}
F_{1,1} &  =
- \frac 32 \int (\partial_x\varepsilon)^2\partial_x \Phi  
- \frac 12 \int \varepsilon^2 \partial_x \Phi
+ \frac 12 \int \left(3\varepsilon^2 V + \frac 43 \varepsilon^3\right)\partial_x \Phi
+ \frac 12 \int \varepsilon^2 \partial_x^3\Phi 
- \int \Phi \varepsilon^2  \partial_x V \\
& \leq - \frac 32 \int (\partial_x\varepsilon)^2\partial_x \Phi  
- \frac 12\left(1-(4\rho)^2-C \|\varepsilon\|_{H^1} \right) \int \varepsilon^2 \partial_x \Phi
- \int \varepsilon^2 (\mu_1 \partial_x \qtun + \mu_2 \partial_x \qtde)
\\& +C (|\mu_1|+|\mu_2|)\,  \|\varepsilon\|_{L^2}^2  e^{-\rho \pyy} + Ce^{ -\frac 34 \pyy}\|\varepsilon\|_{L^2}^2 ;
 \end{split}\end{equation*}
and
\begin{equation*}
	|F_{1,2}|\leq C \|E\|_{L^2} \|\varepsilon\|_{L^2}. 
\end{equation*}
For $F_{1,3}$, we use  \eqref{eq:de1}, \eqref{eq:Qm}, \eqref{eq:or} and \eqref{eq:dc10}, so that
\begin{align*}
& \left | [ (1-\lambda \partial_x^2) \Lambda \qtud ] \Phi - [(1-\lambda \partial_x^2) \Lambda \qtud ] \mu_j \right|
\\
&\leq C e^{-\frac 34 {|x-\yy_j(t)|} } |\Phi-\mu_j|
\leq C (|\mu_1|+|\mu_2|) e^{-\frac 12 {|x-\yy_j(t)|}  } e^{-2 \rho {\pyy} },
\end{align*}
and
\begin{equation} 
	F_{1,3}=  \sum_{j=1,2}(\dot \mu_j - {\cal M}_j ) \int \varepsilon^2 \Lambda \qtud +
	e^{-2 \rho \pyy }(|\mu_1|+|\mu_2|+\|\varepsilon\|_{L^2})  O(\|\varepsilon\|_{L^2}) 
	\sum_{j=1,2}|\dot \mu_j - {\cal M}_j |.
\end{equation}
Similarly,
\begin{equation} 
	F_{1,4} =\sum_{j=1,2}(\mu_j -\dot \yy_j - {\cal N}_j ) \int \varepsilon^2 \partial_x \qtud
 +
	e^{-2 \rho \pyy}  (|\mu_1|+|\mu_2|+\|\varepsilon\|_{L^2}) O(\|\varepsilon\|_{L^2})   \sum_{j=1,2}|\mu_j -\dot \yy_j - {\cal N}_j |.
\end{equation}
Using \eqref{eq:ga}, \eqref{eq:F1} follows.

\medskip

Let $z_1, z_2$ such that $z_1 - \lambda \partial_x^2 z_1 = \varepsilon$, 
$z_2 - \lambda \partial_x^2 z_2 = 2V \varepsilon + \varepsilon^2$. Then, using the equation of $\varepsilon$:
\begin{equation*}\begin{split}
F_2
& = \lambda \int \partial_x^2 (1-\lambda \partial_x^2)^{-1} (\partial_x^2 \varepsilon -\varepsilon + 2V \varepsilon+ \varepsilon^2)
 \varepsilon \partial_x \Phi  + \lambda \int \partial_x(1-\lambda \partial_x^2)^{-1} E \varepsilon \partial_x \Phi \\
& + \lambda \sum_{j=1,2} (\dot \mu_j - {\cal M}_j) \int \partial_x \frac {\partial V}{\partial \mu_j}  \varepsilon \partial_x \Phi
 - \lambda \sum_{j=1,2} (\mu_j - \dot \yy_j - {\cal N}_j) \int \partial_x \frac {\partial V}{\partial y_j}
\varepsilon \partial_x \Phi\\
& = F_{2,1}+F_{2,2}+F_{2,3}+F_{2,4}.
\end{split}\end{equation*}
\begin{equation*}\begin{split}
F_{2,1}
& = - \lambda \int \partial_x (\partial_x^2 z_1 + z_2) \partial_x \varepsilon \partial_x \Phi
- \lambda \int \partial_x(\partial_x^2 z_1 + z_2) \varepsilon \partial_x^2  \Phi
- \lambda \int \partial_x^2 z_1 \varepsilon \partial_x \Phi  
 \\
&\leq \frac 12 \int \left[ \lambda^2   (\partial_x^3 z_1)^2  +   (\partial_x \varepsilon)^2  + 2 \lambda^2 (\partial_x^2 z_1)^2  +\frac 12   \varepsilon^2  + 4 \lambda^2   (\partial_x z_2)^2  +  \frac 14  (\partial_x \varepsilon)^2\right] \partial_x \Phi \\
& + 2 \rho \int\left[ \lambda^2   (\partial_x^3 z_1)^2  +   \varepsilon^2    + 4 \lambda^2   (\partial_x z_2)^2  +  \frac 14   \varepsilon^2
\right] \partial_x \Phi,
\end{split}\end{equation*}
by using Cauchy Schwarz inequality and \eqref{eq:ph}.
For $\rho$ small enough, using Claim \ref{CL:ab} and \eqref{eq:dc1}, we obtain
\begin{equation*}\begin{split}
|F_{2,1}|
& \leq  \frac 32 \int (\partial_x \varepsilon)^2 \partial_x \Phi + \left(\frac 14 + 8\rho +C \|\varepsilon\|_{H^1}+ C e^{-\pyy} \right)   \int \varepsilon^2 \partial_x \Phi+ C (|\mu_1|+|\mu_2|) \|\varepsilon\|_{L^2}^2 e^{- 2 \rho \pyy}\\
 & \leq \frac 32 \int (\partial_x \varepsilon)^2 \partial_x \Phi + \frac 38 \int \varepsilon^2 \partial_x \Phi
+ C \Theta.
\end{split}\end{equation*}
The term $F_{2,2}$ is estimated as $F_{1,2}$ and by arguments previously used, we also obtain
$ 
	|F_{2,3}| + |F_{2,4}| \leq C\Theta.
$
Thus,  \eqref{eq:F2} is proved.

\medskip

Next,
\begin{equation*}\begin{split} 
F_3=
& \frac 12 \int (\dot \mu_1 \varphi + \dot \mu_2 (1-\varphi)) [\lambda (\partial_x \varepsilon)^2 + \varepsilon^2],
\end{split}\end{equation*}
so that by $|{\cal M}_j|\leq C e^{-\pyy}$ and \eqref{eq:ga},
\begin{equation}
|F_3|\leq  C (e^{-\pyy}+ \sum_{j=1,2} |\dot \mu_j - {\cal M}_j|)  \|\varepsilon\|_{H^1}^2\leq C \Theta.
\end{equation}

Finally, by \eqref{eq:de3},
\begin{equation*}\begin{split} 
F_4  & = -\int \partial_t V \varepsilon^2 = - \int \sum_{j=1,2}\left(\dot\mu_j \Lambda \qtud  - \dot \yy_j \partial_x \qtud 
\varepsilon^2\right) + O(e^{- \pyy}) \|\varepsilon\|_{L^2}^2   \\
&= \sum_{j=1,2} \mu_j \int \varepsilon^2 \partial_x \qtud 
- \sum_{j=1,2} (\dot \mu_j - {\cal M}_j ) \int\varepsilon^2 \Lambda \qtud \\& 
- \sum_{j=1,2} ( \mu_j - \dot \yy_j - {\cal N}_j) \int\varepsilon^2 \partial_x \qtud
 + O(e^{- \pyy}) \|\varepsilon\|_{L^2}^2 ,
\end{split}\end{equation*} 
which proves \eqref{eq:F4}.

\medskip

\textbf{$\bullet$ Case $\mu_1(t)\leq \mu_2(t)$.}\quad 
Since $\mu_2(t)\geq \mu_1(t)$ we have
$\frac 1{(1+\mu_1(t))^2}\geq  \frac {1}{(1+\mu_2(t))^2}$,
$\partial_x \Phi_1\geq 0$ and $\partial_x \Phi_2\leq 0$. Note also that by explicit computations, for $\mu_j$ small enough:
\begin{equation}\label{eq:cP}
	\left| \partial_x  \Phi_2  + \frac 12 \partial_x\Phi_1 \right|\leq C(|\mu_1|+|\mu_2|) \partial_x\Phi_1. 
\end{equation}

Let us compute $\frac d{dt}{\cal F}_-(t)$:
\begin{equation*}\begin{split}
\frac 12 \frac d{dt}{\cal F}_-(t) & = \int \partial_t \varepsilon\left(-\partial_x^2 \varepsilon + \varepsilon - ((\varepsilon +V)^2-V^2)\right)\Phi_1 - \int \partial_t \varepsilon \partial_x \varepsilon \partial_x \Phi_1\\
& + \int (1-\lambda \partial_x^2) \partial_t \varepsilon  \,  \varepsilon \Phi_2
-\lambda \int  \partial_x \partial_t \varepsilon  \,  \varepsilon \partial_x\Phi_2
\\
& + \frac 12 \int \left\{ \Big[ (\partial_x \varepsilon)^2 + \varepsilon^2 - \frac 23 \left((\varepsilon+V)^3 - V^3 - 3 V^2\varepsilon \right)\Big]\partial_t \Phi_1
+  [\lambda (\partial_x \varepsilon)^2 + \varepsilon^2] \partial_t \Phi_2\right\}  \\
& - \int \partial_t V \left((\varepsilon+V)^2 - V^2 - 2V\varepsilon \right) \Phi_1
= G_1+G_2+G_3+G_4+G_5+G_6.
\end{split}\end{equation*}
Let $z_3$ be such that $(1-\lambda \partial_x^2) z_3=  \varepsilon-\partial_x^2 \varepsilon$ and $z_4$ be such that $(1-\lambda \partial_x^2) z_4=-2 V \varepsilon- \varepsilon^2$,
so that by \eqref{eq:ep}
\begin{equation}\label{eq:ee}
	\partial_t \varepsilon= \partial_x z_3 + \partial_x z_4
	- (1-\lambda \partial_x^2)^{-1} E
	- \sum_{j=1,2}\left[(\dot \mu_j-{\cal M}_j) \frac {\partial V}{\partial \mu_j}
	- (\mu_j-\dot \yy_j-{\cal N}_j) \frac {\partial V}{\partial y_j}\right].
\end{equation}
Then,
\begin{equation*}\begin{split} 
G_1&= \int \partial_x (z_3+z_4) (z_3+z_4 - \lambda \partial_x^2 (z_3+z_4)) \Phi_1  \\
&- \int (1-\lambda \partial_x^2)^{-1} E \left(-\partial_x^2 \varepsilon + \varepsilon - ((\varepsilon +V)^2-V^2)\right)\Phi_1 \\
&- \sum_{j=1,2} (\dot \mu_j-{\cal M}_j) \int  \frac {\partial V}{\partial \mu_j} \left(-\partial_x^2 \varepsilon + \varepsilon - ((\varepsilon +V)^2-V^2)\right)\Phi_1\\
&+ \sum_{j=1,2}  (\mu_j-\dot \yy_j-{\cal N}_j) \int \frac {\partial V}{\partial y_j} \left(-\partial_x^2 \varepsilon + \varepsilon - ((\varepsilon +V)^2-V^2)\right)\Phi_1
\\ &=G_{1,1}+G_{1,2}+G_{1,3}+G_{1,4}.
\end{split}\end{equation*}

\begin{equation*}\begin{split} 
G_2&= - \int \partial_x (z_3+z_4) \partial_x\varepsilon \partial_x\Phi_1  
+ \int (1-\lambda \partial_x^2)^{-1} E \partial_x\varepsilon \partial_x\Phi_1  \\
&+ \sum_{j=1,2} (\dot \mu_j-{\cal M}_j) \int  \frac {\partial V}{\partial \mu_j} \partial_x\varepsilon \partial_x\Phi_1 
- \sum_{j=1,2}  (\mu_j-\dot \yy_j-{\cal N}_j) \int \frac {\partial V}{\partial y_j} \partial_x\varepsilon \partial_x\Phi_1
\\& =G_{2,1}+G_{2,2}+G_{2,3}+G_{2,4}.
\end{split}\end{equation*}
\begin{equation*}\begin{split} 
G_3&=  \int \partial_x (-\partial_x^2 \varepsilon + \varepsilon 
- ((V+\varepsilon)^2 - V^2) )\varepsilon  \Phi_2
- \int E \varepsilon \Phi_2 
\\ & - \sum_{j=1,2} (\dot \mu_j-{\cal M}_j) \int (1-\lambda \partial_x^2) \frac {\partial V}{\partial \mu_j}\varepsilon  \Phi_2 
+ \sum_{j=1,2}  (\mu_j-\dot \yy_j-{\cal N}_j) \int (1-\lambda \partial_x^2)\frac {\partial V}{\partial y_j} \varepsilon  \Phi_2
\\& =G_{3,1}+G_{3,2}+G_{3,3}+G_{3,4}.
\end{split}\end{equation*}
\begin{equation*}\begin{split} 
G_4&= -\lambda \int \partial^2_x (z_3+z_4) \varepsilon \partial_x\Phi_2
+\lambda \int (1-\lambda \partial_x^2)^{-1} \partial_x E \varepsilon \partial_x\Phi_2  \\
&+\lambda  \sum_{j=1,2} (\dot \mu_j-{\cal M}_j) \int \partial_x \frac {\partial V}{\partial \mu_j}\varepsilon \partial_x\Phi_2
-\lambda \sum_{j=1,2}  (\mu_j-\dot \yy_j-{\cal N}_j) \int \partial_x\frac {\partial V}{\partial y_j} \varepsilon \partial_x\Phi_2
\\ &=G_{4,1}+G_{4,2}+G_{4,3}+G_{4,4}.
\end{split}\end{equation*}
Note that the terms $G_{1,2}$, $G_{2,2}$, $G_{3,2}$ and $G_{4,2}$ are readily controled by 
$C \|E\|_{L^2} \|\varepsilon\|_{L^2}.$

\medskip

Now, we focus on $G_{1,1}+G_{2,1}+G_{3,1}+G_{4,1}$. We denote by $G_7$  the quadratic parts of $G_{1,1}+G_{2,1}+G_{3,1}+G_{4,1}$, i.e. the terms coming from the linear part of the equation.
We have
\begin{equation*}\begin{split}
   G_7& =
\int \partial_x z_3 (z_3-\lambda \partial_x^2 z_3)\Phi_1 
- \int \partial_x z_3 \partial_x \varepsilon \partial_x \Phi_1
\\ & + \int \partial_x (-\partial_x^2 \varepsilon +\varepsilon) \varepsilon \Phi_2
- \lambda  \int \partial^2_x  z_3 \varepsilon \partial_x\Phi_2
\end{split}\end{equation*}

Then, using \eqref{eq:ph}, \eqref{eq:cP}
\begin{equation*}\begin{split}
& G_7  =  
\int (z_3-\lambda \partial_x^2 z_3) \partial_x z_3 \Phi_1 +
(1-\frac \lambda  2) \int  z_3 \partial_x^2 \varepsilon  \partial_x \Phi_1 \\ & + \lambda \int   \partial_x^2 z_3 \varepsilon \partial_x ( \frac 12 \Phi_1 -\Phi_2 ) + (1-\frac \lambda  2) \int z_3 \partial_x \varepsilon \partial_x^2 \Phi_1
+\frac \lambda 2 \int \partial_x z_3 \varepsilon \partial_x^2 \Phi_1 \\&
+ \int  \partial_x \left(-\partial_x^2 \varepsilon + \varepsilon \right)  \varepsilon \Phi_2  \\
&\leq  - \frac 12 \int z_3^2 \partial_x \Phi_1 + \frac \lambda 2 \int (\partial_x z_3)^2 \partial_x \Phi_1\\
& +   (1-\frac \lambda 2) \int z_3 (\varepsilon-z_3 + \lambda \partial_x^2 z_3)  \partial_x \Phi_1 
+     2 \int  \varepsilon (\varepsilon-z_3-\partial_x^2 \varepsilon)  \partial_x \Phi_2 
\\& -\frac 32 \int (\partial_x \varepsilon)^2 \partial_x \Phi_2 + \frac 12 \int \varepsilon^2 \partial_x^3 \Phi_2
-\frac 12 \int \varepsilon^2 \partial_x \Phi_2 +C\rho  \int ((\partial_x \varepsilon)^2+\varepsilon^2 + (\partial_x z_3)^2 + z_3^2)
\partial_x \Phi_1
\end{split}\end{equation*}

Integrating by parts, using \eqref{eq:ph}, \eqref{eq:cP} and then choosing $\rho$ small enough, we find
\begin{equation*}\begin{split}
& G_7 \leq 
\left(-\frac 32 +\frac \lambda 2\right) \int z_3^2 \partial_x \Phi_1 
+ \left( -\frac \lambda 2 + \frac {\lambda^2} 2\right) \int (\partial_x z_3)^2 \partial_x \Phi_1
- \frac 14 \int (\partial_x \varepsilon)^2 \partial_x \Phi_1
\\ & - \frac 34 \int \varepsilon^2 \partial_x \Phi_1 
+ \left(2- \frac \lambda 2\right) \int z_3 \varepsilon \partial_x \Phi_1
+C\rho  \int ((\partial_x \varepsilon)^2+\varepsilon^2 + (\partial_x z_3)^2 + z_3^2)
\partial_x \Phi_1\\
& \leq - \frac 1{12} \int \varepsilon^2 \partial_x\Phi_1 -\frac 18  \int (\partial_x \varepsilon)^2 \partial_x\Phi_1.
\end{split}\end{equation*}

The nonlinear terms in $G_{1,1}+G_{2,1}+G_{3,1}+G_{4,1}$ which contain $\partial_x \Phi_1$ or $\partial_x \Phi_2$ are treated by perturbation (for $\epsilon$ small and $\pyy$ large) exactly as in the previous case, using the signed rest 
$ - \frac 1{12} \int \varepsilon^2 \partial_x\Phi_1 -\frac 18  \int (\partial_x \varepsilon)^2 \partial_x\Phi_1$ obtained above. 

In $G_{3,1}$,  we are left with one cubic term $ - \int \varepsilon^2 \partial_x V \Phi_2$, 
which cannot be controlled by $\Theta$, nor by the rest term above, since $\Phi_2$ does not appear with a derivative.
Thus,  computing the main order of this term, and estimating the rest by $\Theta$, we obtain
\begin{equation}\label{eq:per}
 G_{1,1}+G_{2,1}+G_{3,1}+G_{4,1}   \leq  - \sum_{j=1,2}  \frac {\mu_j(t)}{(1+\mu_1(t))^2} \int \varepsilon^2 \partial_x \qtud  +  C \Theta.
\end{equation}

After some computations, similarly as before, we obtain
\begin{equation}\label{eq:g3478}\begin{split}
	& |G_{2,3}|+|G_{2,4}|+|G_{3,3}|+|G_{3,4}|+|G_{4,3}|+|G_{4,4}|\leq C \Theta, \\
	& G_{1,3}+G_{1,4}= \sum_{j=1,2} (\dot \mu_j-\mathcal{M}_j) \nu_j \int \Lambda \qtud \varepsilon^2
	+ \sum_{j=1,2} (\mu_j - \dot y_j -\mathcal{N}_j ) \nu_j  \int   \partial_x \qtud \varepsilon^2 + O(\Theta).
\end{split}\end{equation}

The term $G_5$ is treated exactly as the term $F_3$  so that
$|G_5|\leq C \Theta$.
Finally, using \eqref{eq:de3}, the term $G_6$  writes
\begin{equation}\label{eq:G6}\begin{split}
G_6       = - \int \partial_t V \varepsilon^2 \Phi_1  &=
	\sum_{j=1,2} \frac {\mu_j(t)}{(1+\mu_1(t))^2} \int \varepsilon^2 \partial_x \qtud
 	- \sum_{j=1,2} (\dot \mu_j-\mathcal{M}_j) \frac {1}{(1+\mu_1(t))^2} \int \Lambda \qtud\varepsilon^2
	\\ &- \sum_{j=1,2} (\mu_j - \dot y_j - \mathcal{N}_j ) \frac {1}{(1+\mu_1(t))^2}  \int   \partial_x \qtud \varepsilon^2 + O(\Theta).
\end{split}\end{equation}
In conclusion, combining \eqref{eq:per}, \eqref{eq:g3478} and \eqref{eq:G6}, we finish the proof of Proposition~\ref{PR:cFG}.

\subsection{Proof of Proposition \ref{PR:STAB}.}
By classical arguments (based on the implicit function theorem -- see Lemma \ref{PR:de}, Lemma \ref{le:huit} and \cite{MMT}), there exists $\omega_1>0$, $\bar y_0>1$ such that if
\begin{equation}\label{eq:close}
\text{$\inf_{x_1-x_2 >\bar y_0 }\|u(t)-  Q_{-\mu_0}(x-x_1)-  Q_{\mu_0}(x-x_2)\|_{H^1}\leq \omega_1   $}
\end{equation}
then  $u(t)$ can be decomposed as follows
\begin{equation}
\label{eq:d0}
	u(t,x)=\qoun(t,x)+\qode(t,x)+{\overline \varepsilon}(t,x),
\end{equation}
where
\begin{equation}
\label{eq:pm}
	\qoun(t,x)= Q_{- \mu_0}(x-y_1(t)), \quad \qode(t,x)= Q_{ \mu_0}(x-y_2(t))
\end{equation}
and $y_j(t)$ are $C^1$ functions uniquely chosen so that
\begin{equation}
\label{eq:oh}
\int {\overline \varepsilon}(t,x) \partial_x \qoud(t,x) dx = 0.
\end{equation}
Moreover, $\|\overline \varepsilon\|_{H^1}\leq C \omega_1$.

Let 
 $\pyy(t)=y_1(t)-y_2(t).$
By \eqref{eq:Qm}, the functions ${\overline \varepsilon}(t,x)$ and $y_j(t)$ satisfy the following equation
\begin{equation}\label{eq:et}\begin{split}
	& (1-\lambda \partial_x^2) \partial_t {\overline \varepsilon} 
	+ \partial_x(\partial_x^2 {\overline \varepsilon} - {\overline \varepsilon} + 2 (\qoun+\qode) {\overline \varepsilon} + {\overline \varepsilon}^2 ) = 
	-2 \partial_x (\qoun \qode) \\
&\quad + (-\mu_0 - \dot \yy_1 ) (1-\lambda\partial_x^2) \partial_x  \qoun
+ (\mu_0 - \dot \yy_2 ) (1-\lambda\partial_x^2) \partial_x    \qode,
\end{split}\end{equation}
and as in Lemma \ref{PR:de}, we obtain
\begin{equation}
\label{eq:eQp}
\left| (-1)^j  \mu_0 -  \dot y_j \right|\leq
C [ \|{\overline \varepsilon}\|_{H^1} +  e^{ -\frac 34 \pyy}],
\end{equation}
where we have used (see \eqref{eq:QMM} and \eqref{eq:QQ})
\begin{equation}\label{eq:inter}
	\int \qode(t) \qoun(t) \leq C e^{-\frac 34 \pyy(t)}.
	\end{equation}

\medskip

\noindent \emph{Proof of \eqref{eq:stab+}.}\quad 
For  $C_*>2$  to be chosen later, 
assume \eqref{eq:clo} and define 
\begin{equation*}\begin{split}
T^*  =  \sup \Big\{ & t_0<T<- (\rho\mu_0)^{-1} |\log \mu_0|
\hbox{ such that, for all $t_0<t<T$,
$u(t)$ satisfies \eqref{eq:close},} \\ & \text{ $\|{\overline \varepsilon}(t)\|_{H^1} \leq C_*  \omega\mu_0 + C_* e^{- 4 \rho {\mu_0}  |t|}$  and $\pyy(t) >  \frac 32 \mu_0 |t|$}\Big\}.
\end{split}\end{equation*}
Note that for $C_*$ large enough, $T^*$ is well-defined by \eqref{eq:clo} and by continuity of $u(t)$ in $H^1$.

We prove that  $T^*=- (\rho\mu_0)^{-1}|\log \mu_0|$, for $C_*$ large enough, assuming by contradiction 
that $T^* < - (\rho\mu_0)^{-1} |\log \mu_0|$ and working on the interval $[t_0,T^*]$.

First, we claim the following control of the scaling directions of ${\overline \varepsilon}(t)$.

\begin{claim}\label{eq:scg}
For all $t\in [t_0,T^*]$,
	\begin{equation}\label{eq:tQ}
		\left|\int {\overline \varepsilon}(t,x) (1-\lambda \partial_x^2) \qoud(t,x) dx \right|\leq 
		C [\mu_0^{-1} \sup_{[t_0,t]}  \|{\overline \varepsilon}\|_{H^1}^2  + \sup_{[t_0,t]} e^{-\frac 12 \pyy} + \mu_0 \omega].
	\end{equation}
\end{claim}
\begin{proof}[Proof of Claim \ref{eq:scg}]
 Indeed, \eqref{eq:scg} is  obtained by expanding $u(t)=\qoun(t)+\qode(t)+{\overline \varepsilon}(t)$ in the conservation laws \eqref{eq:i5} and \eqref{eq:i6} (i.e.  $M(u(t_0))=M(u(t))$ and $\mathcal{E}(u(t_0))=\mathcal{E}(u(t))$) using \eqref{eq:Qm}, \eqref{eq:inter} and  \eqref{eq:clo}:
\begin{align*}
 M(u(t_0)) & = M(\qoun(t_0)) + M(\qode(t_0))+ 2 \int \overline \varepsilon(t_0) (1-\lambda \partial_x^2 )\qoun(t_0)\\ &
+  2 \int \overline \varepsilon(t_0) (1-\lambda \partial_x^2 )\qode(t_0) +O(e^{-\frac 34 y(t_0)})+O(\|\overline \varepsilon(t_0)\|^2_{H^1})	 \\
& = M(u(t))  =  M(\qoun(t)) + M(\qode(t))+ 2 \int \overline \varepsilon(t) (1-\lambda \partial_x^2 )\qoun(t)\\ &
+  2 \int \overline \varepsilon(t) (1-\lambda \partial_x^2 )\qode(t) +O(e^{-\frac 34 y(t)})+O(\|\overline \varepsilon(t)\|^2_{H^1});
\end{align*}
\begin{align*}
 \mathcal{E}(u(t_0)) & = \mathcal{E}(\qoun(t_0)) + \mathcal{E}(\qode(t_0))+ 2 \mu_0 \int \overline \varepsilon(t_0) (1-\lambda \partial_x^2 )\qoun(t_0) \\ &
-  2\mu_0 \int \overline \varepsilon(t_0) (1-\lambda \partial_x^2 )\qode(t_0) +O(e^{-\frac 34 y(t_0)})+O(\|\overline \varepsilon(t_0)\|^2_{H^1})	 \\
& = \mathcal{E}(u(t))  =  \mathcal{E}(\qoun(t)) + \mathcal{E}(\qode(t))+ 2 \mu_0 \int \overline \varepsilon(t) (1-\lambda \partial_x^2 )\qoun(t) \\ &
-  2 \mu_0\int \overline \varepsilon(t) (1-\lambda \partial_x^2 )\qode(t) +O(e^{-\frac 34 y(t)})+O(\|\overline \varepsilon(t)\|^2_{H^1});	
\end{align*}
Using $M(\qoud(t_0))=M(\qoud(t))$,  $\mathcal{E}(\qoud(t_0))=\mathcal{E}(\qoud(t))$, and
$\|\overline \varepsilon(t_0)\|_{H^1} \leq C \mu_0 \omega$, we find \eqref{eq:tQ}.
\end{proof}

Now, we use a functional $\overline {\cal F}$ similar to ${\cal F}_-$. Let
\begin{equation}\label{eq:etG}\begin{split}
	{\overline {\cal F}}(t) & =
		\int \left[(\partial_x {\overline \varepsilon})^2 + {\overline \varepsilon}^2 - \frac 23 \left(({\overline \varepsilon}+\qoun+\qode)^3 - (\qoun+\qode)^3 - 3 (\qoun+\qode)^2{\overline \varepsilon} \right)\right]\overline\Phi_1(x)
		\\ &
		+\int \left[\lambda (\partial_x {\overline \varepsilon})^2 + {\overline \varepsilon}^2\right]\overline\Phi_2(x),
\end{split}\end{equation}
where, $\varphi$ being defined in \eqref{eq:ph},
$$\overline \Phi_1(x)=\frac {\varphi(x)}{(1-\mu_0)^2}+ \frac {1-\varphi(x)}{(1+\mu_0)^2} ,\quad 
\overline \Phi_2(x)=\frac {-\mu_0 \varphi(x)}{(1-\mu_0)^2}+ \frac {\mu_0(1-\varphi(x))}{(1+\mu_0)^2} .$$
We perform similar (and simpler) computations as the ones of Propositions \ref{PR:cFG} and \ref{PR:cFG}
(scaling parameters and $\Phi_{j}$ are time independent here). We obtain, for some $\rho>0$ small enough
\begin{equation}\label{eq:ctG}
\frac d{dt}{\overline {\cal F}} (t)    \leq C \|{\overline \varepsilon}\|_{L^2} \left( e^{- 2 \rho \pyy } \|{\overline \varepsilon}\|_{L^2} + e^{-\frac 34 \pyy}\right).
\end{equation}
From this point, the end of the proof is the same as the one of Proposition 3.2 in \cite{MMkdv4} and it is omitted.

\medskip

\noindent\emph{Proof of \eqref{eq:stab-}.} It is completely similar.

\medskip

\noindent\emph{Proof of \eqref{eq:as2}.} The asymptotic stability is a consequence of results in \cite{MM1}, \cite{Mi2}, \cite{Di}, \cite{DiMa} and \cite{Ma}.

\subsection{Proof of Proposition \ref{PR:SHARP}}
 
 For $X_1,X_2 \in \RR$, let $U_{X_1,X_2}$ be the unique  solution of {\eqref{eq:BBM}} such that
	\begin{equation}\label{eq:X1X2}
		\lim_{t\to -\infty} \|U_{X_1,X_2}(t) - Q_{-\mu_0}(x + \mu_0 t -X_1) - Q_{\mu_0}(x- \mu_0 t-X_2) \|_{H^1} =0.
	\end{equation}
	Then, for any $Y_1,Y_2\in \RR$, one has
	\begin{equation}\label{eq:Y1Y2}\begin{split}
		& U_{Y_1,Y_2}(t,x) = U_{X_1,X_2}(t-T_0,x-X_0) \\ & \text{where} \quad X_0 = \frac 12 (Y_1-X_1) + \frac 12 (Y_2-X_2),
		\ T_0 = \frac 1{2\mu_0}\left( (X_2-Y_2) -   (X_1-Y_1)\right).
	\end{split}\end{equation}
	In particular, the map $(X_1,X_2)\mapsto U_{X_1,X_2}$ is smooth and
	\begin{equation}\label{eq:map}
		\frac {\partial U_{X_1,X_2}}{\partial X_j} = (-1)^{j}\frac 1{2 \mu_0} \frac {\partial U_{X_1,X_2}}{\partial t} 
		- \frac 12 \frac {\partial U_{X_1,X_2}}{\partial x}.
	\end{equation}

We  assume ${{\mathcal T}_1}\in (-T,T)$, the case $|{\mathcal T}_1| > T$ being similar, and
we prove the stability result for $t\in (-\infty,{{\mathcal T}_1}]$, the stability proof for $t>{{\mathcal T}_1}$ following from similar arguments.

For $C_1>2$ to be chosen, we define
\begin{equation}\label{eq:sh5} \begin{split}
	T^* = \inf \{ t \leq {{\mathcal T}_1} \ ; & \ \text{such that for all $t\leq t'\leq {{\mathcal T}_1}$, } \\
 	&	\inf_{X_1,X_2} \|u(t')-U_{X_1,X_2}(t')\|_{H^1}\leq C_1 \omega \mu_0 \}.
\end{split}\end{equation}
By the assumption on $u({{\mathcal T}_1})$, $C_1>2$ and continuity of $u(t)$ in $H^1$, $T^*<{{\mathcal T}_1}$ is well-defined. We prove that $T^*=-\infty$ by using a contradiction argument~: we assume   $T^*>-\infty$ and we obtain a contradiction by strictly improving the estimate of $\inf_{X_1,X_2}\|u-U_{X_1,X_2}\|_{H^1}$ on $t\in  [T^*,{{\mathcal T}_1}]$.

By Proposition \ref{pr:st}, $U_{X_1,X_2}$ is close for all time to the sum of two distant solitons.
Thus, on $[T^*,{{\mathcal T}_1}]$, for $\omega$ small enough, we can use modulation theory (as in Lemma \ref{PR:de}) to obtain $(X_1(t), X_2(t))\in \RR^2$, such that
\begin{equation}\label{eq:modW}\begin{split}
	&	u(t,x)={\tilde U}(t,x)+ \tilde { \varepsilon}(t,x), \quad {\tilde U}(t,x)=U_{X_1(t),X_2(t)}(t,x),\\
	&	\int \tilde { \varepsilon}(t,x) (1-\lambda \partial_x^2) \frac {\partial \tilde U}{\partial X_j} =0 \quad (j=1,2). 
\end{split}\end{equation}
Moreover, $\tilde { \varepsilon}$ satisfies $\|\tilde { \varepsilon}(t)\|_{H^1} \leq C C_1 \omega \mu_0$, and
\begin{equation}\label{eq:sh8}
	(1-\lambda \partial_x^2)\partial_t \tilde { \varepsilon} + \partial_x (\partial_x^2 \tilde { \varepsilon} -\tilde { \varepsilon} + 2 {\tilde U} \tilde { \varepsilon} + \tilde { \varepsilon}^2)
	+ \sum_{j=1,2} \dot X_j (1-\lambda\partial_x^2 ) \frac {\partial \tilde U}{\partial X_j} = 0,
\end{equation}
and
\begin{equation}\label{eq:sh9}
	|\dot X_1|+|\dot X_2|\leq C \| \tilde { \varepsilon}\|_{H^1}.
\end{equation}

Note that there exists $\tilde \mu_j(t)$ and $\tilde y_j(t)$ such that for all $t$:
\begin{equation}\label{eq:sh6}
	\|\tilde U(t) - \sum_{j=1,2} Q_{\tilde \mu_j(t)}(x-\tilde y_j(t))\|_{H^1}\leq CY_0 e^{-\YYzz}.
\end{equation}
Moreover, as in  Proposition \ref{pr:st}, there exists $t_0$ such that $\tilde \mu_{1}(t)>\tilde \mu_{2}(t)$ if $t>t_0$ and
$\tilde \mu_{1}(t)<\tilde \mu_{2}(t)$ if $t<t_0$. We assume that $t_0 < T^*$.

To control $\tilde { \varepsilon}(t)$ on $[t_0,{{\mathcal T}_1}]$ (i.e. to prove that $T^*<t_0$), we use the functional
$$
	\tilde{\mathcal{F}}(t)= \int \left[(\partial_x \tilde { \varepsilon})^2 + \tilde { \varepsilon}^2 -\frac 23 ((\tilde { \varepsilon}+{\tilde U})^3 - {\tilde U}^3 - 3 {\tilde U}^2 \tilde { \varepsilon}) \right]\tilde \Phi_1
	+ \int \left[\lambda (\partial_x \tilde { \varepsilon})^2 + \tilde { \varepsilon}^2\right]\tilde \Phi_2,
$$
for $\tilde \Phi_j$ defined from $\tilde \mu_j(t)$ as in \eqref{eq:cG0}.

We follow the same computations as in the proof of Propositions \ref{PR:cFG} and \eqref{PR:STAB}, except that here there is no error term $E(t,x)$, and no scaling parameter; thus we get
\begin{claim}\label{B6} For all $t\in [\max(T^*,t_0),{{\mathcal T}_1}]$,
\begin{equation}\label{eq:cGW}
	\frac d{dt}\tilde{\mathcal{F}}(t)
	\geq - C \|\tilde { \varepsilon}(t)\|_{L^2}^2 e^{-(\frac 12 + \rho) \YYzz} \quad (C>0),
\end{equation}
\begin{equation}\label{eq:cGWW}
	\tilde{\mathcal{F}}(t)\geq \lambda \|\tilde { \varepsilon}(t)\|_{L^2}^2  +
	C \sum_{j=1,2} \left(\int  \tilde \varepsilon (t) (1-\lambda \partial_x^2) Q_{\tilde\mu_j(t)}(x-\tilde y_j(t))\right)^2\quad (\lambda>0).
\end{equation}
\begin{equation}\label{eq:cGWWW}
	\left| \int  \tilde\varepsilon (t) (1-\lambda \partial_x^2) Q_{\tilde\mu_j(t)}(x-\tilde y_j(t)) \right|\leq 
	C \omega \mu_0.
\end{equation}
\end{claim}
We omit the proof of Claim \ref{B6} since it is the same as the proof of Claim B.5 in \cite{MMkdv4}.

From \eqref{eq:cGW}, \eqref{eq:cGWW} and \eqref{eq:cGWWW}, and a continuity argument we deduce that $T^*<t_0$.
Note that the estimates on $|\dot X|$ and $|\dot T|$ come from \eqref{eq:sh9} and \eqref{eq:Y1Y2}.

Finally, to treat  the case $T^*<t_0$, i.e. to prove that  $T^*=-\infty$, one uses another functional, similar to $\mathcal{F}_+$. This completes the proof of Proposition~\ref{PR:SHARP}.

\section{Appendix}\label{reduction}
We write the transformation from equation \eqref{eq:BBMc} to  equation {\eqref{eq:BBM}}.
Note that solitons for   equation \eqref{eq:BBMc} are of the form ($c>1$)
\begin{equation}\label{eq:sol}
	\mathcal{R}_{c,x_0}(t,x)= \mathcal{Q}_c(x-ct-x_0)\quad
	\text{where}\quad \mathcal{Q}_c(x) = (c-1) Q\left(\sqrt{\frac {c-1}c} x \right).
\end{equation}

For $1<c_1<c_2$  close, let $\mathcal{U}(t,x)$ be the unique solution of \eqref{eq:BBMc} (see \cite{DiMa}) such that
\begin{equation}\label{eq:ori}
	\lim_{t\to -\infty}
	\|\mathcal{U}(t) - \mathcal{Q}_{c_1}(.-c_1t-x_1)- \mathcal{Q}_{c_2}(.-c_2t-x_2)\|_{H^1}=0.
\end{equation}
Let
\begin{align}
&  \bar c=\frac 12 (c_1+c_2), \quad \lambda	= \frac {\bar c-1}{\bar c},
\quad \mu_0=\frac {c_2-\bar c}{\bar c-1}= \frac {\bar c-c_1}{\bar c-1},\label{eq:i1}	\\
&  x'=\lambda^{1/2} \left(x- \frac t{1-\lambda}\right),\quad t'=\frac {\lambda^{3/2}}{1-\lambda} t,\quad 
U(t',x')=\frac {1-\lambda}\lambda \mathcal{U}(t,x). \label{eq:i2}
\end{align}
Then, $U(t,x)$ satisfies \eqref{eq:BBM} and it is the unique solution of \eqref{eq:BBM} such that
$$
\lim_{t\to -\infty}
	\|U(t) -  Q_{-\mu_0}(.+\mu_0 t-y_1)- Q_{\mu_0}(.-\mu_0 t-y_2)\|_{H^1}=0,
$$
where
$$
y_1= x_1 \sqrt{\lambda} ,\quad y_2= x_2 \sqrt{\lambda}.
$$
Indeed, for $c>1$, $x_0\in \RR$, a soliton  $\mathcal{R}_{c,x_0}(t,x)=\mathcal{Q}_{c}(x-ct-x_0)$ of \eqref{eq:BBMc} transforms by \eqref{eq:i2} into 
\begin{equation}\label{eq:i3}
	R_{c,y}  (t',x') = \left(\frac {1-\lambda}\lambda \right) (c-1) Q\left(\sqrt{\frac{c-1}{c}} \left(\lambda^{-\frac 12} x'+\lambda^{-\frac 32}    \left[1- c(1-\lambda)\right]t'  -x_0\right) \right).
\end{equation}
Setting
\begin{equation}\label{eq:i4}
	\mu=\mu(c)=\frac 1 \lambda (c(1-\lambda) -1)= \frac {c-\bar c}{\bar c-1},
	\quad y_0=\lambda^{\frac 12} x_0,
\end{equation}
one checks
$$
R_{c,y}(t',x')=(1+\mu) Q\left(\sqrt{\frac {1+\mu}{1+\lambda \mu}} \left(x'-\mu t' - y_0\right)\right)
= Q_{\mu} (x' - \mu t' - y_0).
$$

\end{document}